\documentclass[sn-mathphys,Numbered]{sn-jnl}

\usepackage{graphicx}%
\usepackage{multirow}%
\usepackage{amsmath,amssymb,amsfonts}%
\usepackage{amsthm}%
\usepackage{mathrsfs}%
\usepackage[title]{appendix}%
\usepackage{xcolor}%
\usepackage{textcomp}%
\usepackage{manyfoot}%
\usepackage{booktabs}%
\usepackage{algorithm}%
\usepackage{algorithmicx}%
\usepackage{algpseudocode}%
\usepackage{listings}%
\usepackage{cprotect}
\usepackage{fancyvrb}

\newtheorem{theorem}{Theorem}
\newtheorem{proposition}[theorem]{Proposition}
\theoremstyle{thmstyletwo}

\theoremstyle{thmstylethree}
\newtheorem{definition}{Definition}

\begin{document}

\title[Article Title]{Formalization of Filter-Method Construction for Real Numbers in Coq}

\author{\fnm{Guowei} \sur{Dou}}\email{dgw@bupt.edu.cn}

\author{\fnm{Wensheng} \sur{Yu}}\email{wsyu@bupt.edu.cn}

\affil{\orgdiv{Beijing Key Laboratory of Space-ground Interconnection and Convergence},
\orgname{School of Electronic Engineering}, \orgname{Beijing University of Posts and Telecommunications}, \orgaddress{\city{Beijing}, \postcode{100876}, \country{China}}}

\abstract{
This paper presents a formalization for a theory of constructing real numbers using concepts related to filters.
The formalization is implemented using the Coq proof assistant and is grounded in axiomatic set theory.
The theory encompasses various notions such as filters, ultrafilters, principal ultrafilters, and arithmetical ultrafilters.
A special filter, namely the non-principal arithmetical ultrafilter, whose existence requires the Continuum Hypothesis to be proven,
is employed to establish several non-standard number sets: $^{*}\mathbf N$, $^{*}\mathbf Z$ and $^{*}\mathbf Q$.
The set of real numbers, $\mathbf R$, is subsequently derived by equivalently classifying a specific subset of $^{*}\mathbf Q$.
The entire theory is thoroughly formalized, with each detail verified to ensure rigor and precision.
This work contributes a new selection of foundational materials for the formalization of mathematical theories.
}


\maketitle

\section{Introduction}\label{sec1}

The expansion of number systems consistently evolves through the development of mathematical theories\cite{Kli1972}.
Throughout history, the deepening understanding of the completeness of real numbers is evident in the arithmetization of analysis\cite{Lan1966,Wan2018}.
This process is closely connected to overcoming three well-known crises in the history of mathematics:
the discovery of irrational numbers, the rigorization of analysis, and the resolution of paradoxes in set theory\cite{Wen2018,Kli1972}.

Occurring in the 19th century, the development of real number theory marks a pivotal phase in the arithmetization of analysis.
Through this advancement, humanity effectively addressed the initial mathematical crisis and laid a mathematical groundwork for resolving the second crisis.

In the initial stages of formulating real number theory, methods can be broadly categorized into two approaches:
the Dedekind cut method introduced by Dedekind and the Cauchy sequence method proposed by Cantor,
both of which have become well-established in mathematical discourse.
These two enduring and classic theories continue to be widely accepted and applied in contemporary mathematical research.

Nevertheless, it appears that human curiosity concerning the real numbers is not quite quenched with just these two constructions
and there is still new methods of constructing real numbers proposed up to the 21th century\cite{Weiss2015}.
Meanwhile, fueled by applications in automated theorem proving and verification, where one must represent the real numbers in a computer,
nuances of the differences between various constructions of the reals become very pronounced\cite{Weiss2015}.
As a result, the formal verification of these theories on real numbers represents both a profound exploration of the fundamental theories of mathematical analysis
and a significant enrichment of the foundational materials required for the formalization of mathematics.
This enrichment, in turn, offers a broader array of choices for future formalization efforts.

In this paper, we present the formal implementation of a theory for constructing real numbers based on the Coq proof assistant
and grounded in Morse-Kelley axiomatic set theory\cite{Kel1955}.
The origin of this theory can be traced back to Professor Wang's 1998 paper\cite{Wan1998} at the University of Science and Technology of China,
and it was comprehensively introduced in one of his monographs\cite{Wan2018} published in 2018.
The theory employs a special filter, specifically the non-principal arithmetical ultrafilter, to facilitate the construction of real numbers.
In this work, we refer to it as the "Filter-Method Construction for Reals (FMCR)."

Compared to classical real number theories, FMCR possesses several distinctive features:

1. Infinity numbers\footnote
{The term ``infinity number" may not be a conventional mathematical term.
The use of this expression in this paper is intended to emphasize that in FMCR,
``infinity" is not just a representation of an infinite concept but rather a class of numbers constructed in practice.
The term ``infinitesimal number," which will be mentioned later in the text, shares a similar nature.}
make their entry into mathematics as the identity of arithmetical ultrafilters.
With the aid of non-principal arithmetical ultrafilters, the set of natural numbers can be extended to the non-standard natural number set ${^{*}}{\mathbf N}$,
encompassing infinity natural numbers that essentially correspond to non-principal arithmetical ultrafilters.
Additionally, following classical methods of extending natural numbers to rational numbers,
${^{*}}{\mathbf N}$ can be further extended to the non-standard integer set ${^{*}}{\mathbf Z}$ inclusive of infinity integers,
and the non-standard rational number set ${^{*}}{\mathbf Q}$, which encompasses both infinity and infinitesimal numbers.
The structure of real numbers is derived by equivalently classifying a specific subset of ${^{*}}{\mathbf Q}$.

2. The establishment of this method requires additional set theory hypotheses as premises.
The existence of non-principal arithmetical ultrafilters has not been proven within general axiomatic set theories\cite{Wan2016,Wan2018}.
In 2003, Bartoszynski and Shelah asserted that the existence of non-principal arithmetical ultrafilters cannot be proven solely within the framework of ZFC set theory\cite{BS2003}.
Although this conclusion was retracted due to some errors, it is generally believed to be correct\cite{Wan2016,Wan2018}.
Proving the existence of non-principal arithmetical ultrafilters necessitates seeking support from additional set theory hypotheses.
Among the various set theory hypotheses capable of deducing the existence of non-principal arithmetical ultrafilters,
the Continuum Hypothesis (CH) is undoubtedly well-known and regarded as a safe assumption\cite{Bel2005,Coh1966}.
From the standpoint of mathematical foundations, this holds extraordinary significance,
as the existing set theory system may be inadequate for the future development of mathematics\cite{Wan2016,Wan2018}.

3. This method could potentially be extended to establish the foundation of non-standard analysis -- non-standard real numbers, also known as hyper-reals.
The comprehension of the essence of non-standard real numbers represents a recent milestone in the investigation of number systems
during the modern mathematical era of the 20th century.
Robinson once quoted G\"odel's statements in the preface of his masterpiece {\it Nonstandard Analysis}\cite{Rob1974}:
``Arithmetic starts with the integers and proceeds by successively enlarging the number system by rational and negative numbers, irrational numbers, etc.
But the next quite natural step after the reals, namely the introduction of infinitesimals, has simply been omitted.
I think, in coming centuries it will be considered a great oddity in the history of mathematics that
the first exact theory of infinitesimals was developed 300 years after the invention of the differential calculus."
``There are good reasons to believe that non-standard analysis, in some version or other, will be the analysis of the future."
Introducing ``ideal elements", such as infinitesimals, into real numbers was precisely the original concept proposed by Leibniz
when he invented calculus three hundred years ago.
However, this idea was not widely accepted at the time due to its deviation from tradition and its seemingly illogical nature\cite{Wan2018}.
The majority of literature on non-standard analysis typically utilizes the ultrapower method to provide a model for hyper-reals.
Although FMCR does not directly construct hyper-reals, its underlying concept involving non-principal arithmetical ultrafilters
might serve as an alternative to the ultrapower method for hyper-real construction.
The author team of this paper is prepared to explore this avenue in conjunction with formal verification.

It is evident that FMCR establishes number systems from the perspective of contemporary mathematics,
inheriting and surpassing ideas found in classical methods of constructing number systems.
This approach enables us to reexamine real numbers through the lens of modern mathematics.

This paper seeks to showcase our formal implementation of FMCR, create a real number library distinct from classical construction methods,
and offer additional starting points and possibilities for future formal verification efforts in complex mathematical theories.
We have selected the Coq proof assistant as our verification tool and set theory as the formal language for description.

The structure of this paper is organized as follows:

\begin{itemize}
\item Section \ref{sec2} provides an introduction to some closely related background knowledge,
including (\ref{sec2-1}) real numbers in formalization libraries, (\ref{sec2-2}) concepts concerning non-principal arithmetical ultrafilters,
and (\ref{sec2-3}) a brief discussion about type theory and set theory.

\item Section \ref{sec3} summarizes and extends our previous formalization work on Morse-Kelley (MK) axiomatic set theory.
Subsections \ref{sec3-1} and \ref{sec3-2} present preparatory logics and the structure of MK,
while \ref{sec3-3} and \ref{sec3-4} introduce fundamental definitions and the axiom system.
Subsections \ref{sec3-5} and \ref{sec3-6} present the formalization of non-negative integers and their algebraic operations.

\item Section \ref{sec4} focuses on the formalization of notions about filters,
covering (\ref{sec4-1}) filters and ultrafilters, (\ref{sec4-2}) ultrafilter transformation, and (\ref{sec4-3}) arithmetical ultrafilter.

\item Section \ref{sec5}, as the core of this paper, details our implementation of real numbers.
Subsections include the construction of non-standard number sets ${^{*}}{\mathbf N}$ (\ref{sec5-1}),
${^{*}}{\mathbf Z}$ (\ref{sec5-2}), ${^{*}}{\mathbf Q}$ (\ref{sec5-3}), and the real number set ${\mathbf R}$ (\ref{sec5-4}).

\item Section \ref{sec6} is conclusion.
\end{itemize}

\section{Preliminary}\label{sec2}

\subsection{Real Numbers in Formalization Libraries}\label{sec2-1}

The formal verification of mathematical theorems has made significant progress in recent years with the rapid development of computer science\cite{Cas2021,Gon2008,Gon2013,GAAB2013,HABD2015},
especially with the emergence of proof assistants such as Coq\cite{BC2004}, Isabelle/HOL\cite{NPW2002}, Lean\cite{MU2021}, and others.
The formalization of real numbers, the fundamental elements for mathematical research, is undoubtedly a crucial step.
Various implementations of real numbers exist across different proof assistants\cite{BLM2016},
and as our work is built upon Coq, we will outline some of these implementations within the Coq environment.

Firstly, the Coq standard library offers comprehensive implementations for real numbers.
It initially covers axiomatic reals, followed by Dedekind cut reals and Cauchy sequence reals.
Additionally, it provides an interface for converting between Dedekind cut and Cauchy sequence reals.
The verification of the isomorphism between any two constructive real number structures is also ensured.
For more details, please refer to https://coq.inria.fr/distrib/current/stdlib/.

In addition to Coq standard library, some formalization projects offer implementations of real numbers based on Coq as well.
We will now introduce C-CoRN, MathClasses, MathComp and Gaia.

C-CoRN and MathClasses are projects from Radboud University Nijmegen.
C-CoRN, consisting of a library of constructive algebra and analysis, grew out of the FTA-project (implemented at Radboud University Nijmegen as well),
where a constructive proof of the Fundamental Theorem of Algebra was formalized in Coq\cite{CGW2004}.
Built upon a set of axioms for constructive reals utilized in the FTA-project,
the C-CoRN library includes not only a concrete model of real numbers (namely, the standard construction as Cauchy sequences of rationals)
but also a formal proof that any two real number structures are equivalent, and some alternative sets of axioms that can be used to define them\cite{CGW2004,GN2002}.

MathClasses builds a standard algebraic hierarchy and a set of mathematically sound abstract interfaces for different kinds of numbers based on
a new set of type classes for mathematics compared to first-class type classes in the Coq system\cite{SW2010,SW2011}.
MathClasses is to some degree an extension of C-CoRN since its implementation of real numbers is the result of optimizing O'Connor's\cite{KS2011,KS2013}.
Using the notion of monad (a notion in category theory),
O'Connor provided a new constructive implementation of real numbers\cite{O'Co2007} and formalized a theory of complete metric spaces\cite{O'Co2008}, extending the C-CoRN library.

MathComp (Mathematical Components) is the name of a library of formalized mathematics for the Coq system.
Using the proof language of SSReflect\cite{GM2010,GM2016},
it covers a variety of topics from the theory of basic data structures (e.g., numbers, lists, finite sets) to advanced results in various flavors of algebra\cite{MT2022}.
The MathComp project does not primarily focus on the construction of real numbers;
however, its initial version included an implementation for constructing reals from Cauchy sequences
(refer to https://github.com/math-comp/math-comp/blob/archive/mathcomp/real\_closed/cauchyreals.v).

Gaia, standing for Geometry, Algebra, Informatics and Applications,
is a project started in 2009 from INRIA (National Institute for Research in Digital Science and Technology, France).
It aims to formalize, on a computer, the fundamental notions of mathematics with the guideline of {\it Elements of Mathematics} of Bourbaki
(refer to http://www-sop.inria.fr/marelle/gaia/).
The Gaia Library is quite different from those above because its implementation relies on a theory of set rather than a theory of type\cite{Gri2010}.
It firstly formally describes a set of axioms proposed by Carlos Simpson\cite{Sim2004a,Sim2004b} that extends the logic of Coq and is sufficient to prove Bourbaki,
then implements the theory of set of Bourbaki based on Simpson's axioms, naturally including the set of natural numbers\cite{Gri2010}.
The sets of integers and rationals are derived through equivalence classification and the fraction method, respectively.
Subsequently, the set of real numbers is implemented using the Dedekind cut.\cite{Gri2016}

The foundation of our work is grounded in a set theory (resembling Gaia in this respect),
specifically the Morse-Kelley axiomatic set theory initially published in \cite{Kel1955}.
Guided by Wang's theory, we employ the filter tool to expand the set of natural numbers
into a set of non-standard natural numbers with precisely the same algebraic properties as their standard counterparts.
Building upon this, sets of non-standard integers, non-standard rational numbers, and standard real numbers are all established.
The formalization of this process will be presented in section \ref{sec4} and \ref{sec5}; Morse-Kelley set theory in section \ref{sec3}.

\subsection{Non-principal Arithmetical Ultrafilter}\label{sec2-2}

The study of filters needs to be conducted within the context of set theories.
The earliest concept related to filter is ultrafilter,
which was apparently enunciated first by Frigyes Riesz in an address\cite{Rie1909} that unfortunately did not receive at the time the attention it deserved\cite{CN1974}.
An ultrafilter is a truth-value assignment to the family of subsets of a set, and a method of convergence to infinity.
From the first property (logical) arises its connection with two-valued logic and model theory;
from the second property (convergence) arises its connection with topology and set theory\cite{CN1974}.

Ultrafilters can be classified into principal ultrafilters and non-principal ultrafilters.
The former are relatively straightforward to construct, while the construction of the latter is not as succinct or direct.
In 1929, Stanis{\l}aw Marcin Ulam proved the existence of non-principal ultrafilters over $\omega$ by transfinite induction
according to a well-order, that needs Axiom of Chioce to accomplish\cite{Ula1929}, of the power set of $\omega$.

The concept of a filter, more concise and abstract, was introduced after the ultrafilter in \cite{Car1937a,Car1937b} by Henri Paul Cartan
in 1937 and subsequently adopted by Bourbaki in their book {\it General Topology}\cite{Bou1995}.
Following this, as research deepened, an increasing number of concepts and theories related to filters were proposed.

In 1979, Daguenet-Teissier proved that every 3-arrow ultrafilter has the following property\cite{Dag1979,Wan1998}:
\begin{equation}
\forall f,g \in {\omega}^{\omega},\ f\langle F\rangle = g\langle F\rangle\ \Longleftrightarrow\ f =_F g,
\footnote{Note that our notations are different from those in \cite{Dag1979,Wan1998},
which will be explained along with the specific meaning of this propoty later in section \ref{sec4}.}
\end{equation}
\noindent where $F$ represents a 3-arrow ultrafilter.
But 3-arrow ultrafilters are at least as hard to find as P-points (which are also a type of ultrafilter),
whose existence is independent of ZFC set theory\cite{Wan1998}.

Subsequently, drawing from Daguenet-Teissier's research, Fangting Wang observed ultrafilters with property (1) and highlighted that certain ones,
specifically non-principal ultrafilters satisfying (1), could be employed to establish arithmetic models, including a real number model,
with simple structures and superior properties.
This observation was made in \cite{Wan1998} in 1998, and consequently, they were named ``arithmetical ultrafilters" in \cite{Wan2000} in 2000.

It is straightforward to prove that every principal ultrafilter is an arithmetical ultrafilter (i.e., satisfies property (1)),
but the existence of non-principal arithmetical ultrafilters is challenging to prove within the framework of general set theories.
In one of Wang's publications \cite{Wan2016},
he introduces various methods to prove the existence of non-principal ones using different set theory hypotheses independent of ZFC.
Since our implementation of reals relies on the existence of non-principal arithmetical ultrafilters,
we have conducted a formal verification of it using the Continuum Hypothesis with Wang's guidance in \cite{Wan2016}.
This ensures the consistency of our implementation of reals.

In a parallel but independent exploration from Wang's findings, Mauro Di Nasso and Marco Forti from Italy also recognized arithmetical ultrafilters
in their attempts to isolate a few topological properties that could characterize the non-standard models\cite{FNB2002,NF2005a}.
Nasso and Forti named them ``hausdorff ultrafilters"\cite{NF2006}.
They gave a simple characterization of non-standard extensions in purely topological terms
and established connections with hausdorff ultrafilters.
According to Nasso and Forti's study, every hausdorff non-standard extension\cite{FNB2002} $^{*}X$
can be identified with an invariant subspace of the Stone-\v{C}ech compactification $\beta X$ of the discrete space $X$.
In this way, the non-standard elements of $\beta X$ may be viewed as non-principal ultrafilters over $X$.\cite{NF2005a}
A concrete instance is the extension from the model of natural numbers to its non-standard model $^{*}{\mathbf N}$\cite{NF2005b},
which serves as the first step of FMCR from natural numbers to reals.

Consequently, as a matter of fact, FMCR of Wang is precisely a specific application of the hausdorff non-standard extension introduced by Nasso and Forti.
It is thoroughly presented in Wang's book \cite{Wan2018}, which serves as the guiding reference for our implementation.
The primary objective of our formalization is to establish a new real number structure and we do not emphasize the topological properties of non-standard extensions.
Instead, our focus is on the formal verification of the algebraic properties of each number structure, particularly those of the non-standard ones
(i.e., $^{*}\mathbf N$, $^{*}\mathbf Z$, $^{*}\mathbf Q$).

\subsection{Set Theory and Type Theory}\label{sec2-3}

After the discovery of Russell's Paradox in naive set theory in the late 19th century,
many mathematicians delved into an in-depth study of mathematical foundations to resolve this logical contradiction.
Among numerous theories related, the two most representative ones are axiomatic set theory and type theory.
Even Kurt G\"{o}del, when presenting his celebrated incompleteness theorems in 1931, gave the following assessment:
The most comprehensive formal systems that have been set up hitherto are the system of {\it Principia Mathematica} (referring type theory)
on the one hand and the Zermelo-Fraenkel axiom system of set theory on the other\cite{AH2014}.

The method of axiomatic set theory is to provide a set of axioms that sets should satisfy, and on this basis study the properties of sets.
In 1908, Zermelo proposed the first axiomatic system of set theory.
Through the improvement of Fraenkel and Skolem, the famous Zermelo-Fraenkel set theory (ZFC) that includes the axiom of choice was formed\cite{BF1958}.
In 1920, Von Neumann proposed his axiomatic system.
It was revised by Bernays in 1937 and further simplified by G\"{o}del in 1940\cite{Jec2003}.
This is the von Neumann-Bernays-G\"{o}del set theory (NBG).

Another important axiomatic set theory is Morse-Kelley set theory (MK),
which was first proposed by Wang\cite{Wan1949} in 1949 and was formally published in Kelley's {\it General Topology}\cite{Kel1955} in 1955.
In {\it General Topology}, Kelley clearly pointed out that
MK is a variant of systems of Skolem and of Morse and owes much to the NBG system.
Like NBG, MK admits classes (whose quantity is more extensive than that of sets) as fundamental objects.
Differently, a finite axiom system is abandoned and the development is based on eight axioms and one axiom-scheme
(i.e., all statements of a specific prescribed form are accepted as axioms)\cite{Kel1955}.
Thus MK cannot be finitely axiomatized and it is strictly stronger than both NBG and ZFC.
In fact, NBG and ZFC can be proven consistent in MK\cite{BF1958}; NBG is a conservative extension of ZFC and MK is a proper extension of NBG.
Monk, Rubin, and Mendelson submit that MK does what is expected of a set theory while being less cumbersome than ZFC and NBG\cite{Men1997,Mon1969,Rub1967}.

The type theoretic approach to foundations of mathematics is based on rather different primitives than axiomatic set theory.
In general, type theory categorizes all objects according to a predetermined hierarchical structure to break free from the Russell's Paradox.
At its core lies the simply typed lambda calculus, an abstract theory of functions and computation proposed by Church in 1940\cite{Chu1940},
which can be thought of as a simple programming language\cite{Kir2018}.

The design of lambda calculus was originally intended to use functions rather than sets to study mathematical foundations.
In this system, Church accurately formulated the syntax including quantifiers,
solving inconsistencies in the original lambda calculus, based on which many later type systems such as Martin-L\"{o}f type theory and homotopy type theory were established.
This provides various theoretic foundations for current proof assistants including Coq that our work is based on. \cite{JLWZ2020En}

Compared to set theory, type theory can not only serve as the foundation of mathematics but also
be applied in many other fields such as compiler verification, operating-system microkernel verification, and circuit design verification. \cite{JLWZ2020En}

From the perspective of the formal verification of mathematical theory, both set theory and type theory are acceptable,
and some formalization work also demonstrates compatibility between this two\cite{Bar2010,Wer1997}. Thus, we have chosen the path of set theory.
On the one hand, in recent years, our team has completed a comprehensive formalization of MK\cite{SY2020,YSF2020} and is constantly improving it;
on the other hand, filter is a notion in set theory, and as illustrated in section \ref{sec1},
FMCR requires the assistance of additional set theory hypotheses and could potentially be extended to non-standard analysis,
which is also an in-depth exploration for set theory itself.

\section{Morse-Kelley Axiomatic Set Theory}\label{sec3}

The implementation of MK presented in \cite{SY2020,YSF2020} is an early version, and since then we have streamlined and optimized the code.
Building upon these improvements, we will now introduce essential components of MK including axiom system, fundamental definitions
and theorems that will be utilized in subsequent sections.

\subsection{Preparatory Logic}\label{sec3-1}

Introducing fundamental logical symbols in the code brings it closer to mathematical language, thereby enhancing readability.
The notation ``$\forall$" and ``$\exists$" respectively represent the Universal Quantifier ``forall" and Existential Quantifier ``exists"
(``$\exists$!" stands for ``exists unique").
Lambda Calculus, implemented through ``{\tt fun}" in Coq system, is useful for defining the ``classifier" in MK and
we introduce the notation ``$\lambda$" for it.

\begin{lstlisting}
Notation "$\forall$ x .. y , P" := (forall x, .. (forall y, P) ..)
  (at level 200, x binder, y binder, right associativity,
  format "'[ ' $\forall$ x .. y ']' , P").

Notation "$\exists$ x .. y , P" := (exists x, .. (exists y, P) ..)
  (at level 200, x binder, y binder, right associativity,
  format "'[ ' $\exists$ x .. y ']' , P").

Notation "$\exists$! x .. y , P" := (exists ! x, .. (exists ! y, P) ..)
  (at level 200, x binder, y binder, right associativity,
  format "'[ ' $\exists$! x .. y ']' , P").

Notation "'$\lambda$' x .. y , t" := (fun x => .. (fun y => t) ..)
  (at level 200, x binder, y binder, right associativity,
  format "'[ ' '$\lambda$' x .. y ']' , t").
\end{lstlisting}

MK is grounded in Classical Logic, whereas the Coq system adopts Intuitionistic Logic.
Consequently, we need to acknowledge the validity of the Law of Excluded Middle,
which allows for the immediate verification of certain simple logic propositions.

\begin{lstlisting}
Axiom classic : $\forall$ (P : Prop), P \/ $\sim$P.

Proposition peirce : $\forall$ (P : Prop), ($\sim$P -> P) -> P.
Proof. intros; destruct (classic P); auto. Qed.

Proposition NNPP : $\forall$ (P : Prop), $\sim$$\sim$P <-> P.
Proof. split; intros; destruct (classic P); tauto. Qed.

Proposition notandor : $\forall$ (P Q : Prop),
  ($\sim$(P /\ Q) <-> ($\sim$P) \/ ($\sim$Q)) /\ ($\sim$(P \/ Q) <-> ($\sim$P) /\ ($\sim$Q)).
Proof. intros; destruct (classic P); tauto. Qed.

Proposition inp : $\forall$ {P Q: Prop}, (P -> Q) -> ($\sim$Q) -> ($\sim$P).
Proof. intros; intro; auto. Qed.
\end{lstlisting}

In fact, Law of Excluded Middle (i.e., ``{\tt classic}" in code)
and propositions above can be acquired by importing ``Classical" in standard library of Coq, but it is not necessary for our implementation.
In ``Classical" library, there are much more propositions and theorems about Classical Logic but we do not need to use all of them.
The above four propositions are already adequate for the formalization of MK and FMCR.

\subsection{Structure of MK}\label{sec3-2}

Compared to ZFC that only pays attention to the objects of ``set", MK acknowledges the more extensive objects of ``class".
The term ``class" does not appear in any axiom, definition or theorem in MK,
but the primary interpretation of these statements is as assertions about classes.
This is a term that is directly recognized without definition.
From the aspect of type theory, ``class" is the type of all mathematical objects (sets, functions, non-set objects, etc.).
In other words, every mathematical object is a class.

There are two primitive constants besides the term ``class".
The first is ``$\in$", which is read ``is a member of" or ``belongs to".
This means that one mathematical object (class) can be an element of another one.
The second is denoted, rather strangely, ``$\{\cdot\cdot\ :\ \cdots\}$" and is read ``the class of all $\cdot\cdot$ such that $\cdots$".
It is the classifier and represents a class consisting of classes that satisfy a specific property.
For example, ``$\{x\ :\ x\notin y\}$" represents the class that consists of all members not belonging to $y$.

The constants ``class", ``$\in$" and ``$\{\cdot\cdot\ :\ \cdots\}$" form the basic structure of set language of MK,
through which the axioms, definitions and theorems in MK could be described and proved.
Thus we use the keyword ``{\tt \color[RGB]{255,69,0} Class}" in Coq to encapsulate the three into a typeclass,
enhancing the clarity and conciseness of the formalized MK system.

\begin{lstlisting}
Class mk_structure := {
  $\mbox{\color[RGB]{0,0,0}Class}$ : Type;
  In : $\mbox{\color[RGB]{0,0,0}Class}$ -> $\mbox{\color[RGB]{0,0,0}Class}$ -> Prop;
  Classifier : ($\mbox{\color[RGB]{0,0,0}Class}$ -> Prop) -> $\mbox{\color[RGB]{0,0,0}Class}$
  }.

Parameter MKS: mk_structure.

Notation "x $\in$ y" := (@ In MKS x y) (at level 70).
Notation "x $\notin$ y" := ($\sim$ (x $\in$ y)) (at level 10).
Notation "\{ P \}" := (@ Classifier MKS P) (at level 0).
\end{lstlisting}

In the formalization, the term ``class" is denoted as ``{\tt Class}"{\footnote
{Note that it is different from and not related to the keyword ``{\tt \color[RGB]{255,69,0}Class}" in Coq system.}}.
We choose to let ``{\tt Class}" live in ``{\tt \color[RGB]{0,0,255} Type}", which is the topmost sort in Coq.
The constant ``$\in$" is denoted as ``{\tt In}" and ``$\{\cdot\cdot\ :\ \cdots\}$" as ``{\tt Classifier}".

The type of ``{\tt In}" is actually ``{\tt Class -> Class -> {\color[RGB]{0,0,255}Prop}}", signifying that it takes two classes as input and produces a proposition.
For instance, given classes $x$ and $y$, ``$x\in y$" represents the proposition that $x$ belongs to $y$.

The type of ``{\tt Classifier}" is slightly more intricate.
In its type ``{\tt (Class -> {\color[RGB]{0,0,255}Prop}) -> Class}", the first ``{\tt Class}" corresponds to the first placeholder in the classifier constant,
intended to be filled by a variable representing a member of the classifier.
The ``{\tt {\color[RGB]{0,0,255}Prop}}" corresponds to the second placeholder, meant for a proposition (whether correct or not).
The last ``{\tt Class}" indicates that the classifier is indeed a class, but we do not know if there are members in it.
Taken together, ``{\tt (Class -> {\color[RGB]{0,0,255}Prop})}" implies that the proposition in the second placeholder takes the variable from the first as a parameter.
Thus in the notation ``{\tt \verb+\{+ P \verb+\}+}", the type of ``{\tt P}", representing a proposition with a parameter, is precisely ``{\tt Class -> {\color[RGB]{0,0,255}Prop}}".

\subsection{Elementary Definitions of MK}\label{sec3-3}

The definition of set is one of the most crucial in MK.
In MK, a class that is a member of another class is termed a set. A class that is not a set is a proper class.

\begin{definition}
$x$ is a set if and only if for some $y$, $x\in y$.
\end{definition}

The definition constrains the concept of sets by excluding certain classes that are deemed ``too large"
(e.g., $\{u : u = u\}$, $\{u : u \notin u\}$) to avoid the Russell's Paradox.
The formalization in Coq is as follows:

\begin{lstlisting}
Definition Ensemble x := $\exists$ y, x $\in$ y.
\end{lstlisting}

There is no type differences among classes, sets and members in formalized system.
Every mathematical object (either set or proper class) has type ``{\tt Class}", which mostly could be defaulted in code
as the type inference mechanism in Coq can infer it.
These help the system avoid the problem of nested type mismatch in Coq formalization\cite{SY2020}.

In addition to the definition of set, MK inherits most elementary definitions and operations concerning sets
(e.g., intersection, union, complement, pairs) from naive set theory and extends them to all classes.
For example, $x \cap y$ in MK represents the intersection of classes $x$ and $y$, either of which could be a set or not.
Table \ref{tab1} lists the fundamental definitions in MK along with their notations in Coq and mathematical meanings that are frequently used in latter sections.

\begin{sidewaystable}
\renewcommand\arraystretch{0.5}
\caption{Definitions and Notations in Coq}\label{tab1}
\begin{tabular*}{\textwidth}{@{\extracolsep\fill}lcccccc}
\toprule%
Definition in Coq \footnotemark[1] & Notation in Coq & Mathematical Meaning & Mathematical Symbol \\

\midrule
\verb+Union x y := \{+ $\lambda$ \verb+z, z+ $\in$ \verb+x \/ z+ $\in$ \verb+y \}+  &
\verb+x+ $\cup$ \verb+y+             &                   union of $x$ and $y$            &     $x\ \cup\ y$ \\\specialrule{0em}{3pt}{3pt}

\verb+Intersection x y := \{+ $\lambda$ \verb+z, z+ $\in$ \verb+x /\ z+ $\in$ \verb+y \}+  &
\verb+x+ $\cap$ \verb+y+             &                   intersection of $x$ and $y$     &     $x\ \cap\ y$  \\\specialrule{0em}{3pt}{3pt}

\verb+Complement x := \{+ $\lambda$ \verb+y, y+ $\notin$ \verb+x \}+         &
$\neg$ \verb+x+                      &       complement of $x$                           &      $\neg\ x$ \\\specialrule{0em}{3pt}{3pt}

\verb+Setminus x y := x+ $\cap$ \verb+(+$\neg$ \verb+y)+  &
\verb+x+ $\sim$ \verb+y+             &     difference of $x$ and $y$ or                  & $x\ \sim\ y$ \\
&&complement of $y$ relative to $x$  &  \\\specialrule{0em}{3pt}{3pt}

$\Phi$ \verb+:= \{+ $\lambda$ \verb+x, x <> x \}+  &
--                                   & void class as well as empty set                   & $\emptyset$ \\\specialrule{0em}{3pt}{3pt}

$\mu$ \verb+:= \{+ $\lambda$ \verb+x, x = x \}+  &
--                                   & universe, the class including all sets           & $\mathcal{U}$ \\\specialrule{0em}{3pt}{3pt}

\verb+Element_I x := \{+ $\lambda$ \verb+z,+ $\forall$ \verb+y, y+ $\in$ \verb+x -> z+ $\in$ \verb+y \}+  &
$\cap$\verb+x+                       & class of the intersection of the members of $x$     & $\bigcap x$ \\\specialrule{0em}{3pt}{3pt}

\verb+Element_U x := \{+ $\lambda$ \verb+z,+ $\exists$ \verb+y, z+ $\in$ \verb+y /\ y+ $\in$ \verb+x \}+  &
$\cup$\verb+x+                       & class of the union of the members of $x$     & $\bigcup x$ \\\specialrule{0em}{3pt}{3pt}

\verb+Included x y :=+ $\forall$ \verb+z, z+ $\in$ \verb+x -> z+ $\in$ \verb+y+     &
\verb+x+ $\subset$ \verb+y+          &   $x$ is a subclass of $y$ or $x$ is contained in $y$     & $x\ \subset\ y$,\ $x\ \subseteq\ y$ \\\specialrule{0em}{3pt}{3pt}

\verb+PowerClass x := \{+ $\lambda$ \verb+y, y+ $\subset$ \verb+x \}+  &
\verb+pow(x)+                        & powerclass of $x$                            & $2^{x}$ \\\specialrule{0em}{3pt}{3pt}

\verb+Singleton x := \{+ $\lambda$ \verb+z, x+ $\in$ $\mu$ \verb+-> z = x \}+  &
\verb+[x]+                           & singleton class of $x$                       & $\{x\}$ \\\specialrule{0em}{3pt}{3pt}

\verb+Unordered x y := [x]+ $\cup$ \verb+[y]+  &
\verb+[x|y]+                         & unordered pair of $x$ and $y$                & $\{xy\}$ \\\specialrule{0em}{3pt}{3pt}

\verb+Ordered x y := [[x]|[x|y]]+  &
\verb+[x,y]+                         & ordered pair of $x$ and $y$                  & $(x,y)$ \\\specialrule{0em}{3pt}{3pt}

\verb+First z :=+ $\cap\cap$\verb+z+  &
--                                   & the first coordinate of $z$                  & $1^{st}$coord $z$ \\\specialrule{0em}{3pt}{3pt}

\verb+Second z := (+$\cap\cup$\verb+z)+$\cup$\verb+(+$\cup\cup$\verb+z)+ $\sim$ \verb+(+$\cup\cap$\verb+z)+  &
--                                   & the second coordinate of $z$                  & $2^{st}$coord $z$ \\\specialrule{0em}{3pt}{3pt}

\verb+Relation r :=+ $\forall$ \verb+z, z+ $\in$ \verb+r ->+ $\exists$ \verb+x y, z = [x,y]+  &
--                                   & $r$ is a relation iff its members are ordered pairs   & -- \\\specialrule{0em}{3pt}{3pt}

\verb+Composition r s :=+ &    \verb+r+ $\circ$ \verb+s+   &   composition of $r$ and $s$    &     $r\ \circ\ s$   \\
\verb+  \{\+ $\lambda$ \verb+x z,+ $\exists$ \verb+y, [x,y]+ $\in$ \verb+s /\ [y,z]+ $\in$ \verb+r \}\+ \footnotemark[2]  &&&  \\\specialrule{0em}{3pt}{3pt}

\verb+Inverse r := \{\+ $\lambda$ \verb+x y, [y,x]+ $\in$ \verb+r \}\+  &
\verb+r+$^{-1}$                      & relation inverse to $r$                    & $r^{-1}$ \\\specialrule{0em}{3pt}{3pt}

\verb+Function f := Relation f+      &    --               &    $f$ is a function     &     --   \\
\verb+  /\ (+$\forall$ \verb+x y z, [x,y]+ $\in$ \verb+f -> [x,z]+ $\in$ \verb+f -> y = z)+  &&&  \\\specialrule{0em}{3pt}{3pt}

\verb+Domain f := \{+ $\lambda$ \verb+x,+ $\exists$ \verb+y, [x,y]+ $\in$ \verb+f \}+  &
\verb+dom(f)+                        & domain of the class $f$ \footnotemark[3]        & domain $f$ \\\specialrule{0em}{3pt}{3pt}

\verb+Range f := \{+ $\lambda$ \verb+y,+ $\exists$ \verb+x, [x,y]+ $\in$ \verb+f \}+  &
\verb+ran(f)+                        & range of the class $f$ \footnotemark[4]        & range $f$ \\\specialrule{0em}{3pt}{3pt}

\verb+Value f x :=+ $\cap$\verb+(\{+ $\lambda$ \verb+y, [x,y]+ $\in$ \verb+f \})+  &
\verb+f[x]+                          & value of $f$ at $x$ or image of $x$ under $f$     & $f(x)$ \\\specialrule{0em}{3pt}{3pt}

\verb+Function1_1 f := Function f /\ Function (f+$^{-1}$\verb+)+  &
--                                   & $f$ is a 1-1 function                          & -- \\\specialrule{0em}{3pt}{3pt}

\verb+Exponent y x :=+               &         --                 &       class consisting of functions              &    $y^{x}$    \\
\verb+  \{+ $\lambda$ \verb+f, Function f /\ dom(f) = x /\ ran(f)+ $\subset$ \verb+y \}+  &   & whose domain is $x$ and range is contained in $y$ &  \\\specialrule{0em}{3pt}{3pt}

\verb+Cartesian x y := \{\+ $\lambda$ \verb+u v, u+ $\in$ \verb+x /\ v+ $\in$ \verb+y \}\+  &
\verb+x+ $\times$ \verb+y+           & cartesian product of $x$ and $y$               & $x\ \times\ y$ \\\specialrule{0em}{3pt}{3pt}

\verb+Restriction f x := f+ $\cap$ \verb+(x+ $\times\ \mu$\verb+)+  &
\verb+f|(x)+                         & restriction of $f$ to $x$                      & $f|x$ \\\specialrule{0em}{3pt}{3pt}

\verb+ChoiceFunction c := Function c+      &    --               &    $f$ is a choice function     &     --   \\
\verb+  /\ (+$\forall$ \verb+x, x+ $\in$ \verb+dom(c) -> c[x]+ $\in$ \verb+x)+  &&&  \\\specialrule{0em}{3pt}{3pt}

\verb+PlusOne x := x+ $\cup$ \verb+[x]+   &
--                                   &  successor of $x$                              & $x'$   \\\specialrule{0em}{3pt}{3pt}

\verb+Equivalent x y :=+      &    \verb+x+ $\approx$ \verb+y+             &    there exists a 1-1 function between $x$ and $y$     &     $x \approx y$   \\
\verb+  +$\exists$ \verb+f, Function1_1 f /\ dom(f) = x /\ ran(f) = y+  &&&  \\\specialrule{0em}{3pt}{3pt}

\botrule
\end{tabular*}

\cprotect
\footnotetext[1]{All definitions are defined with the keyword ``{\color[RGB]{255,69,0}\verb+Definition+}".}
\cprotect
\footnotetext[2]{The notation ``\verb+\{\+ $\cdots$ \verb+\}\+" is denoted for the classifiers whose members are ordered pairs.\\
{\color[RGB]{147,65,200}\verb+  Notation+}\verb+ "\{\ P \}\" := (\{+ $\lambda$ \verb+z,+ $\exists$ \verb+x y, z = [x,y] /\ P x y \}) (at level 0).+\\
\verb+ + where the type of {\tt P} is ``\verb+Class -> Class -> +{\color[RGB]{0,0,255}\verb+Prop+}". }
\footnotetext[3,4]{In MK the definition of domain and range can be specific to any classes, not limited to functions.}

\end{sidewaystable}

\subsection{Axioms of MK}\label{sec3-4}

There are nine axioms\footnote{Logically, classification axiom-scheme ({\tt AxiomII}) actually represents infinite axioms as introduced in \ref{sec2-3}.
More specifically, every time the classification axiom-scheme is applied, it will have different meanings due to the difference in {\tt P}.
Here in formalization we consider it as a single axiom.} in MK.
Similar to the constants ``class", ``$\in$" and ``$\{\cdot\cdot\ :\ \cdots\}$",
we package the nine axioms into a typeclass to clearly structuralize the formalized MK system.

\begin{lstlisting}
Class MK_Axioms := {
  A_I : $\forall$ x y, x = y <-> ($\forall$ z, z $\in$ x <-> z $\in$ y);
  A_II : $\forall$ b P, b $\in$ \{ P \} <-> Ensemble b /\ (P b);
  A_III : $\forall$ {x}, Ensemble x -> $\exists$ y, Ensemble y /\ ($\forall$ z, z $\subset$ x -> z $\in$ y);
  A_IV : $\forall$ {x y}, Ensemble x -> Ensemble y -> Ensemble (x $\cup$ y);
  A_V : $\forall$ {f}, Function f -> Ensemble dom(f) -> Ensemble ran(f);
  A_VI : $\forall$ x, Ensemble x -> Ensemble ($\cup$ x);
  A_VII : $\forall$ x, x <> $\Phi$ -> $\exists$ y, y $\in$ x /\ x $\cap$ y = $\Phi$;
  A_VIII : $\exists$ y, Ensemble y /\ $\Phi$ $\in$ y /\ ($\forall$ x, x $\in$ y -> (x $\cup$ [x]) $\in$ y);
  A_IX : $\exists$ c, ChoiceFunction c /\ dom(c) = $\mu$ $\sim$ [$\Phi$]
  }.

Parameter MK_Axiom : MK_Axioms.

Notation AxiomI := (@ A_I MK_Axiom).
Notation AxiomII := (@ A_II MK_Axiom).
Notation AxiomIII := (@ A_III MK_Axiom).
Notation AxiomIV := (@ A_IV MK_Axiom).
Notation AxiomV := (@ A_V MK_Axiom).
Notation AxiomVI := (@ A_VI MK_Axiom).
Notation AxiomVII := (@ A_VII MK_Axiom).
Notation AxiomVIII := (@ A_VIII MK_Axiom).
Notation AxiomIX := (@ A_IX MK_Axiom).
\end{lstlisting}

These axioms (from I to IX) are: axiom of extent, classification axiom-scheme,
axiom of subsets, axiom of amalgamation, axiom of regularity, axiom of infinity, axiom of choice.

{\bf Axiom of Extent} For each $x$ and each $y$ it is true that $x = y$ if and only if for each $z$, $z \in x$ when and only when $z \in y$.

The axiom of extent extends the logical meaning of equality to the relation between classes and their members,
which is frequently used to prove the equality of two classes.

{\bf Classification Axiom-Scheme} For each $x$, $x \in \{u : P\ u\}$ if and only if $x$ is a set and $P\ x$.

In naive set theory, every object satisfying the proposition $P$ could be a member of $\{u : P\ u\}$,
which will lead to some contradictory sets such as $\{u : u \notin u\}$.
You will find $\{u : u \notin u\} \in \{u : u \notin u\}$ if and only if $\{u : u \notin u\} \notin \{u : u \notin u\}$.
This is the famous Russell's Paradox.
But in MK, with the classification axiom-scheme, $\{u : u \notin u\}$ is a proper class rather than set.
Consequently, only $\{u : u \notin u\} \notin \{u : u \notin u\}$ holds and this is consistent and permitted.

{\bf Axiom of Subsets} If $x$ is a set there is a set $y$ such that for each $z$, if $z \subset x$, then $z \in y$.

{\bf Axiom of Union} If $x$ is a set and $y$ is a set so is $x \cup y$.

These two axioms enable us to construct some new sets from certain sets $x$ and $y$, such as $\{x\}$, $2^{x}$, $x \cup y$, etc.

{\bf Axiom of Substitution} If $f$ is a function and domain $f$ is a set, then range $f$ is a set.

{\bf Axiom of Amalgamation} If $x$ is a set so is $\bigcup x$.

These two axioms identify more additional sets such as $\{x\} \times y$, $x \times y$, etc, along with functions that are sets.
Actually, with the axioms of subsets, union, substitution and amalgamation,
a substantial number of sets can be derived from a limited set of predetermined ones.

{\bf Axiom of Regularity} If $x \neq \emptyset$ there is a member $y$ of $x$ such that $x \cap y = \emptyset$.

It is a priori possible that there are classes $x$ and $y$ such that $x$ is the only member of $y$ and $y$ is the only member of $x$.
More generally, it is possible that there is a class $z$ whose members exist by taking in each other's laundry,
in the sense that every member of $z$ consists of members of $z$.
The axiom of regularity denies this possibility by requiring that each non-void class $z$ has at least one member whose elements do not belong to $z$.
Thus the proposition $x \notin x$ holds, and there will not be a situation like $x \in y \wedge y \in x$.

Note that, up to this point, every concrete object is still a class rather than set.
Even the most simple class $\emptyset$ cannot be proven a set.
The next axiom will ensure the existence of sets.

{\bf Axiom of Infinity} There exists a $y$, $y$ is a set, $\emptyset \in y$ and $x \cup \{x\} \in y$ whenever $x \in y$.

The axiom of infinity asserts the unconditional existence of two sets: the infinite inductive set $y$ and one of its member void class $\emptyset$
($\emptyset$ is a set since $\emptyset \in y$).
With these two sets, all other sets characterized by previous axioms can be constructed.
Especially, the non-negative integers of MK that is to be used in FMCR can be constructed on the basis of the axiom of infinity.

{\bf Axiom of Choice} There is a choice function $c$ whose domain is $\mathcal{U} \sim \{\emptyset\}$.

The axiom of choice is an axiom regarding the existence of mappings between two sets, first proposed by Zermelo in 1904 and used to prove the Well-ordering Theorem\cite{Zer2010}.
It plays a crucial role in modern mathematics and is closely related to many profound mathematical conclusions.
In MK, the axiom of choice is used to prove the Well-ordering Theorem, Hausdorff Maximal Principle and some theorems involving cardinal numbers.

\subsection{Non-Negative Integers}\label{sec3-5}

The definition of non-negative integers necessitates the concepts of ordinals and well-ordering,
which have been formalized and can be found in \cite{SY2020,YSF2020}.

\begin{definition}
$x$ is an integer(non-negative integer) if and only if $x$ is an ordinal and $E^{-1}$ well-orders $x$.
\end{definition}

$E$ is the relation $\{(u,v) : u \in v\}$ that rules the order of non-negative integers and
the successor of a integer $x$ is defined as $x \cup \{x\}$.
All non-negative integers constitute the class $\omega$.

\begin{lstlisting}
Definition E := \{\ $\lambda$ x y, x $\in$ y \}\.
Definition Integer x := Ordinal x /\ WellOrdered (E$^{-1}$) x.
Definition PlusOne x := x $\cup$ [x].
Definition $\omega$ := \{ $\lambda$ x, Integer x \}.
\end{lstlisting}

In MK, it has been proven that $\omega$ is a set.
Although not explicitly stated, $\omega$ is actually a constructive set of natural numbers
as it satisfies the Peano's Axioms, which are proven as theorems in MK.

\begin{lstlisting}
Theorem MKT134 : $\forall$ {x}, x $\in$ $\omega$ -> (PlusOne x) $\in$ $\omega$.
Theorem MKT135a : $\Phi$ $\in$ $\omega$.
Theorem MKT135b : $\forall$ x, x $\in$ $\omega$ -> $\Phi$ <> PlusOne x.
Theorem MKT136 : $\forall$ x y, x $\in$ $\omega$ -> y $\in$ $\omega$ -> PlusOne x = PlusOne y -> x = y.
Theorem MKT137 : $\forall$ x, x $\subset$ $\omega$ -> $\Phi$ $\in$ x -> ($\forall$ u, u $\in$ x -> (PlusOne u) $\in$ x)
  -> x = $\omega$.
\end{lstlisting}

With these now we can describe the structure of $\omega$\ :
$$\emptyset\ \in\ \{\emptyset\}\ \in\ \big\{\emptyset,\ \ \{\emptyset\}\big\}\ \in\ \Big\{\emptyset,\ \ \{\emptyset\},\ \ \big\{\emptyset, \{\emptyset\}\big\}\Big\}\ \cdots$$
This is consistent with our intuitive comprehension of natural numbers:
$$0\ <\ 1\ <\ 2\ < 3\ \cdots$$
The order ``less($<$)" is defined as ``in($\in$)", and $0$ as ``$\emptyset$" who has no members,
$1$ as the successor of $0$ ``$\emptyset \cup \{\emptyset\}$"($= \{\emptyset\}$) who has only one member and so forth.

\begin{lstlisting}
Notation "0" := $\Phi$ : $\omega$_scope.
Notation "1" := [$\Phi$] : $\omega$_scope.
Notation "2" := ([$\Phi$] $\cup$ [[$\Phi$]]) : $\omega$_scope.

Fact in_$\omega$_0 : 0 $\in$ $\omega$.
Fact in_$\omega$_1 : 1 $\in$ $\omega$.
Fact in_$\omega$_2 : 2 $\in$ $\omega$.
\end{lstlisting}

Thus a natural numbers system with a set structure is built, after which our attention turns to operations on $\omega$.

\subsection{Algebraic Operations on $\omega$}\label{sec3-6}

The operations among natural numbers may seem familiar to us, and we often do not pay much attention to how they are defined.
But in formal verification, every detail matters.
Operations among natural numbers are usually defined recursively.
For example, The sum of two natural numbers $m$ and $n$ can be defined as follows:

\begin{equation}
\left\{
\begin{array}{l}
m + 0 = m \\
m + n' = (m = n)'.
\end{array}
\right.
\end{equation}

Thus $1 + 2$ can be written as $1 + 1' = (1 + 1)' = ((1 + 0)')' = (1')'$.
In Coq system, the keyword {\tt Fixpoint} provides a convenient recursion mechanism for type theory,
while we have to establish a set-based recursion method for MK, and the Recursion Theorem can accomplish this task.

\begin{theorem}[Recursion Theorem]
Let $A$ be a set, $h$ be a member of $A^{A}$, namely a function whose domain is $A$ and range is contained in $A$, and $x$ be a member of $A$,
then there exists a unique function $f \in A^{\omega}$ that satisfies:
\begin{equation}\label{equa3}
\left\{
\begin{array}{l}
f(0) = x \\
\forall n \in \omega, f(n') = h(f(n)).
\end{array}
\right.
\end{equation}
\end{theorem}

Recursive condition (\ref{equa3}) uniquely determines a function $f \in A^{\omega}$.
Thus the value of $f$ at $n'$ is dependent on the value at $n$, which aligns with the feature of recursion.
The formalization of this theorem is as follows:

\begin{lstlisting}
Theorem Recursion : $\forall$ x A h, Ensemble A
  -> Function h -> dom(h) = A -> ran(h) $\subset$ A -> x $\in$ A
  -> ($\exists$! f, Function f /\ dom(f) = $\omega$ /\ ran(f) $\subset$ A
    /\ f[0] = x /\ ($\forall$ n, n $\in$ $\omega$ -> f[PlusOne n] = h[(f[n])])).
\end{lstlisting}

In Recursion Theorem, let $A$ be $\omega$, $x$ be a member $m$ in $\omega$, and $h$ be the function as:
$$h = \{(u,v) : u \in \omega \wedge v = u'\}\ (\forall n \in \omega, h(n) = n'),$$
then there uniquely exists a function $f$ in ${\omega}^{\omega}$ that satisfies:

\begin{equation}\label{equa4}
\left\{
\begin{array}{l}
f(0) = m \\
\forall n \in \omega, f(n') = h(f(n)) = (f(n))',
\end{array}
\right.
\end{equation}

\noindent where the unique function $f$ is connected to the determined natural number $m$
and every natural number can uniquely correspond such a function. Thus here comes the next theorem:

\begin{theorem}[Addition on $\omega$]
For every $m \in \omega$, there exists a unique function $f_m \in {\omega}^{\omega}$ that satisfies:
\begin{equation}\label{equa5}
\left\{
\begin{array}{l}
f_m(0) = m \\
\forall n \in \omega, f_m(n') = (f_m(n))',
\end{array}
\right.
\end{equation}
\noindent where $f_m$ can be named ``addition function of $m$" and $f_m(n)$ be denoted as ``$m + n$",
then equation (\ref{equa5}) can be written as:
\begin{equation}\label{equa6}
\left\{
\begin{array}{l}
m + 0 = m \\
\forall n \in \omega, m + n' = (m + n)'.
\end{array}
\right.
\end{equation}
\end{theorem}

We present here the entire process of the formalization of addition on $\omega$:

\begin{lstlisting}
Theorem Pre_PlusFunction : $\forall$ m, m $\in$ $\omega$
  -> $\exists$! f, Function f /\ dom(f) = $\omega$ /\ ran(f) $\subset$ $\omega$
    /\ f[0] = m /\ ($\forall$ n, n $\in$ $\omega$ -> f[PlusOne n] = PlusOne (f[n])).

Corollary PlusFunction_uni : $\forall$ m, m $\in$ $\omega$
  -> $\exists$ f, Ensemble f /\ (\{ $\lambda$ f, Function f
    /\ dom(f) = $\omega$ /\ ran(f) $\subset$ $\omega$ /\ m $\in$ $\omega$ /\ f[0] = m
    /\ ($\forall$ n, n $\in$ $\omega$ -> f[PlusOne n] = PlusOne f[n]) \}) = [f].

Definition PlusFunction m := $\cap$(\{ $\lambda$ f, Function f
  /\ dom(f) = $\omega$ /\ ran(f) $\subset$ $\omega$ /\ m $\in$ $\omega$ /\ f[0] = m
  /\ ($\forall$ n, n $\in$ $\omega$ -> f[PlusOne n] = PlusOne f[n]) \}).

Definition Plus m n := (PlusFunction m)[n].
Notation "m + n" := (Plus m n) : $\omega$_scope.
\end{lstlisting}

In our implementation, ``{\tt Pre\_PlusFunction}" and ``{\tt PlusFunction\_uni}" guarantee the uniqueness of addition function,
which makes it feasible to define the addition function of $m$ as ``$\bigcap\{f : f\ is\ an\ addition\ function\ of\ m\}$"
(the intersection of all addition functions, actually only one, of $m$).
Without using axiom of choice, this perfectly solves the problem of addition functions not being directly constructed.
The formal definition ``{\tt PlusFunction m}" exactly refers to the addition function of $m$.

The idea of defining multiplication is consistent with that of addition.
In Recursion Theorem, let $A$ be $\omega$, $x$ be a member $m$ in $\omega$, and $h$ be the function as:
$$h = \{(u,v) : u \in \omega \wedge v = u + m\}\ (\forall n \in \omega, h(n) = n + m),$$
then there uniquely exists a function $f$ in ${\omega}^{\omega}$ that satisfies:

\begin{equation}\label{equa7}
\left\{
\begin{array}{l}
f(0) = 0 \\
\forall n \in \omega, f(n') = h(f(n)) = f(n) + m,
\end{array}
\right.
\end{equation}

\noindent where the unique $f$ is connected to the determined natural number $m$
and every natural number can uniquely correspond such a function. Thus here comes:

\begin{theorem}[Multiplication on $\omega$]
For every $m \in \omega$, there exists a unique function $f_m \in {\omega}^{\omega}$ that satisfies:
\begin{equation}\label{equa8}
\left\{
\begin{array}{l}
f_m(0) = 0 \\
\forall n \in \omega, f_m(n') = f_m(n) + m,
\end{array}
\right.
\end{equation}
\noindent where $f_m$ can be named ``multiplication function of $m$" and $f_m(n)$ be denoted as ``$m \cdot n$",
then equation (\ref{equa8}) can be written as:
\begin{equation}\label{equa9}
\left\{
\begin{array}{l}
m \cdot 0 = 0 \\
\forall n \in \omega, m \cdot n' = (m \cdot n) + m.
\end{array}
\right.
\end{equation}
\end{theorem}

The process of the formalization of multiplication is exactly consistent with that of addition:

\begin{lstlisting}
Theorem Pre_MultiFunction : $\forall$ m, m $\in$ $\omega$
  -> ($\exists$! f, Function f /\ dom(f) = $\omega$ /\ ran(f) $\subset$ $\omega$
    /\ f[0] = 0 /\ ($\forall$ n, n $\in$ $\omega$ -> f[PlusOne n] = f[n] + m)).

Corollary MultiFunction_uni : $\forall$ m, m $\in$ $\omega$
  -> $\exists$ f, Ensemble f /\ \{ $\lambda$ f, Function f
    /\ dom(f) = $\omega$ /\ ran(f) $\subset$ $\omega$ /\ m $\in$ $\omega$ /\ f[0] = 0
    /\ ($\forall$ n, n $\in$ $\omega$ -> f[PlusOne n] = f[n] + m) \} = [f].

Definition MultiFunction m := $\cap$(\{ $\lambda$ f, Function f
  /\ dom(f) = $\omega$ /\ ran(f) $\subset$ $\omega$ /\ m $\in$ $\omega$ /\ f[0] = 0
  /\ ($\forall$ n, n $\in$ $\omega$ -> f[PlusOne n] = f[n] + m) \}).

Definition Mult m n := (MultiFunction m)[n].
Notation "m $\cdot$ n" := (Mult m n)(at level 10) : $\omega$_scope.
\end{lstlisting}

Regarding the properties of operations, such as the commutative law, associative law, etc., we typically accept them directly.
However, in formal implementation, each of these properties needs to be explicitly described and verified.
The artificial description of these properties is simple and familiar to people, and we will provide some formalization of them here:

\begin{lstlisting}
Fact $\omega$_Plus_in_$\omega$ : $\forall$ m n, m $\in$ $\omega$ -> n $\in$ $\omega$ -> (m + n) $\in$ $\omega$.

Fact Plus_Property1_a : $\forall$ m, m $\in$ $\omega$ -> m + 0 = m.
Fact Plus_Property2_a : $\forall$ m n, m $\in$ $\omega$ -> n $\in$ $\omega$
  -> m + (PlusOne n) = PlusOne (m + n).
Fact Plus_Property1_b : $\forall$ m, m $\in$ $\omega$ -> 0 + m = m.
Fact Plus_Property2_b : $\forall$ m n, m $\in$ $\omega$ -> n $\in$ $\omega$
  -> (PlusOne m) + n = PlusOne (m + n).

Fact Plus_Commutation : $\forall$ m n, m $\in$ $\omega$ -> n $\in$ $\omega$ -> m + n = n + m.
Fact Plus_Association : $\forall$ m n k, m $\in$ $\omega$ -> n $\in$ $\omega$ -> k $\in$ $\omega$
  -> (m + n) + k = m + (n + k).
Fact Plus_Cancellation : $\forall$ m n k, m $\in$ $\omega$ -> n $\in$ $\omega$ -> k $\in$ $\omega$
  -> m + n = m + k <-> n = k.
Fact Plus_PrOrder : $\forall$ m n k, m $\in$ $\omega$ -> n $\in$ $\omega$ -> k $\in$ $\omega$
  -> (m + n) $\in$ (m + k) <-> n $\in$ k.

Fact $\omega$_Mult_in_$\omega$ : $\forall$ m n, m $\in$ $\omega$ -> n $\in$ $\omega$ -> (m $\cdot$ n) $\in$ $\omega$.

Fact Mult_Property1_a : $\forall$ m, m $\in$ $\omega$ -> m $\cdot$ 0 = 0.
Fact Mult_Property2_a : $\forall$ m n, m $\in$ $\omega$ -> n $\in$ $\omega$
  -> m $\cdot$ (PlusOne n) = (m $\cdot$ n) + m.
Fact Mult_Property1_b : $\forall$ m, m $\in$ $\omega$ -> 0 $\cdot$ m = 0.
Fact Mult_Property2_b : $\forall$ m n, m $\in$ $\omega$ -> n $\in$ $\omega$
  -> (PlusOne m) $\cdot$ n = (m $\cdot$ n) + n.

Fact Mult_Commutation : $\forall$ m n, m $\in$ $\omega$ -> n $\in$ $\omega$ -> m $\cdot$ n = n $\cdot$ m.
Fact Mult_Distribution : $\forall$ m n k, m $\in$ $\omega$ -> n $\in$ $\omega$ -> k $\in$ $\omega$
  -> m $\cdot$ (n + k) = (m $\cdot$ n) + (m $\cdot$ k).
Fact Mult_Association : $\forall$ m n k, m $\in$ $\omega$ -> n $\in$ $\omega$ -> k $\in$ $\omega$
  -> (m $\cdot$ n) $\cdot$ k = m $\cdot$ (n $\cdot$ k).
Fact Mult_Cancellation : $\forall$ m n k, m $\in$ $\omega$ -> n $\in$ $\omega$ -> k $\in$ $\omega$
  -> m <> 0 -> m $\cdot$ n = m $\cdot$ k -> n = k.
Fact Mult_PrOrder : $\forall$ m n k, m $\in$ $\omega$ -> n $\in$ $\omega$ -> k $\in$ $\omega$
  -> m <> 0 -> (m $\cdot$ n) $\in$ (m $\cdot$ k) <-> n $\in$ k.
\end{lstlisting}

As a straightforward application of the operations above,
concepts of even and odd numbers are formalized in our implementation.
We give here the formal description of these two concepts.

\begin{lstlisting}
Definition Even m := $\exists$ k, k $\in$ $\omega$ /\ m = 2 $\cdot$ k.
Definition $\omega$_E := \{ $\lambda$ u, Even u \}.
Corollary $\omega$_E_is_Set : Ensemble ($\omega$_E).

Definition Odd m := $\exists$ k, k $\in$ $\omega$ /\ m = (2 $\cdot$ k) + 1.
Definition $\omega$_O := \{ $\lambda$ u, Odd u \}.
Corollary $\omega$_O_is_Set : Ensemble ($\omega$_O).
\end{lstlisting}

\noindent It is trivial to verify that the classes of even and odd numbers defined with the classifier $\{\cdot\cdot : \cdots\}$ are both sets,
which reflects the most significant feature of MK that every concept starts from a class.
Besides, it is also verified in our implementation that $\{u : u\ is\ odd\} \cup \{u : u\ is\ even\} = \omega$, $\{u : u\ is\ odd\} \cap \{u : u\ is\ even\} = \emptyset$
and $\{u : u\ is\ odd\} \approx \{u : u\ is\ even\} \approx \omega$.

\section{Filters and Arithmetical Ultrafilters}\label{sec4}

Starting from this section, most objects we discuss can be proven sets.
If there are no specific instructions we will directly use the term ``set" rather than ``class" to describe them.

\subsection{Filters and Ultrafilters}\label{sec4-1}

Filters actually can be defined over any sets.
As introduced in introduction and \ref{sec2-2}, the first step of FMCR is extending $\omega$ to ${^{*}{\mathbf N}}$ with the use of non-principal arithmetical ultrafilters,
thus we mainly pay attention to the filters over $\omega$.

A filter is a family of subsets of $\omega$, namely a subset of $2^{\omega}$.
From the perspective of the properties of natural numbers, a filter over $\omega$ is essentially a compatible combination of certain properties of natural numbers.
In other words, a filter actually combines some subsets of $2^{\omega}$ through some compatibility, creating certain mutual relationships between these originally independent subsets.
The specific definition of a filter is as follows:

\begin{definition}[Filter]
Assume that $F$ is a family of subsets of $\omega$ (i.e., $F \subset 2^{\omega}$) and satisfies:\\
1) $\emptyset \notin F,\ \omega \in F,$\\
2) if $a,b \in F$, then $a \cap b \in F$,\\
3) if $a \subset b \subset \omega$ and $a \in F$, then $b \in F.$\\
$F$ is called a filter over $\omega$.
\end{definition}

The notion of ultrafilters, proposed by Frigyes Riesz\cite{Rie1909},
actually predates the concept of filters and includes one additional condition.

\begin{definition}[Ultrafilter]
A filter $F$ over $\omega$ is an ultrafilter if it satisfies:
$$\forall a,\ a \subset \omega\ \Longrightarrow\ a \in F\ \vee\ (\omega \sim a) \in F.$$
\end{definition}

According to the definition of filter and ultrafilter, we can immediately construct a simple filter $\{\omega\}$.
It is a filter but not an ultrafilter, called the ``trivial filter". The formalization is as follows:

\begin{lstlisting}
Definition Filter_On_$\omega$ F := F $\subset$ pow($\omega$) /\ $\Phi$ $\notin$ F /\ $\omega$ $\in$ F
  /\ ($\forall$ a b, a $\in$ F -> b $\in$ F -> (a $\cap$ b) $\in$ F)
  /\ ($\forall$ a b, a $\subset$ b -> b $\subset$ $\omega$ -> a $\in$ F -> b $\in$ F).

Definition ultraFilter_On_$\omega$ F := Filter_On_$\omega$ F
  /\ ($\forall$ a, a $\subset$ $\omega$ -> a $\in$ F \/ ($\omega$ $\sim$ a) $\in$ F).

Example Example1 : Filter_On_$\omega$ ([$\omega$]) /\ $\sim$ ultraFilter_On_$\omega$ ([$\omega$]).
\end{lstlisting}

Every ultrafilter possesses maximality, corresponding to a maximally compatible combination of certain properties of natural numbers.
More specifically, for every ultrafilter $F$, if a filter $G$ contains $F$ (i.e., $F \subset G$) then $G = F$ must hold,
indicating that an ultrafilter cannot be extended into a larger one.
This leads to another equivalent definition for ultrafilters.

\begin{definition}[Maximal Filter]
A filter $F$ over $\omega$ is a maximal filter if it satisfies:
$$\forall G,\ G\ \mbox{is a filter over}\ \omega,\ F \subset G\ \Longrightarrow\ G = F.$$
\end{definition}

\noindent The definition of ultrafilter and maximal filter are actually equivalent.

\begin{lstlisting}
Definition maxFilter_On_$\omega$ F := Filter_On_$\omega$ F
  /\ ($\forall$ G, Filter_On_$\omega$ G -> F $\subset$ G -> G = F).

Corollary ultraFilter_Equ_maxFilter :
  $\forall$ F, ultraFilter_On_$\omega$ F <-> maxFilter_On_$\omega$ F.
\end{lstlisting}

Before further discussing ultrafilters, we need to review the definitions of finite and infinite in MK.
A set is finite, which generally means that there exists a 1-1 function between this set and a natural number.
For example, the set $\{2,3\}$ is finite because $\{(0,2),(1,3)\}$ is a 1-1 function between $\{2,3\}$ and $\{0,1\}$.
Note that $\{0,1\}$ (i.e., $\{\emptyset,\{\emptyset\}\}$) exactly represents 2 in MK.

However, it is a little inconvenient to describe finite through 1-1 functions for the reason that
one must find an exact natural number for the target set.
Thus in MK, the cardinality function $P = \{(u,v) : u \approx v \wedge v \in C\}$, whose domain is universe $\mathcal{U}$ and range is $C$, is introduced.
$C$ is the class of all cardinal numbers, it is a proper class.
In short, the concept of a cardinal number is a measure of the size of a set.
For finite sets, the cardinal number can accurately describe the number of their elements, and every natural number is a cardinal number itself.
Thus, for each set $x$, $P(x)$ represents the cardinal number equivalent to $x$; it is the cardinality of $x$,
and $x$ is finite means that $P(x) \in \omega$, namely the number of members in $x$ could be accurately depicted with a natural number.

\begin{lstlisting}
Definition C := \{ $\lambda$ x, Cardinal_Number x \}.
Definition P := \{\ $\lambda$ x y, x $\approx$ y /\ y $\in$ C \}\.
Definition Finite x := P[x] $\in$ $\omega$.
\end{lstlisting}

Now let us return to the theme, with the definition of finite, we can give the definition of free ultrafilter.

\begin{definition}[Free Ultrafilter]
An ultrafilter $F$ over $\omega$ is a free ultrafilter if it satisfies:
$$\forall a,\ a\ \mbox{is a finite subset of}\ \omega\ \Longrightarrow\ a \notin F.$$
Non-free ultrafilters represent the ultrafilters that are not free ones.
\end{definition}

\noindent According to this definition, every finite subset of $\omega$ is not a member of any free ultrafilters.
That is to say, the members of a free ultrafilter are all infinite sets.

\begin{lstlisting}
Definition free_ultraFilter_On_$\omega$ F := ultraFilter_On_$\omega$ F
  /\ ($\forall$ a, a $\subset$ $\omega$ -> Finite a -> a $\notin$ F).

Corollary free_ultraFilter_infinite : $\forall$ F, free_ultraFilter_On_$\omega$ F
  -> ($\forall$ x, x $\in$ F -> $\sim$ Finite x).
\end{lstlisting}

Free ultrafilters are usually hard to construct.
The non-principal arithmetical ultrafilters used to establish reals are free ultrafilters.
To find free ultrafilters, a specific filter is required:
$$F_{\sigma} = \{a : a \subset \omega \wedge \omega \sim a\ \mbox{is finite} \}.$$
$F_{\sigma}$ is called the Fr\'{e}chet filter;
it is a filter rather than an ultrafilter but free ultrafilters are closely related to it.

\begin{proposition}
For an ultrafilter $F$, $F$ is a free ultrafilter if and only if $F_{\sigma}$ is contained in $F$.
\end{proposition}

\begin{lstlisting}
Definition F$\sigma$ := \{ $\lambda$ a, a $\subset$ $\omega$ /\ Finite ($\omega$ $\sim$ a) \}.

Corollary F$\sigma$_is_just_Filter : Filter_On_$\omega$ F$\sigma$ /\ $\sim$ ultraFilter_On_$\omega$ F$\sigma$
  /\ ($\forall$ a, a $\subset$ $\omega$ -> Finite a -> a $\notin$ F$\sigma$).

Theorem FT2 : $\forall$ F, ultraFilter_On_$\omega$ F -> free_ultraFilter_On_$\omega$ F <-> F$\sigma$ $\subset$ F.
\end{lstlisting}

With the use of the Extension Principle of Filter
(a consequence of Axiom of Choice, also called the Ultrafilter Theorem in \cite{Tho1973})
stating that every filter over $\omega$ can be extended to an ultrafilter,
we can extend the Fr\'{e}chet filter to obtain a free ultrafilter.

Another type of ultrafilter is the principal ultrafilter
and there exists a close correspondence between principal ultrafilters and natural numbers.

\begin{definition}[Principal Ultrafilter]
For every $n \in \omega$, the following set
$$\{a : a \subset \omega \wedge n \in a\},$$
denoted as $F_n$, can be proven to be an ultrafilter. Each $F_n$, corresponded to the natural number $n$, is called a principal ultrafilter.
An ultrafilter that is not principal is referred to as a non-principal ultrafilter.
\end{definition}

\noindent Principal ultrafilters are directly constructed and formalizing this definition involves two steps.
The first step is to construct the set $F_n$, and the second is to verify that $F_n$ is indeed an ultrafilter.

\begin{lstlisting}
Definition F n := \{ $\lambda$ a, a $\subset$ $\omega$ /\ n $\in$ a \}.

Corollary Fn_Corollary1 : $\forall$ n, F n = $\Phi$ <-> n $\notin$ $\omega$.
Corollary Fn_Corollary2 : $\forall$ n, ultraFilter_On_$\omega$ (F n) <-> n $\in$ $\omega$.
\end{lstlisting}

\noindent In mathematics, when discussing principal ultrafilters, we typically assume that the $n$ of $F_n$ is a natural number.
However, in code, the type of ``{\tt n}" that serves as a parameter of ``{\tt F n}" is ``{\tt Class}", indicating that $n \in \omega$ may not hold.
This embodies the rigorousness in computer formalization.
Thus we present the situation where $n \notin \omega$ even if it is rarely discussed.

From every natural number $n$, its corresponding principal ultrafilter $F_n$ can be constructed,
and $F_n$ is a compatible combination of all properties of $n$.
Thus two different natural numbers will correspond to different principal ultrafilters of their own,
and every member of a principal ultrafilter must include its corresponding natural number.
We directly show here the formalization of these illuminations:

\begin{lstlisting}
Example Example2_a : $\forall$ n m, n $\in$ $\omega$ -> m $\in$ $\omega$ -> n <> m <-> F n <> F m.
Example Example2_b : $\forall$ a n, a $\subset$ $\omega$ -> n $\in$ $\omega$ -> a $\in$ (F n) <-> n $\in$ a.
\end{lstlisting}

In fact, the set of all principal ultrafilters can serve as the set of natural numbers,
as will be presented in \ref{sec5-1-4}.

The next theorem indicates that ultrafilters can be precisely classified into free ultrafilters and principal ultrafilters.

\begin{theorem}
Every principal ultrafilter is a non-free ultrafilter, and conversely, every non-free ultrafilter is a specific principal ultrafilter.
\end{theorem}

\begin{lstlisting}
Theorem FT3_a : $\forall$ n, n $\in$ $\omega$
  -> ultraFilter_On_$\omega$ (F n) /\ $\sim$ free_ultraFilter_On_$\omega$ (F n).
Theorem FT3_b : $\forall$ F0, ultraFilter_On_$\omega$ F0
  -> $\sim$ free_ultraFilter_On_$\omega$ F0 -> ($\exists$ n, n $\in$ $\omega$ /\ F0 = F n).
\end{lstlisting}

\noindent Thus, principal ultrafilters and non-free ultrafilters are two identical concepts;
similarly, non-principal ultrafilters and free ultrafilters are equivalent.

\subsection{Ultrafilter Transformation}\label{sec4-2}

The transformation of ultrafilters is a method for studying their properties and is crucial for the definition of non-principal arithmetical ultrafilter.
It is the result of the impact of functions and some related properties of functions should be introduced first.

\begin{definition}[Image Set]\label{imageset}
Let $f$ be a function and $A$ be subset of domain $f$, then the image set of $f$ at $A$ is denoted as:
$$f\lceil A\rfloor = \{u : \exists x, u = f(x) \wedge x \in A\}.$$
\end{definition}

\begin{definition}[Preimage Set]\label{preimageset}
Let $f$ be a function and $A$ be subset of range $f$, then the preimage set of $f$ at $A$ is denoted as:
$$f^{-1}\lceil A\rfloor = \{u : u \in domain\ f \wedge f(u) \in A\}.$$
\end{definition}

The image set at $A$ is the set of values of $f$ at members in $A$.
However, the condition $A \subset$ domain $f$ in definition \ref{imageset} is not strictly necessary;
$A$ could be any set such that the image set at $A$ is the set of values at members in $A\ \cap$ domain $f$.
The same principle applies to the definition of preimage set.
Thus the formalization is as follows:

\begin{lstlisting}
Definition ImageSet f A := \{ $\lambda$ u, $\exists$ m, u = f[m] /\ m $\in$ A \}.
Notation "f $\lceil$ A $\rfloor$" := (ImageSet f A)(at level 5) : filter_scope.
Corollary ImageSet_Corollary : $\forall$ f A, Function f
  -> A $\cap$ dom(f) <> $\Phi$ <-> f$\lceil$A$\rfloor$ <> $\Phi$.

Definition PreimageSet f A := \{ $\lambda$ u, u $\in$ dom(f) /\ f[u] $\in$ A \}.
Notation "f $^{-1}$$\lceil$ A $\rfloor$" := (PreimageSet f A)(at level 5) : filter_scope.
Corollary PreimageSet_Corollary : $\forall$ f A, Function f
  -> A $\cap$ ran(f) <> $\Phi$ <-> f$^{-1}$$\lceil$A$\rfloor$ <> $\Phi$.
\end{lstlisting}

\noindent It can be observed that $f\lceil A\rfloor$ is not empty if and only if $A\ \cap$ domain $f$ is not empty,
and $f^{-1}\lceil A\rfloor$ is not empty if and only if $A\ \cap$ range $f$ is not.
These facts can simplify formalization, because $f\lceil A\rfloor$ and $f^{-1}\lceil A\rfloor$ still make sense
when $A$ is not a subset of domain $f$ or range $f$,
and restrictions on the relationship between $A$ and domain $f$ or range $f$ can be removed in some specific situations.
After all, fewer restrictions on description usually mean more convenient application in computer formalization.
Following this idea, here are some properties concerning the image and preimage sets,
most of which do not require $A$ (and $B$) to be a subset of domain $f$ or range $f$, and even $f$ to be a function.

\begin{lstlisting}
Fact ImageSet_Property1 : $\forall$ f A B, f$^{-1}$$\lceil$A $\cup$ B$\rfloor$ = f$^{-1}$$\lceil$A$\rfloor$ $\cup$ f$^{-1}$$\lceil$B$\rfloor$.
Fact ImageSet_Property2 : $\forall$ f A B, f$^{-1}$$\lceil$A $\cap$ B$\rfloor$ = f$^{-1}$$\lceil$A$\rfloor$ $\cap$ f$^{-1}$$\lceil$B$\rfloor$.
Fact ImageSet_Property3 : $\forall$ f A B, f$\lceil$A $\cap$ B$\rfloor$ $\subset$ (f$\lceil$A$\rfloor$ $\cap$ f$\lceil$B$\rfloor$).
Fact ImageSet_Property4 : $\forall$ f A B, f$^{-1}$$\lceil$A $\sim$ B$\rfloor$ = f$^{-1}$$\lceil$A$\rfloor$ $\sim$ f$^{-1}$$\lceil$B$\rfloor$.
Fact ImageSet_Property5 : $\forall$ f A B, (f$\lceil$A$\rfloor$ $\sim$ f$\lceil$B$\rfloor$) $\subset$ f$\lceil$A $\sim$ B$\rfloor$.
Fact ImageSet_Property6 : $\forall$ f A, A $\subset$ dom(f) -> A $\subset$ f$^{-1}$$\lceil$f$\lceil$A$\rfloor$$\rfloor$.
Fact ImageSet_Property7 : $\forall$ f A, Function1_1 f -> A $\subset$ dom(f) -> A = f$^{-1}$$\lceil$f$\lceil$A$\rfloor$$\rfloor$.
Fact ImageSet_Property8_a : $\forall$ f A, f$\lceil$f$^{-1}$$\lceil$A$\rfloor$$\rfloor$ $\subset$ A.
Fact ImageSet_Property8_b : $\forall$ f A, Function f -> A $\subset$ ran(f) -> f$\lceil$f$^{-1}$$\lceil$A$\rfloor$$\rfloor$ = A.
\end{lstlisting}

To conveniently describe ultrafilters, the set $\beta\omega$ that consists of all ultrafilters over $\omega$ is introduced,
such that $F$ is an ultrafilter if and only if $F \in \beta\omega$.
$\beta\omega$ is an important set, within which some research involving filters\cite{FNB2002,NF2005a,NF2005b,Wan1998} is carried out.

\begin{lstlisting}
Definition $\beta$$\omega$ := \{ $\lambda$ u, ultraFilter_On_$\omega$ u \}.
Corollary $\beta$$\omega$_is_Set : Ensemble $\beta$$\omega$.
Corollary $\beta$$\omega$_Corollary : $\forall$ F, F $\in$ $\beta$$\omega$ <-> ultraFilter_On_$\omega$ F.
\end{lstlisting}

The transformation of ultrafilters is also a method for studying $\beta\omega$.
It operates similarly to a function on $\beta\omega$, which can map one ultrafilter to another one.

\begin{definition}[Ultrafilter Transformation]\label{ultrafiltertrans}
Let $F$ be an ultrafilter, $f$ is a function in ${\omega}^{\omega}$, then the following set
$$f\langle F\rangle = \{u : u \subset \omega \wedge f^{-1}\lceil u\rfloor \in F\}$$
can be proven an ultrafilter. $f\langle F\rangle$ is called a transformation of ultrafilter $F$ under function $f$.
\end{definition}

\begin{lstlisting}
Definition Transform f F := \{ $\lambda$ u, u $\subset$ $\omega$ /\ f$^{-1}$$\lceil$u$\rfloor$ $\in$ F \}.
Notation "f $\langle$ F $\rangle$" := (Transform f F)(at level 5) : filter_scope.

Theorem FT5 : $\forall$ F f, F $\in$ $\beta$$\omega$ -> Function f
  -> dom(f) = $\in$ -> ran(f) $\subset$ $\omega$ -> f$\langle$F$\rangle$ $\in$ $\beta$$\omega$.
\end{lstlisting}

According to definition \ref{ultrafiltertrans},
any ultrafilter can be transformed to another one and this is achieved using functions on $\omega$ rather than functions on $\beta\omega$.
Consequently, ultrafilter transformation will exhibit certain properties analogous to those of functions on natural numbers.
These properties can be observed in transformations associated with principal ultrafilters due to their close correspondence with natural numbers.
The next two theorems illustrate the transformations related to principal ultrafilters.

\begin{theorem}\label{the6}
Let $F$ be a principal ultrafilter, $f$ a function in ${\omega}^{\omega}$,
then $f\langle F \rangle$ is a principle ultrafilter as well.
\end{theorem}

\begin{theorem}\label{the7}
Let $F$ be any of an ultrafilter, $f$ a function whose domain is $\omega$ and value is a constant $m \in \omega$,
then $f\langle F \rangle$ is the principal ultrafilter $F_m$.
\end{theorem}

Theorem \ref{the6} states that any principal ultrafilter remains a principal ultrafilter after any transformation.
Meanwhile, theorem \ref{the7} demonstrates that every ultrafilter can be transformed to any other principal ultrafilter with the assistance of constant functions.
Their formalization is presented here:

\begin{lstlisting}
Theorem FT6 : $\forall$ f n, Function f -> dom(f) = $\omega$ -> ran(f) $\subset$ $\omega$
  -> n $\in$ $\omega$ -> f$\langle$F n$\rangle$ = F f[n].
Theorem FT7 : $\forall$ F0 f m, F0 $\in$ $\beta$$\omega$ -> Function f -> dom(f) = $\omega$
  -> m $\in$ $\omega$ -> ($\forall$ n, n $\in$ $\omega$ -> f[n] = m) -> f$\langle$F0$\rangle$ = F m.
\end{lstlisting}

\noindent Indeed, the function $f$ in theorem \ref{the7} does not necessarily need to be a constant function.
If $f$ merely takes a constant value $m$ on a certain member $b$ of $F$ (i.e., $b\in F \wedge \forall n \in b, f[n] = m$), the theorem still holds.
Hence, we can formalize it as:

\begin{lstlisting}
Theorem FT8 : $\forall$ F0 f m, F0 $\in$ $\beta$$\omega$ -> Function f -> dom(f) = $\omega$ -> ran(f) $\subset$ $\omega$
  -> m $\in$ $\omega$ -> ($\exists$ b, b $\in$ F0 /\ ($\forall$ n, n $\in$ b -> f(n) = m)) -> f$\langle$F0$\rangle$ = F m.
\end{lstlisting}

As a method for studying ultrafilters, ultrafilter transformation exhibits various properties\cite{Wan2016}.
For arithmetical ultrafilter and FMCR, our formalization above is adequate.

\subsection{Arithmetical Ultrafilter}\label{sec4-3}

\begin{definition}\label{Fequ}
Let $F$ be an ultrafilter over $\omega$,
then two functions $f$ and $g$ in ${\omega}^{\omega}$ are $F$-equivalent if and only if
$$\{u : u \in \omega \wedge f(u) = g(u)\}\ \in\ F,$$
namely $f$ and $g$ are equal on a set in $F$, which is denoted as $f =_F g$.
\end{definition}

\begin{lstlisting}
Definition NearlyEqual f F g := Function f /\ Function g
  /\ dom(f) = $\omega$ /\ dom(g) = $\omega$ /\ ran(f) $\subset$ $\omega$ /\ ran(g) $\subset$ $\omega$
  /\ F $\in$ $\beta$$\omega$ /\ \{ $\lambda$ u, u $\in$ $\omega$ /\ f[u] = g[u] \} $\in$ F.
Notation "f =_ F g" := (NearlyEqual f F g)
  (at level 10, F at level 9) : filter_scope.
\end{lstlisting}

Definition \ref{Fequ} establishes a connection between the equivalence relation in $\omega$ and that in $\beta\omega$,
offering a means to explore the arithmetical relationship between them.
This connection paves the way for the definition of arithmetical ultrafilter.

\begin{definition}\label{arithf}
An ultrafilter $F$ over $\omega$ is an arithmetical ultrafilter if and only if
\begin{equation}\label{equa10}
\forall f,g \in {\omega}^{\omega},\ f\langle F\rangle = g\langle F\rangle\ \Longleftrightarrow\ f =_F g.
\end{equation}
\end{definition}

\noindent The definition of arithmetical ultrafilter stipulates that $f\langle F\rangle = g\langle F\rangle$ and $f =_F g$ are equivalent in equation (\ref{equa10}),
while it is can be verified that for every ultrafilter $F$, $f =_F g$ implies $f\langle F\rangle = g\langle F\rangle$.
Consequently, our formal definition for arithmetical ultrafilter is slightly simplified:

\begin{lstlisting}
Theorem FT9 : $\forall$ f F g, f =_F g -> f$\langle$F$\rangle$ = g$\langle$F$\rangle$.
Definition Arithmetical_ultraFilter F0 := F0 $\in$ $\beta$$\omega$
  /\ $\forall$ f g, Function f -> Function g -> dom(f) = $\omega$ -> dom(g) = $\omega$
    -> ran(f) $\subset$ $\omega$ -> ran(g) $\subset$ $\omega$ -> f$\langle$F0$\rangle$ = g$\langle$F0$\rangle$ -> f =_F0 g.
\end{lstlisting}

The next two theorems illustrate some basic properties of arithmetical ultrafilters.

\begin{theorem}\label{tre8}
For every $n \in \omega$, $F_n$ is an arithmetical ultrafilter.
\end{theorem}

\begin{theorem}\label{tre9}
Let $F$ be an arithmetical ultrafilter, $f$ be a function in ${\omega}^{\omega}$, then $f\langle F\rangle$ is an arithmetical ultrafilter as well.
\end{theorem}

\begin{lstlisting}
Theorem FT10 : $\forall$ n, n $\in$ $\omega$ -> Arithmetical_ultraFilter (F n).
Theorem FT11 : $\forall$ F0 f, Arithmetical_ultraFilter F0 -> Function f
  -> dom(f) = $\omega$ -> ran(f) $\subset$ $\omega$ -> Arithmetical_ultraFilter (f$\langle$F0$\rangle$).
\end{lstlisting}

\noindent Theorem \ref{tre8} illustrates that all principal ultrafilters are arithmetical ultrafilters
and theorem \ref{tre9} demonstrates that any arithmetical ultrafilter remains an arithmetical ultrafilter after undergoing any transformation.
This raises the question of whether there exist non-principal arithmetical ultrafilters.
In 2003, Bartoszynski and Shelah claimed that the existence of non-principal arithmetical ultrafilters could not be proven solely with ZFC set theory\cite{BS2003}.
Although the conclusion was later retracted due to some errors, it is generally believed to be correct\cite{Wan2016,Wan2018}.
To find non-principal arithmetical ultrafilters, additional hypotheses of set theory are required.

\begin{figure}[h]%
\centering
\includegraphics[width=1\textwidth]{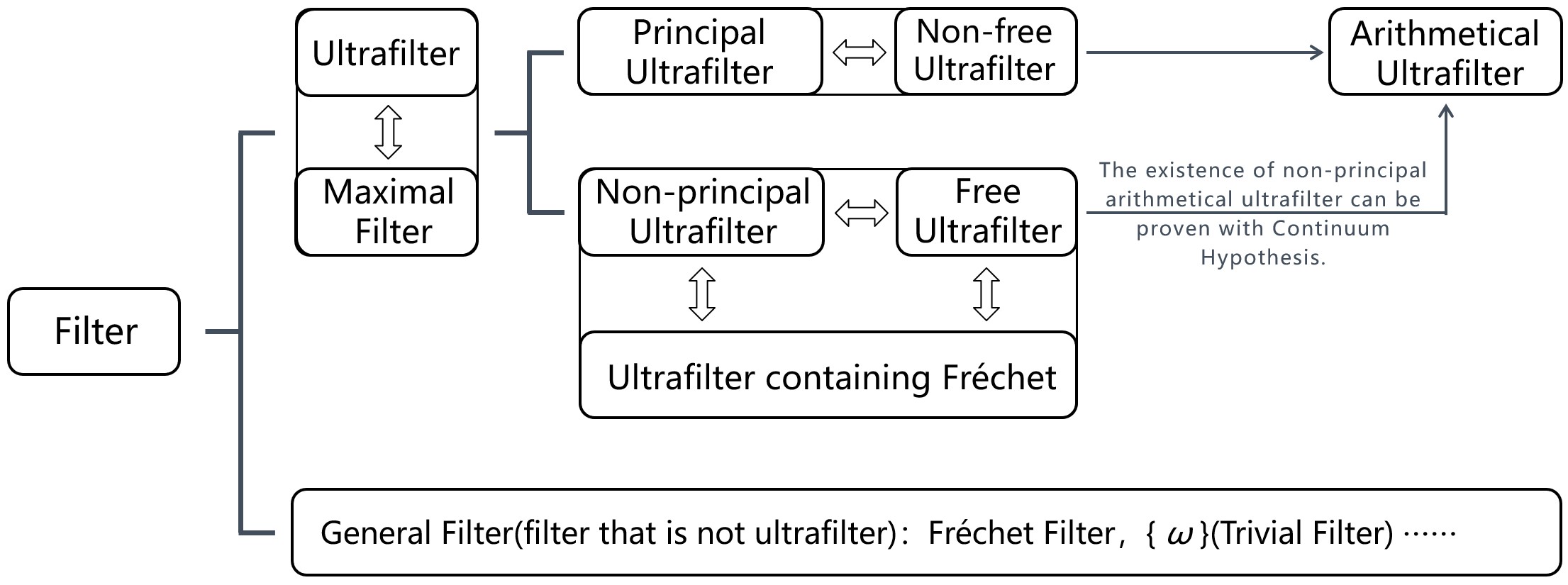}
\caption{relationships among filters}\label{fig1}
\end{figure}

Up to this point, we have recognized numerous concepts related to filters, and figure \ref{fig1} presents the relationships among them.
For further concepts on filters and arithmetical ultrafilters, please refer to \cite{CN1974,Wan2016}.

\section{Implementation of Reals}\label{sec5}

This section is the core of FMCR, the formalization of the entire process of $\omega$ to $^{*}{\mathbf N}$, $^{*}{\mathbf N}$ to $^{*}{\mathbf Z}$, $^{*}{\mathbf Z}$ to $^{*}{\mathbf Q}$,
and finally $^{*}{\mathbf Q}$ to $\mathbf R$ will be comprehensively and meticulously presented.

\subsection{Extension from $\omega$ to $^{*}{\mathbf N}$}\label{sec5-1}

\subsubsection{Construction of $^{*}{\mathbf N}$}\label{sec5-1-1}

For every non-principal arithmetical ultrafilter $F0$\footnote {Note that $F0$ and $F_0$(principal ultrafilter corresponding to $0$) are two different symbols.},
the set $^{*}{\mathbf N}$ is constructed as follows:
$$^{*}{\mathbf N} = \{u : \exists f, f \in \omega^{\omega} \wedge u = f\langle F0\rangle\}.$$
It can be observed, according to theorem \ref{tre9}, that $^{*}{\mathbf N}$ is the set consisting of all possible transformations from $F0$ and, of course, a subset of $\beta\omega$.

\begin{lstlisting}
Definition N' F0 :=  \{ $\lambda$ u, $\exists$ f, Function f
  /\ dom(f) = $\omega$ /\ ran(f) $\subset$ $\omega$ /\ u = f$\langle$F0$\rangle$ \}.

Corollary N'_Subset_$\beta\omega$ : $\forall$ F0, F0 $\in$ $\beta\omega$ -> (N' F0) $\subset$ $\beta\omega$.
Corollary N'_is_Set : $\forall$ F0, F0 $\in$ $\beta\omega$ -> Ensemble (N' F0).
\end{lstlisting}

In formal descriptions, $F0$ is set as a parameter of $^{*}{\mathbf N}$ rather than declared as a global parameter like these:

\begin{lstlisting}
Parameter F0 : $\mbox{\color[RGB]{0,0,0}Class}$.
Parameter NPAUF : Arithmetical_ultraFilter F0 /\ ($\forall$ m, F0 <> F m).
\end{lstlisting}

\noindent where $F0$ is admitted as a specific non-principal arithmetical ultrafilter using the key word {\tt Parameter},
with which $^{*}{\mathbf N}$ can be formally described directly without any parameters as follows:

\begin{lstlisting}
Definition N' :=  \{ $\lambda$ u, $\exists$ f, Function f
  /\ dom(f) = $\omega$ /\ ran(f) $\subset$ $\omega$ /\ u = f$\langle$F0$\rangle$ \}.
\end{lstlisting}

\noindent Even though this description seems clearer and more brief, we adopt the parameterized description of $^{*}{\mathbf N}$({\tt N' F0} in codes).
Two main reasons support this choice: firstly, from a mathematical perspective,
such $^{*}{\mathbf N}$ can be constructed from any non-principal arithmetical ultrafilter $F0$,
and $F0$ dose not need to be set as a fixed constant in code;
secondly, during the formalization process, we observed that the discussion regarding $^{*}{\mathbf N}$ remains meaningful
under the condition that $F0$ is any arithmetical ultrafilter, and sometimes even just an ultrafilter.
And the same rationale applies to the formalization of $^{*}{\mathbf Z}$, $^{*}{\mathbf Q}$ and ${\mathbf R}$.

Therefore, in our formalization, the conditions for $F0$ are tailored to specific situations to streamline the premises related to $F0$ in formal descriptions.
For instance, when verifying that $^{*}{\mathbf N}$ is a subset of $\beta\omega$, $F0$ is only required to ba an ultrafilter, not necessarily an arithmetical ultrafilter.

\begin{lstlisting}
Corollary N'_Subset_$\beta\omega$ : $\forall$ F0, F0 $\in$ $\beta\omega$ -> (N' F0) $\subset$ $\beta\omega$.
Corollary N'_is_Set : $\forall$ F0, F0 $\in$ $\beta\omega$ -> Ensemble (N' F0).
\end{lstlisting}

\noindent To verify that every principal ultrafilter is in $^{*}{\mathbf N}$,
which is crucial to finding connections between $^{*}{\mathbf N}$ and $\omega$ in \ref{sec5-1-4},
$F0$ is also only required to ba an ultrafilter.

\begin{lstlisting}
Fact Fn_in_N' : $\forall$ n F0, F0 $\in$ $\beta\omega$ -> n $\in$ $\omega$ -> (F n) $\in$ (N' F0).
\end{lstlisting}

\noindent More obvious and specific examples will be seen in following sections.

\subsubsection{Order on $^{*}{\mathbf N}$}\label{sec5-1-2}

\begin{definition}[Order on $^{*}{\mathbf N}$]
Let $^{*}{\mathbf N}$ be constructed with non-principal arithmetical ultrafilter $F0$,
$f, g \in {\omega}^{\omega}$, $u, v \in \omega$ and $u = f\langle F0\rangle$, $v = g\langle F0\rangle$,
then the order of $u$ and $v$ is defined as follows:
$$u <_{^{*}{\mathbf N}} v\ \Longleftrightarrow\ f\langle F0\rangle\ <_{^{*}{\mathbf N}}\ g\langle F0\rangle\ \Longleftrightarrow\ \{n : f(n) \in g(n) \} \in F0.$$
\end{definition}

It can be verified that the definition of this order relation is reasonable, because for any functions $f_1,f_2,g_1,g_2 \in {\omega}^{\omega}$, if
$$u = f_1\langle F0\rangle = f_2\langle F0\rangle,\ v = g_1\langle F0\rangle = g_2\langle F0\rangle,$$
then
$$u <_{^{*}{\mathbf N}} v\ \Longleftrightarrow\ f_1\langle F0\rangle\ <_{^{*}{\mathbf N}}\ g_1\langle F0\rangle\ \Longleftrightarrow\
f_2\langle F0\rangle\ <_{^{*}{\mathbf N}}\ g_2\langle F0\rangle.$$
Namely, the order relation on $^{*}{\mathbf N}$ is independent of the selection of functions.

\begin{lstlisting}
Definition N'_Ord F0 := \{\ $\lambda$ u v, $\forall$ f g, Function f -> Function g
  -> dom(f) = $\omega$ -> dom(g) = $\omega$ -> ran(f) $\subset$ $\omega$ -> ran(g) $\subset$ $\omega$
  -> u = f$\langle$F0$\rangle$ -> v = g$\langle$F0$\rangle$ -> \{ $\lambda$ w, f[w] $\in$ g[w] \} $\in$ F0 \}\.
Notation "u <_ F0 v" := ([u,v] $\in$ (N'_Ord F0))
  (at level 10, F0 at level 9) : n'_scope.

Corollary N'_Ord_Corollary : $\forall$ F0 f g, Arithmetical_ultraFilter F0
  -> Function f -> Function g -> dom(f) = $\omega$
  -> dom(g) = $\omega$ -> ran(f) $\subset$ $\omega$ -> ran(g) $\subset$ $\omega$
  -> f$\langle$F0$\rangle$ <_F0 g$\langle$F0$\rangle$ <-> \{ $\lambda$ w, f[w] $\in$ g[w] \} $\in$ F0.
Corollary N'_Ord_Reasonable : $\forall$ F0 f1 f2 g1 g2, Arithmetical_ultraFilter F0
  -> Function f1 -> Function f2 -> Function g1 -> Function g2
  -> dom(f1) = $\omega$ -> dom(f2) = $\omega$ -> dom(g1) = $\omega$ -> dom(g2) = $\omega$
  -> ran(f1) $\subset$ $\omega$ -> ran(f2) $\subset$ $\omega$ -> ran(g1) $\subset$ $\omega$ -> ran(g2) $\subset$ $\omega$
  -> f1$\langle$F0$\rangle$ = f2$\langle$F0$\rangle$ -> g1$\langle$F0$\rangle$ = g2$\langle$F0$\rangle$
  -> f1$\langle$F0$\rangle$ <_F0 g1$\langle$F0$\rangle$ <-> f2$\langle$F0$\rangle$ <_F0 g2$\langle$F0$\rangle$.
\end{lstlisting}

The principal ultrafilters $F_0$ and $F_1$ are two important members of $^{*}{\mathbf N}$,
analogous to the counterparts $0$ and $1$ in $\omega$, which will be illustrated in \ref{sec5-1-3} and \ref{sec5-1-4}.
Here are two trivial order properties about $F_0$ and $F_1$: $F_0 <_{^{*}{\mathbf N}} F_1$ and $F0$ is the first member of $^{*}{\mathbf N}$
(i.e., there is no member $u$ in $^{*}{\mathbf N}$ such that $u <_{^{*}{\mathbf N}} F_0$).

\begin{lstlisting}
Fact N'0_lt_N'1 : $\forall$ F0, Arithmetical_ultraFilter F0 -> (F 0) <_F0 (F 1).
Fact N'0_is_FirstMember : $\forall$ F0 u, Arithmetical_ultraFilter F0
  -> u $\in$ (N' F0) -> $\sim$ u <_F0 (F 0).
\end{lstlisting}

Besides, the order on $^{*}{\mathbf N}$ is a total order, meaning that for any $u,v,w \in {^{*}{\mathbf N}}$, it satisfies:

1) Irreflexivity: $\sim u <_{^{*}{\mathbf N}} u$,

2) Transitivity: $u <_{^{*}{\mathbf N}} v \wedge v <_{^{*}{\mathbf N}} w\ \Longrightarrow\ u <_{^{*}{\mathbf N}} w$,

3) Trichotomy: $u <_{^{*}{\mathbf N}} v \vee v <_{^{*}{\mathbf N}} u \vee u = v$.

\begin{lstlisting}
Fact N'_Ord_irReflex : $\forall$ F0 u v, Arithmetical_ultraFilter F0
  -> u $\in$ (N' F0) -> v $\in$ (N' F0) -> u = v -> $\sim$ u <_F0 v.
Fact N'_Ord_Trans : $\forall$ F0 u v w, F0 $\in$ $\beta\omega$ -> u $\in$ (N' F0) -> v $\in$ (N' F0)
  -> w $\in$ (N' F0) -> u <_F0 v -> v <_F0 w -> u <_F0 w.
Fact N'_Ord_Tri : $\forall$ F0, Arithmetical_ultraFilter F0
  -> Connect (N'_Ord F0) (N' F0).
\end{lstlisting}

In code, ``{\tt Connect}" is the formalization of the concept ``connect" from MK that is equivalent to ``trichotomy".\cite{Kel1955,SY2020}

\subsubsection{Operations on $^{*}{\mathbf N}$}\label{sec5-1-3}

Before defining operations on $^{*}{\mathbf N}$, the operations among functions in ${\omega}^{\omega}$ is required.

\begin{definition}[Addition among Functions]
For any functions $f$ and $g$ in ${\omega}^{\omega}$, the sum of $f$ and $g$ is:
$$f + g = \{(u,v) : u \in \omega \wedge v = f(u) + g(u)\},$$
where $f(u) + g(u)$ represents the addition on $\omega$.
\end{definition}

\begin{lstlisting}
Definition Function_Plus f g := \{\ $\lambda$ u v, u $\in$ $\omega$ /\ v = f[u] + g[u] \}\.
Notation "f + g" := (Function_Plus f g) : fu_scope.
\end{lstlisting}

According to this definition, $f + g$ is the sum of values of $f$ and $g$ at same preimages.
Thus $f + g$ is also a function in ${\omega}^{\omega}$ and for every $n \in \omega$, $(f + g)(n) = f(n) + g(n)$.

\begin{lstlisting}
Corollary Function_Plus_Corollary1 : $\forall$ f g, Function f -> Function g
  -> dom(f) = $\omega$ -> dom(g) = $\omega$ -> ran(f) $\subset$ $\omega$ -> ran(g) $\subset$ $\omega$
 -> Function (f + g) /\ dom(f + g) = $\omega$ /\ ran(f + g) $\subset$ $\omega$.
Corollary Function_Plus_Corollary2 : $\forall$ f g n, Function f -> Function g
  -> dom(f) = $\omega$ -> dom(g) = $\omega$ -> ran(f) $\subset$ $\omega$ -> ran(g) $\subset$ $\omega$
  -> (f + g)[n] = (f[n] + g[n])$\mbox{\color[RGB]{0,100,0}\%}$$\omega$.
\end{lstlisting}

\noindent Thanks to the scope mechanism in Coq, we can use the same notation ``{\tt \_ + \_}" in different situations
and ``{\tt (f[n] + g[n]){\color[RGB]{0,100,0}\%}$\omega$}" represents the addition on $\omega$, which is realized by:

\begin{lstlisting}
Declare Scope $\omega$_scope.
Delimit Scope $\omega$_scope with $\omega$.
Notation "m + n" := (Plus m n) : $\omega$_scope.
\end{lstlisting}

The multiplication among functions is totally similar to addition, it is the product of two function values in ${\omega}^{\omega}$ at same preimages.

\begin{definition}[Multiplication among Functions]
For any functions $f$ and $g$ in ${\omega}^{\omega}$, the product of $f$ and $g$ is:
$$f \cdot g = \{(u,v) : u \in \omega \wedge v = f(u) \cdot g(u)\},$$
where $f(u) \cdot g(u)$ represents the multiplication on $\omega$.
\end{definition}

\begin{lstlisting}
Definition Function_Mult f g := \{\ $\lambda$ u v, u $\in$ $\omega$ /\ v = f[u] $\cdot$ g[u] \}\.
Notation "f $\cdot$ g" := (Function_Mult f g)(at level 10) : fu_scope.

Corollary Function_Mult_Corollary1 : $\forall$ f g, Function f -> Function g
  -> dom(f) = $\omega$ -> dom(g) = $\omega$ -> ran(f) $\subset$ $\omega$ -> ran(g) $\subset$ $\omega$
  -> Function (f $\cdot$ g) /\ dom(f $\cdot$ g) = $\omega$ /\ ran(f $\cdot$ g) $\subset$ $\omega$.
Corollary Function_Mult_Corollary2 : $\forall$ f g n, Function f -> Function g
  -> dom(f) = $\omega$ -> dom(g) = $\omega$ -> ran(f) $\subset$ $\omega$ -> ran(g) $\subset$ $\omega$
  -> (f $\cdot$ g)[n] = (f[n] $\cdot$ g[n])$\mbox{\color[RGB]{0,100,0}\%}$$\omega$.
\end{lstlisting}

Now the addition and multiplication on $^{*}{\mathbf N}$ can be defined using the corresponding operations among functions.

\begin{definition}[Addition on $^{*}{\mathbf N}$]
Let $^{*}{\mathbf N}$ be constructed with non-principal arithmetical ultrafilter $F0$,
$f, g \in {\omega}^{\omega}$, $u, v \in {^{*}{\mathbf N}}$ and $u = f\langle F0\rangle$, $v = g\langle F0\rangle$,
then the sum of $u$ and $v$ is defined as follows:
$$u +_{^{*}{\mathbf N}} v\ =\ f\langle F0\rangle\ +_{^{*}{\mathbf N}}\ g\langle F0\rangle\ =\ (f + g)\langle F0\rangle,$$
where $f + g$ represents the addition among functions in ${\omega}^{\omega}$.
\end{definition}

The definition of this addition is reasonable, because for any functions $f_1,f_2,g_1,g_2 \in {\omega}^{\omega}$, if
$$u = f_1\langle F0\rangle = f_2\langle F0\rangle,\ v = g_1\langle F0\rangle = g_2\langle F0\rangle,$$
then it can be verified that
$$u +_{^{*}{\mathbf N}} v\ =\ f_1\langle F0\rangle\ +_{^{*}{\mathbf N}}\ g_1\langle F0\rangle\ =\
f_2\langle F0\rangle\ +_{^{*}{\mathbf N}}\ g_2\langle F0\rangle.$$
Namely, the result of $u +_{^{*}{\mathbf N}} v$ is unique and independent of the selection of functions.

\begin{lstlisting}
Definition N'_Plus F0 u v := $\cap$(\{ $\lambda$ a, $\forall$ f g, Function f
  -> Function g -> dom(f) = $\omega$ -> dom(g) = $\omega$ -> ran(f) $\subset$ $\omega$
  -> ran(g) $\subset$ $\omega$ -> u = f$\langle$F0$\rangle$ -> v = g$\langle$F0$\rangle$ -> a = (f + g)$\langle$F0$\rangle$ \}).
Notation "u +_ F0 v" := (N'_Plus F0 u v)
  (at level 10, F0 at level 9) : n'_scope.

Lemma N'_Plus_Reasonable_Lemma : $\forall$ F0 f g g1, Function f
  -> Function g -> Function g1 -> dom(f) = $\omega$ -> dom(g) = $\omega$
  -> dom(g1) = $\omega$ -> ran(f) $\subset$ $\omega$ -> ran(g) $\subset$ $\omega$ -> ran(g1) $\subset$ $\omega$
  -> g <> g1 -> Arithmetical_ultraFilter F0 -> g$\langle$F0$\rangle$ = g1$\langle$F0$\rangle$
  -> (f + g)$\langle$F0$\rangle$ = (f + g1)$\langle$F0$\rangle$.
Corollary N'_Plus_Reasonable : $\forall$ F0 u v, Arithmetical_ultraFilter F0
  -> u $\in$ (N' F0) -> v $\in$ (N' F0) -> $\exists$ a, Arithmetical_ultraFilter a
    /\ [a] = \{ $\lambda$ a, $\forall$ f g, Function f -> Function g -> dom(f) = $\omega$
    -> dom(g) = $\omega$ -> ran(f) $\subset$ $\omega$ -> ran(g) $\subset$ $\omega$ -> u = f$\langle$F0$\rangle$
    -> v = g$\langle$F0$\rangle$ -> a = (f + g)$\langle$F0$\rangle$ \}.

Corollary N'_Plus_Corollary : $\forall$ F0 f g, Arithmetical_ultraFilter F0
  -> Function f -> Function g -> dom(f) = $\omega$ -> dom(g) = $\omega$
  -> ran(f) $\subset$ $\omega$ -> ran(g) $\subset$ $\omega$ -> f$\langle$F0$\rangle$ +_F0 g$\langle$F0$\rangle$ = (f + g)$\langle$F0$\rangle$.
\end{lstlisting}

In code, similar to formalization of operations on $\omega$ in \ref{sec3-6},
for any $u,v$ in $^{*}{\mathbf N}$ the result of $u +_{^{*}{\mathbf N}} v$ is formalized as the form of ``$\cap\{\cdots\}$" due to its uniqueness.
``{\tt N'\_Plus\_Reasonable}" guarantees that ``{\tt N'\_Plus}" can be defined as the intersection of a singleton that includes the unique result of $u +_{^{*}{\mathbf N}} v$.
Moreover, according to the definition of $^{*}{\mathbf N}$, for every $x$ in $^{*}{\mathbf N}$,
there exists a function $h$ in ${\omega}^{\omega}$ that satisfies $x = h\langle F0\rangle$,
where $F0$ is the non-principal arithmetical ultrafilter used to construct $^{*}{\mathbf N}$.
Thus when $u$ and $v$ are specifically described as $f\langle F0\rangle$ and $g\langle F0\rangle$,
``{\tt N'\_Plus\_Corollary}" guarantees that $u +_{^{*}{\mathbf N}} v\ =\ f\langle F0\rangle\ +_{^{*}{\mathbf N}}\ g\langle F0\rangle\ =\ (f + g)\langle F0\rangle$.
With these, the addition on $^{*}{\mathbf N}$ is entirely defined.

The definition and formalization of multiplication follow the same idea as that of addition.
All we need to do is replacing addition among functions with multiplication.

\begin{definition}[Multiplication on $^{*}{\mathbf N}$]
Let $^{*}{\mathbf N}$ be constructed with non-principal arithmetical ultrafilter $F0$,
$f, g \in {\omega}^{\omega}$, $u, v \in {^{*}{\mathbf N}}$ and $u = f\langle F0\rangle$, $v = g\langle F0\rangle$,
then the product of $u$ and $v$ is defined as follows:
$$u \cdot_{^{*}{\mathbf N}} v\ =\ f\langle F0\rangle\ \cdot_{^{*}{\mathbf N}}\ g\langle F0\rangle\ =\ (f \cdot g)\langle F0\rangle,$$
where $f \cdot g$ represents the multiplication among functions in ${\omega}^{\omega}$.
\end{definition}

The definition of this multiplication is reasonable, because for any functions $f_1,f_2,g_1,g_2 \in {\omega}^{\omega}$, if
$$u = f_1\langle F0\rangle = f_2\langle F0\rangle,\ v = g_1\langle F0\rangle = g_2\langle F0\rangle,$$
then it can be verified that
$$u \cdot_{^{*}{\mathbf N}} v\ =\ f_1\langle F0\rangle\ \cdot_{^{*}{\mathbf N}}\ g_1\langle F0\rangle\ =\
f_2\langle F0\rangle\ \cdot_{^{*}{\mathbf N}}\ g_2\langle F0\rangle.$$
Namely, the result of $u \cdot_{^{*}{\mathbf N}} v$ is unique and independent of the selection of functions.

\begin{lstlisting}
Definition N'_Mult F0 u v := $\cap$(\{ $\lambda$ a, $\forall$ f g, Function f
  -> Function g -> dom(f) = $\omega$ -> dom(g) = $\omega$ -> ran(f) $\subset$ $\omega$
  -> ran(g) $\subset$ $\omega$ -> u = f$\langle$F0$\rangle$ -> v = g$\langle$F0$\rangle$ -> a = (f $\cdot$ g)$\langle$F0$\rangle$ \}).
Notation "u $\cdot$_ F0 v" :=
  (N'_Mult F0 u v)(at level 10, F0 at level 9) : n'_scope.

Lemma N'_Mult_Reasonable_Lemma : $\forall$ F0 f g g1, Function f
  -> Function g -> Function g1 -> dom(f) = $\omega$ -> dom(g) = $\omega$
  -> dom(g1) = $\omega$ -> ran(f) $\subset$ $\omega$ -> ran(g) $\subset$ $\omega$ -> ran(g1) $\subset$ $\omega$
  -> g <> g1 -> Arithmetical_ultraFilter F0 -> g$\langle$F0$\rangle$ = g1$\langle$F0$\rangle$
  -> (f $\cdot$ g)$\langle$F0$\rangle$ = (f $\cdot$ g1)$\langle$F0$\rangle$.
Corollary N'_Mult_Reasonable : $\forall$ F0 u v, Arithmetical_ultraFilter F0
  -> u $\in$ (N' F0) -> v $\in$ (N' F0) -> $\exists$ a, Arithmetical_ultraFilter a
    /\ [a] = \{ $\lambda$ a, $\forall$ f g, Function f -> Function g -> dom(f) = $\omega$
    -> dom(g) = $\omega$ -> ran(f) $\subset$ $\omega$ -> ran(g) $\subset$ $\omega$ -> u = f$\langle$F0$\rangle$
    -> v = g$\langle$F0$\rangle$ -> a = (f $\cdot$ g)$\langle$F0$\rangle$ \}.

Corollary N'_Mult_Corollary : $\forall$ F0 f g, Arithmetical_ultraFilter F0
  -> Function f -> Function g -> dom(f) = $\omega$ -> dom(g) = $\omega$
  -> ran(f) $\subset$ $\omega$ -> ran(g) $\subset$ $\omega$ -> f$\langle$F0$\rangle$ $\cdot$_F0 g$\langle$F0$\rangle$ = (f $\cdot$ g)$\langle$F0$\rangle$.
\end{lstlisting}

Next we present the formalization of properties of the operations defined above.

The operations on $^{*}{\mathbf N}$ are closed in $^{*}{\mathbf N}$.
This means that for any $u,v$ in $^{*}{\mathbf N}$,
both $u +_{^{*}{\mathbf N}} v$ and $u \cdot_{^{*}{\mathbf N}} v$ belong to $^{*}{\mathbf N}$.

\begin{lstlisting}
Fact N'_Plus_in_N' : $\forall$ F0 u v, Arithmetical_ultraFilter F0
  -> u $\in$ (N' F0) -> v $\in$ (N' F0) -> (u +_F0 v) $\in$ (N' F0).
Fact N'_Mult_in_N' : $\forall$ F0 u v, Arithmetical_ultraFilter F0
  -> u $\in$ (N' F0) -> v $\in$ (N' F0) -> (u $\cdot$_F0 v) $\in$ (N' F0).
\end{lstlisting}

The principal ultrafilter $F_0$ and $F_1$ play significant roles in $^{*}{\mathbf N}$, as mentioned in \ref{sec5-1-2}.
Specifically, $F_0$ serves as the zero element of $^{*}{\mathbf N}$ and $F_1$ as the identity element.
In other words, every $u$ in $^{*}{\mathbf N}$ satisfies: $u +_{^{*}{\mathbf N}} F_0 = u$, $u \cdot_{^{*}{\mathbf N}} F_0 = F_0$
and $u \cdot_{^{*}{\mathbf N}} F_1 = u$.

\begin{lstlisting}
Fact N'_Plus_Property : $\forall$ F0 u, Arithmetical_ultraFilter F0
  -> u $\in$ (N' F0) -> u +_F0 (F 0) = u.
Fact N'_Mult_Property1 : $\forall$ F0 u, Arithmetical_ultraFilter F0
  -> u $\in$ (N' F0) -> u $\cdot$_F0 (F 0) = F 0.
Fact N'_Mult_Property2 : $\forall$ F0 u, Arithmetical_ultraFilter F0
  -> u $\in$ (N' F0) -> u $\cdot$_F0 (F 1) = u.
Fact N'_Mult_Property3 : $\forall$ F0 u v, Arithmetical_ultraFilter F0
  -> u $\in$ (N' F0) -> v $\in$ (N' F0) -> u $\cdot$_F0 v = F 0 -> u = F 0 \/ v = F 0.
\end{lstlisting}

Both addition and multiplication on $^{*}{\mathbf N}$ satisfy commutative law, associative law and cancellation law;
besides, multiplication satisfies distributive law.

\begin{lstlisting}
Fact N'_Plus_Commutation : $\forall$ F0 u v, Arithmetical_ultraFilter F0
  -> u $\in$ (N' F0) -> v $\in$ (N' F0) -> u +_F0 v = v +_F0 u.
Fact N'_Plus_Association : $\forall$ F0 u v w, Arithmetical_ultraFilter F0
  -> u $\in$ (N' F0) -> v $\in$ (N' F0) -> w $\in$ (N' F0)
  -> (u +_F0 v) +_F0 w = u +_F0 (v +_F0 w).
Fact N'_Plus_Cancellation : $\forall$ F0 u v w, Arithmetical_ultraFilter F0
  -> u $\in$ (N' F0) -> v $\in$ (N' F0) -> w $\in$ (N' F0)
  -> u +_F0 v = u +_F0 w -> v = w.

Fact N'_Mult_Commutation : $\forall$ F0 u v, Arithmetical_ultraFilter F0
  -> u $\in$ (N' F0) -> v $\in$ (N' F0) -> u $\cdot$_F0 v = v $\cdot$_F0 u.
Fact N'_Mult_Association : $\forall$ F0 u v w, Arithmetical_ultraFilter F0
  -> u $\in$ (N' F0) -> v $\in$ (N' F0) -> w $\in$ (N' F0)
  -> (u $\cdot$_F0 v) $\cdot$_F0 w = u $\cdot$_F0 (v $\cdot$_F0 w).
Fact N'_Mult_Distribution : $\forall$ F0 u v w, Arithmetical_ultraFilter F0
  -> u $\in$ (N' F0) -> v $\in$ (N' F0) -> w $\in$ (N' F0)
  -> u $\cdot$_F0 (v +_F0 w) = (u $\cdot$_F0 v) +_F0 (u $\cdot$_F0 w).
Fact N'_Mult_Cancellation : $\forall$ F0 w u v, Arithmetical_ultraFilter F0
  -> w $\in$ (N' F0) -> u $\in$ (N' F0) -> v $\in$ (N' F0) -> w <> F 0
  -> w $\cdot$_F0 u = w $\cdot$_F0 v -> u = v.
\end{lstlisting}

There are some other straightforward properties related to operations and order on $^{*}{\mathbf N}$:
for any $u,v,w$ in $^{*}{\mathbf N}$, $u <_{^{*}{\mathbf N}} v$ is equivalent to $w +_{^{*}{\mathbf N}} u\ <_{^{*}{\mathbf N}}\ w +_{^{*}{\mathbf N}} v$,
and to $w \cdot_{^{*}{\mathbf N}} u\ <_{^{*}{\mathbf N}}\ w \cdot_{^{*}{\mathbf N}} v$ when $w \neq F_0$ holds;
moreover, if $u <_{^{*}{\mathbf N}} v$ holds, then there exists an unique member $x$ in $^{*}{\mathbf N}$
satisfying $F_0 <_{^{*}{\mathbf N}} x$ and $u +_{^{*}{\mathbf N}} x = v$.

\begin{lstlisting}
Fact N'_Plus_PrOrder : $\forall$ F0 u v w, Arithmetical_ultraFilter F0
  -> u $\in$ (N' F0) -> v $\in$ (N' F0) -> w $\in$ (N' F0)
  -> u <_F0 v <-> (w +_F0 u) <_F0 (w +_F0 v).
Corollary N'_Plus_PrOrder_Corollary : $\forall$ F0 u v, Arithmetical_ultraFilter F0
  -> u $\in$ (N' F0) -> v $\in$ (N' F0)
  -> u <_F0 v <-> $\exists$! w, w $\in$ (N' F0) /\ (F 0) <_F0 w /\ u +_F0 w = v.

Fact N'_Mult_PrOrder : $\forall$ F0 u v w, Arithmetical_ultraFilter F0
  -> u $\in$ (N' F0) -> v $\in$ (N' F0) -> w $\in$ (N' F0) -> w <> F 0
  -> u <_F0 v <-> (w $\cdot$_F0 u) <_F0 (w $\cdot$_F0 v).
\end{lstlisting}

\subsubsection{$\omega$ in $^{*}{\mathbf N}$}\label{sec5-1-4}

As seen in \ref{sec5-1-3} and \ref{sec3-6},
the operations and order on $^{*}{\mathbf N}$ exhibit nearly identical properties to those on $\omega$.
Actually, $\omega$ can be isomorphically embedded in $^{*}{\mathbf N}$,
that is to say, a subset of $^{*}{\mathbf N}$ isomorphic to $\omega$ can be constructed.

\begin{theorem}\label{omega-iso}
For any $^{*}{\mathbf N}$, the following set:
$${^{*}{\mathbf N}}_{\mathbf N} = \{u : u\ is\ an\ principal\ ultrafilter\}$$
is a subset of $^{*}{\mathbf N}$. Meanwhile, a 1-1 function whose domain is $\omega$ and range is ${^{*}{\mathbf N}}_{\mathbf N}$ can be constructed:
$$\varphi = \{(u,v) : u \in \omega \wedge v = F_u\},$$
and for any $m,n$ in $\omega$, $h$ in ${\omega}^{\omega}$, $\varphi$ satisfies:

1) $\varphi(0) = F_0\ \wedge\ \varphi(1) = F_1$.

2) Order Preservation: $m < n \Longleftrightarrow \varphi(m) <_{^{*}{\mathbf N}} \varphi(n)$.

3) Addition Preservation: $\varphi(m + n) = \varphi(m) +_{^{*}{\mathbf N}} \varphi(n)$.

4) Multiplication Preservation: $\varphi(m \cdot n) = \varphi(m) \cdot_{^{*}{\mathbf N}} \varphi(n)$.

5) Value Preservation: $\varphi(h(m)) = h\langle\varphi(m)\rangle$.
\end{theorem}

It requires several steps to formalize this theorem.
Firstly, we construct the subset ${^{*}{\mathbf N}}_{\mathbf N}$.

\begin{lstlisting}
Definition N'_N := \{ $\lambda$ u, $\exists$ n, n $\in$ $\omega$ /\ u = F n \}.
Corollary N'_N_is_Set : Ensemble (N'_N).
Corollary N'_N_Subset_N' : $\forall$ F0, Arithmetical_ultraFilter F0 -> N'_N $\subset$ (N' F0).
\end{lstlisting}

\noindent Next we construct the 1-1 function $\varphi$.

\begin{lstlisting}
Definition $\varphi$ := \{\ $\lambda$ u v, u $\in$ $\omega$ /\ v = F u \}\.
Corollary $\phi$_is_Function : Function1_1 $\phi$ /\ dom($\phi$) = $\omega$ /\ ran($\phi$) = N'_N.
\end{lstlisting}

\noindent Finally, all the five properties of $\varphi$ can be verified.

\begin{lstlisting}
Corollary $\varphi$_0 : $\varphi$[0] = F 0.
Corollary $\varphi$_1 : $\varphi$[1] = F 1.
Corollary $\varphi$_PrOrder : $\forall$ F0 m n, Arithmetical_ultraFilter F0
  -> m $\in$ $\omega$ -> n $\in$ $\omega$ -> m $\in$ n <-> $\varphi$[m] <_F0 $\varphi$[n].
Corollary $\varphi$_PrPlus : $\forall$ F0 m n, Arithmetical_ultraFilter F0
  -> m $\in$ $\omega$ -> n $\in$ $\omega$ -> $\varphi$[m + n] = $\varphi$[m] +_F0 $\varphi$[n].
Corollary $\varphi$_PrMult : $\forall$ F0 m n, Arithmetical_ultraFilter F0
  -> m $\in$ $\omega$ -> n $\in$ $\omega$ -> $\varphi$[m $\cdot$ n] = $\varphi$[m] $\cdot$_F0 $\varphi$[n].
Corollary $\varphi$_PrComposition : $\forall$ F0 m h, Arithmetical_ultraFilter F0
  -> m $\in$ $\omega$ -> Function h -> dom(h) = $\omega$ -> ran(h) $\subset$ $\omega$ -> $\varphi$[h[m]] = h$\langle$$\varphi$[m]$\rangle$.
\end{lstlisting}

Theorem \ref{omega-iso} shows that the structure of $\omega$ and ${^{*}{\mathbf N}}_{\mathbf N}$
are entirely identical under their own respective operations and order.
In other words, ${^{*}{\mathbf N}}_{\mathbf N}$ can also serve as a set of natural numbers.
Consequently, it is evident that the operations on $^{*}{\mathbf N}$ are closed in ${^{*}{\mathbf N}}_{\mathbf N}$ as well.

\begin{lstlisting}
Fact N'_N_Plus_in_N'_N : $\forall$ F0 u v, Arithmetical_ultraFilter F0
  -> u $\in$ N'_N -> v $\in$ N'_N -> (u +_F0 v) $\in$ N'_N.
Fact N'_N_Mult_in_N'_N : $\forall$ F0 u v, Arithmetical_ultraFilter F0
  -> u $\in$ N'_N -> v $\in$ N'_N -> (u $\cdot$_F0 v) $\in$ N'_N.
\end{lstlisting}

\subsubsection{Infinities in $^{*}{\mathbf N}$}\label{sec5-1-5}

Now that the structure ${^{*}{\mathbf N}}_{\mathbf N}$ isomorphic to $\omega$ can be found in $^{*}{\mathbf N}$,
the question arises: what is the difference between ${^{*}{\mathbf N}}_{\mathbf N}$ and $^{*}{\mathbf N}$?
More intuitively, what distinguishes $^{*}{\mathbf N}$ from the natural number set?

Actually, it can be verified that ${^{*}{\mathbf N}}_{\mathbf N}$ is a proper-subset of $^{*}{\mathbf N}$.

\begin{lstlisting}
Fact N'_N_properSubset_N' : $\forall$ F0, Arithmetical_ultraFilter F0
  -> ($\forall$ n, F0 <> F n) -> N'_N $\subset$ (N' F0) /\ N'_N <> (N' F0).
\end{lstlisting}

\noindent With this consequence, it is feasible to discuss the members in ${^{*}{\mathbf N}} \sim {^{*}{\mathbf N}}_{\mathbf N}$.

\begin{theorem}
For every member $\tau$ in (${^{*}{\mathbf N}} \sim {^{*}{\mathbf N}}_{\mathbf N}$) and every member $u$ in ${^{*}{\mathbf N}}_{\mathbf N}$,
$u <_{^{*}{\mathbf N}} \tau$ holds.
\end{theorem}

\begin{lstlisting}
Theorem FT12 : $\forall$ F0 t u, Arithmetical_ultraFilter F0
  -> ($\forall$ n, F0 <> F n) -> t $\in$ ((N' F0) $\sim$ N'_N) -> u $\in$ N'_N -> u <_F0 t.
\end{lstlisting}

Now we can observe that the members of ${^{*}{\mathbf N}} \sim {^{*}{\mathbf N}}_{\mathbf N}$ are all infinity numbers
because each of them is larger than every member of ${^{*}{\mathbf N}}_{\mathbf N}$ within the order on $^{*}{\mathbf N}$.
Intuitively, the numbers in ${^{*}{\mathbf N}} \sim {^{*}{\mathbf N}}_{\mathbf N}$ are larger than all natural numbers
but share the same operation and order properties as natural numbers.
Therefore, ${^{*}{\mathbf N}}$ is actually a non-standard natural number set that includes infinity natural numbers,
and its structure can be intuitively depicted as:
$$\underbrace{0 < 1 < 2 < \cdots < n <}_{{^{*}{\mathbf N}}_{\mathbf N}} \cdots < \tau < (\tau + 1) < \cdots$$
where $\tau$ represents an infinity natural number.

\subsection{Extension from $^{*}{\mathbf N}$ to $^{*}{\mathbf Z}$}\label{sec5-2}

\subsubsection{Equivalence Relation and Classification}\label{sec5-2-1}

FMCR adopts the idea of equivalence classification to extend number systems
and the relevant definitions and theorems should be formalized first.
According to MK system, the relevant concepts can be defined on general classes.

\begin{definition}[Equivalence Relation]\label{equ-rela}
For each class $a$, let $R$ be a subclass of $a \times a$, $R$ is called an equivalence relation over $a$ if it satisfies:

1) Reflexivity: $\forall x \in a, xRx$.

2) Symmetry: $\forall x,y \in a$, if $xRy$ then $yRx$.

3) Transitivity: $\forall x,y,z \in a$, if $xRy$ and $yRz$ then $xRz$.

\noindent where $xRy$ stands for $(x,y) \in R$, namely $(x,y)$ satisfies $R$-relation.
\end{definition}

\begin{lstlisting}
Definition Reflexivity R a := $\forall$ x, x $\in$ a -> [x,x] $\in$ R.
Definition Symmetry R a := $\forall$ x y, x $\in$ a -> y $\in$ a
  -> [x,y] $\in$ R -> [y,x] $\in$ R.
Definition Transitivity R a := $\forall$ x y z, x $\in$ a -> y $\in$ a -> z $\in$ a
  -> [x,y] $\in$ R -> [y,z] $\in$ R -> [x,z] $\in$ R.

Definition equRelation R a := R $\subset$ a $\times$ a
  /\ Reflexivity R a /\ Symmetry R a /\ Transitivity R a.
\end{lstlisting}

A trivial example of equivalence relation is identity relation.
Because for every $a$ and its members $x$ and $y$, it obviously satisfies
$x = x$ (reflexivity), if $x = y$ then $y = x$ (symmetry), and if $x = y$ and $y = z$ then $x = z$ (transitivity).
The identity relation over $a$ can be set as $\{(u,v) : u \in a \wedge v \in a \wedge u = v\}$.

\begin{lstlisting}
Example Example_equR : $\forall$ a, equRelation (\{\ $\lambda$ u v, u $\in$ a
  /\ v $\in$ a /\ u = v \}\) a.
\end{lstlisting}

\begin{definition}[Equivalence Class]\label{equ-class}
Let $R$ be an equivalence relation over $a$, for each $x$ in $a$, the following class
$$\{u : u \in a \wedge uRx \}$$
is called equivalence class represented by $x$ regarding the $R$-relation, denoted as $[x]_R$.
\end{definition}

To formalize this definition, three parameters are required:
the first one serves as the representation $x$ of this equivalence class, the second the equivalence relation $R$ and the last the class $a$.
And for formalization, a notation must contain all three parameters and cannot be simplified like the mathematical symbol $[x]_R$
where some parameter can even be omitted.
Thus we do not introduce notation for equivalence class since no convenience will arise.

\begin{lstlisting}
Definition equClass x R a := \{ $\lambda$ u, u $\in$ a /\ [u,x] $\in$ R \}.
\end{lstlisting}

Next we introduce two theorems concerning equivalence class.

\begin{theorem}
Let $R$ be an equivalence relation over $a$, for any $x,y$ in $a$, $xRy$ if and only if $[x]_R = [y]_R$.
\end{theorem}

\begin{theorem}
Let $R$ be an equivalence relation over $a$, for any $x,y$ in $a$, $[x]_R \neq [y]_R$ if and only if $[x]_R \cap [y]_R = \emptyset$.
\end{theorem}

\begin{lstlisting}
Theorem equClassT1 : $\forall$ R a x y, equRelation R a -> x $\in$ a -> y $\in$ a
  -> [x,y] $\in$ R <-> equClass x R a = equClass y R a.
Theorem equClassT2 : $\forall$ R a x y, equRelation R a -> x $\in$ a -> y $\in$ a
  -> equClass x R a <> equClass y R a
    <-> (equClass x R a) $\cap$ (equClass y R a) = $\Phi$.
\end{lstlisting}

The two theorems indicate that for any two equivalence classes $u$ and $v$ of $a$, there exist only two situations:
$u$ and $v$ are entirely identical or $u$ and $v$ have no common elements (entirely not identical).
So this is actually a procedure of classifying the elements of $a$ according to some specific equivalence relation.
The next definition describes the result of this classification.

\begin{definition}[Quotient Class]\label{quo-class}
Let $R$ be an equivalence relation over $a$, the class consisting of all equivalence classes regarding $R$ is denoted as:
$$a / R = \{u : \exists x, x \in a \wedge u = [x]_R\},$$
which is called a quotient class (or set if proven) of $a$.
\end{definition}

\begin{lstlisting}
Definition QuotientClass a R := \{ $\lambda$ u, $\exists$ x, x $\in$ a /\ u = equClass x R a \}.
Notation "a / R" := (QuotientClass a R) : ec_scope.
\end{lstlisting}

\subsubsection{Construction of $^{*}{\mathbf Z}$}\label{sec5-2-2}

The process of extending $^{*}{\mathbf N}$ to $^{*}{\mathbf Z}$ is actually
equivalently classifying the elements of ${^{*}{\mathbf N}} \times {^{*}{\mathbf N}}$,
and a specific equivalence relation over ${^{*}{\mathbf N}} \times {^{*}{\mathbf N}}$ must be set first:
$$R_{^{*}{\mathbf N}} = \{((m,n),(p,q)) : m,n,p,q \in {^{*}{\mathbf N}} \wedge m +_{^{*}{\mathbf N}} q = n +_{^{*}{\mathbf N}} p\}.$$

\begin{lstlisting}
Definition R_N' F0 := \{\ $\lambda$ u v, $\exists$ m n p q, m $\in$ (N' F0)
  /\ n $\in$ (N' F0) /\ p $\in$ (N' F0) /\ q $\in$ (N' F0) /\ u = [m,n]
  /\ v = [p,q] /\ m +_F0 q = n +_F0 p \}\.
Corollary R_N'_is_equRelation : $\forall$ F0, Arithmetical_ultraFilter F0
  -> equRelation (R_N' F0) (N' F0)$\times$(N' F0).
\end{lstlisting}

Now that $R_{^{*}{\mathbf N}}$ is an equivalence relation over ${^{*}{\mathbf N}} \times {^{*}{\mathbf N}}$,
it is logical to discuss relationship between two identical equivalence classes regarding $R_{^{*}{\mathbf N}}$:
$$\forall m,n,p,q \in {^{*}{\mathbf N}},\ [(m,n)]_{R_{^{*}{\mathbf N}}} = [(p,q)]_{R_{^{*}{\mathbf N}}}\ \Longleftrightarrow\
m\ +_{^{*}{\mathbf N}}\ q\ =\ n\ +_{^{*}{\mathbf N}}\ p.$$

\begin{lstlisting}
Notation "\[ u \]_ F0" := (equClass u (R_N' F0) (N' F0)$\times$(N' F0))
  (at level 5, u at level 0, F0 at level 0) : z'_scope.

Corollary R_N'_Corollary : $\forall$ F0 m n p q, Arithmetical_ultraFilter F0
  -> m $\in$ (N' F0) -> n $\in$ (N' F0) -> p $\in$ (N' F0) -> q $\in$ (N' F0)
  -> \[[m,n]\]_F0 = \[[p,q]\]_F0 <-> m +_F0 q = n +_F0 p.
\end{lstlisting}

To simplify formalization code, here we introduce the notation ``{\verb+\[ u \]_ F0+" for the equivalence class regarding $R_{^{*}{\mathbf N}}$.
Different from the formalization of definition \ref{equ-class} where no notations of equivalence class are introduced,
here the last two parameters of {\tt equClass} are instantiated, and it only requires two variables for its notation,
which can reduce the code amount to some extent.

Now let us explore from an intuitive perspective how $R_{^{*}{\mathbf N}}$ classifies ${^{*}{\mathbf N}} \times {^{*}{\mathbf N}}$.
According to the equivalence relation, we can immediately find many equivalent elements in ${^{*}{\mathbf N}} \times {^{*}{\mathbf N}}$.
For example, $(0,0),(1,1),(2,2),\cdots,(\tau,\tau),\cdots$ are equivalent; $(0,1),(1,2),(2,3),\cdots,(\tau,\tau + 1),\cdots$ are equivalent; and
so are $(2,0),(3,1),(4,2),\cdots,(\tau + 2,\tau),\cdots$, where $\tau$ represents an infinity natural number.
In a planar dot array composed of elements in ${^{*}{\mathbf N}} \times {^{*}{\mathbf N}}$,
we can connect equivalent points with straight lines and thus obtain figure \ref{fig2}.

\begin{figure}[h]%
\centering
\includegraphics[width=1\textwidth]{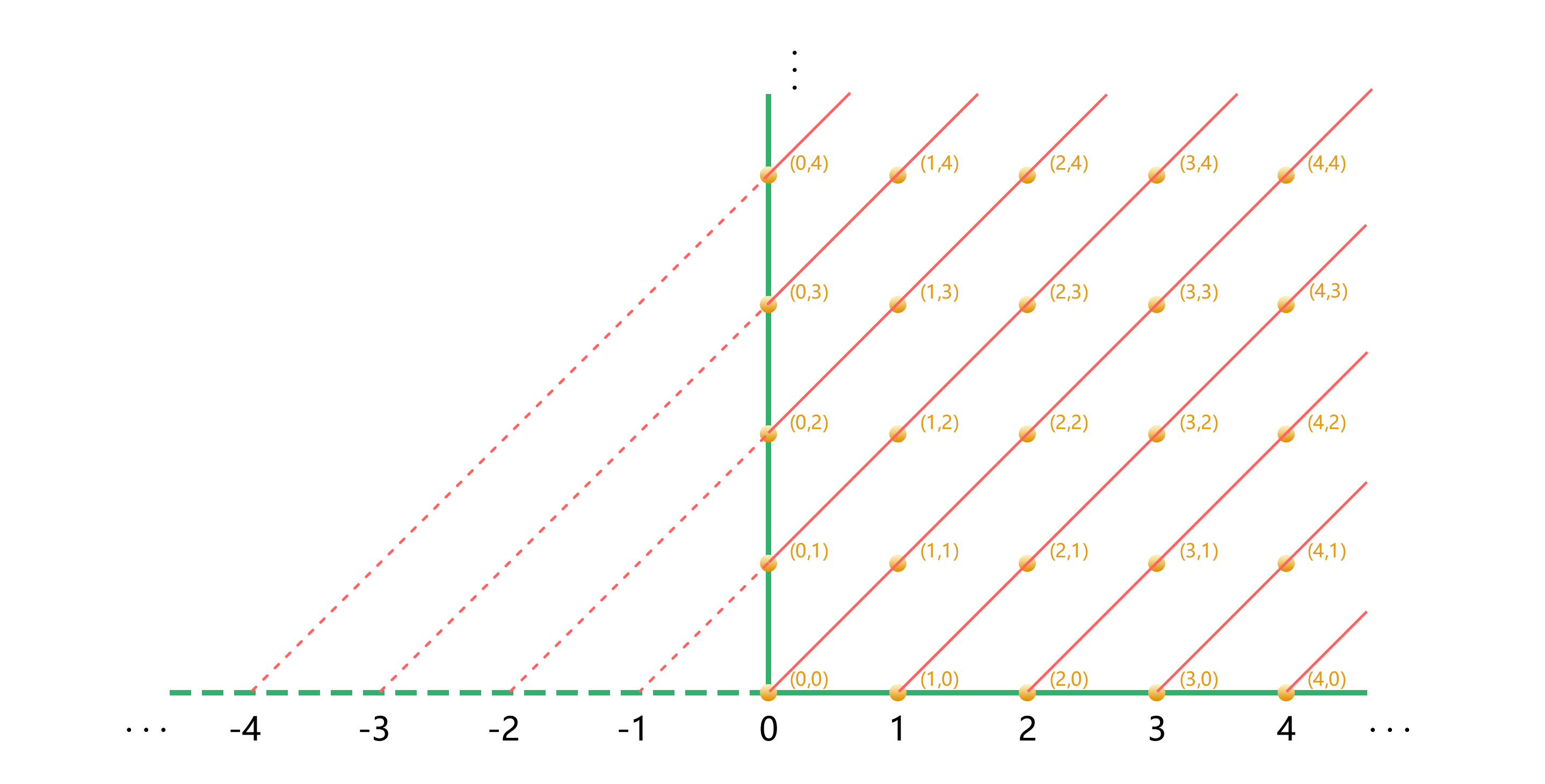}
\caption{equivalence classification of elements in ${^{*}{\mathbf N}} \times {^{*}{\mathbf N}}$}\label{fig2}
\end{figure}

In figure \ref{fig2}, every oblique line represents an equivalence class of ${^{*}{\mathbf N}} \times {^{*}{\mathbf N}}$,
and intuitively the order of these oblique lines on the horizontal axis forms a new number ``sequence"\footnote
{This ``sequence" is extended to infinity numbers.}
and the lines extended by dashed lines extends this ``sequence" to range of negative numbers.
Then the set consisting of all these oblique lines,
namely the quotient set of ${^{*}{\mathbf N}} \times {^{*}{\mathbf N}}$ regarding $R_{^{*}{\mathbf N}}$,
is exactly ${^{*}{\mathbf Z}}$.
$${^{*}{\mathbf Z}} = ({^{*}{\mathbf N}} \times {^{*}{\mathbf N}}) / {R_{^{*}{\mathbf N}}} = \{u : \exists x, u = [x]_{R_{^{*}{\mathbf N}}}\}.$$

\begin{lstlisting}
Definition Z' F0 := (N' F0)$\times$(N' F0) / (R_N' F0).
Corollary Z'_is_Set : $\forall$ F0, F0 $\in$ $\beta\omega$ -> Ensemble (Z' F0).

Corollary Z'_Corollary1 : $\forall$ F0 u, Arithmetical_ultraFilter F0
  -> u $\in$ (Z' F0) -> $\exists$ m n, m $\in$ (N' F0) /\ n $\in$ (N' F0) /\ u = \[[m,n]\]_F0.
Corollary Z'_Corollary2 : $\forall$ F0 m n, Arithmetical_ultraFilter F0
  -> m $\in$ (N' F0) -> n $\in$ (N' F0) -> \[[m,n]\]_F0 $\in$ (Z' F0).
\end{lstlisting}

Like ${^{*}{\mathbf N}}$, ${^{*}{\mathbf Z}}$ has its own zero and identity element,
which can be easily found according to figure \ref{fig2}.
We give the symbol $^{*}$Z$_0$ and $^{*}$Z$_1$ respectively to them:
$${^{*}}\mbox{Z}_0 = \{(u,v) : u,v \in {^{*}{\mathbf N}} \wedge u = v\} = [(F_0,F_0)]_{R_{^{*}{\mathbf N}}},$$
$${^{*}}\mbox{Z}_1 = \{(u,v) : u,v \in {^{*}{\mathbf N}} \wedge u = v +_{^{*}{\mathbf N}} F_1\} = [(F_1,F_0)]_{R_{^{*}{\mathbf N}}},$$
where $F_0$ is the principal ultrafilter corresponding to $0$ and $F_1$ to $1$,
respectively standing for zero and identity element in ${^{*}{\mathbf N}}$.

\begin{lstlisting}
Definition Z'0 F0 := \{\ $\lambda$ u v, u $\in$ (N' F0) /\ v $\in$ (N' F0) /\ u = v \}\.
Corollary Z'0_Corollary : $\forall$ F0, Arithmetical_ultraFilter F0
  -> Z'0 F0 = \[[(F 0),(F 0)]\]_F0.

Definition Z'1 F0 := \{\ $\lambda$ u v, v $\in$ (N' F0) /\ u = v +_F0 (F 1) \}\.
Corollary Z'1_Corollary : $\forall$ F0, Arithmetical_ultraFilter F0
  -> Z'1 F0 = \[[(F 1),(F 0)]\]_F0.
\end{lstlisting}

\noindent And of course $^{*}$Z$_0$ and $^{*}$Z$_1$ are two different members in ${^{*}{\mathbf Z}}$.

\begin{lstlisting}
Corollary Z'0_in_Z' : $\forall$ F0, Arithmetical_ultraFilter F0 -> (Z'0 F0) $\in$ (Z' F0).
Corollary Z'1_in_Z' : $\forall$ F0, Arithmetical_ultraFilter F0 -> (Z'1 F0) $\in$ (Z' F0).

Fact Z'0_isn't_Z'1 : $\forall$ F0, Arithmetical_ultraFilter F0 -> Z'0 F0 <> Z'1 F0.
\end{lstlisting}

\subsubsection{Order on $^{*}{\mathbf Z}$}\label{sec5-2-3}

\begin{definition}[Order on $^{*}{\mathbf Z}$]
Let $u,v$ be members in $^{*}{\mathbf Z}$ and $u = [(m,n)]_{R_{^{*}{\mathbf N}}}, v = [(p,q)]_{R_{^{*}{\mathbf N}}}$
($m,n,p,q \in {^{*}{\mathbf N}}$), then the order of $u$ and $v$ is defined as follows:
$$u <_{^{*}{\mathbf Z}} v\ \Longleftrightarrow\
[(m,n)]_{R_{^{*}{\mathbf N}}} <_{^{*}{\mathbf Z}} [(p,q)]_{R_{^{*}{\mathbf N}}}\ \Longleftrightarrow\
m +_{^{*}{\mathbf N}} q\ <_{^{*}{\mathbf N}}\ n +_{^{*}{\mathbf N}} p.$$
\end{definition}

The definition of this order relation is reasonable,
because for any $m,m_1,n,n_1,p,p_1,q,q_1 \in {^{*}{\mathbf N}}$, if
$$u = [(m,n)]_{R_{^{*}{\mathbf N}}} = [(m_1,n_1)]_{R_{^{*}{\mathbf N}}},\
v = [(p,q)]_{R_{^{*}{\mathbf N}}} = [(p_1,q_1)]_{R_{^{*}{\mathbf N}}},$$
then it can be verified that
$$u <_{^{*}{\mathbf Z}} v\ \Longleftrightarrow\
m +_{^{*}{\mathbf N}} q\ <_{^{*}{\mathbf N}}\ n +_{^{*}{\mathbf N}} p\ \Longleftrightarrow\
m_1 +_{^{*}{\mathbf N}} q_1\ <_{^{*}{\mathbf N}}\ n_1 +_{^{*}{\mathbf N}} p_1.$$
Namely, the order relation on $^{*}{\mathbf Z}$ is independent of the selection of representations of equivalence classes.

\begin{lstlisting}
Definition Z'_Ord F0 := \{\ $\lambda$ u v, $\forall$ m n p q, m $\in$ (N' F0) -> n $\in$ (N' F0)
  -> p $\in$ (N' F0) -> q $\in$ (N' F0) -> u = \[[m,n]\]_F0 -> v = \[[p,q]\]_F0
  -> (m +_F0 q) <_F0 (n +_F0 p) \}\.
Notation "u <_ F0 v" := ([u,v] $\in$ (Z'_Ord F0))
  (at level 10, F0 at level 9) : z'_scope.

Corollary Z'_Ord_reasonable :
  $\forall$ F0 m n p q m1 n1 p1 q1, Arithmetical_ultraFilter F0
  -> m $\in$ (N' F0) -> n $\in$ (N' F0) -> p $\in$ (N' F0) -> q $\in$ (N' F0)
  -> m1 $\in$ (N' F0) -> n1 $\in$ (N' F0) -> p1 $\in$ (N' F0) -> q1 $\in$ (N' F0)
  -> \[[m,n]\]_F0 = \[[m1,n1]\]_F0 -> \[[p,q]\]_F0 = \[[p1,q1]\]_F0
  -> ((m +_F0 q) <_F0 (n +_F0 p) <-> (m1 +_F0 q1) <_F0 (n1 +_F0 p1))$\mbox{\color[RGB]{0,100,0}\%}$n'.
Corollary Z'_Ord_Corollary : $\forall$ F0 m n p q, Arithmetical_ultraFilter F0
  -> m $\in$ (N' F0) -> n $\in$ (N' F0) -> p $\in$ (N' F0) -> q $\in$ (N' F0)
  -> \[[m,n]\]_F0 <_F0 \[[p,q]\]_F0 <-> ((m +_F0 q) <_F0 (n +_F0 p))$\mbox{\color[RGB]{0,100,0}\%}$n'.
\end{lstlisting}

\noindent The delimiting key {\tt n'} indicates that the expressions it explains are connected with $^{*}{\mathbf N}$.

As for the properties of this order, they are same as the order on $^{*}{\mathbf N}$,
namely satisfying irreflexivity, transitivity and trichotomy, and $^{*}$Z$_0$ $<_{^{*}{\mathbf Z}}$ $^{*}$Z$_1$ undoubtedly holds as well.

\begin{lstlisting}
Fact Z'_Ord_irReflex : $\forall$ F0 u v, Arithmetical_ultraFilter F0
  -> u $\in$ (Z' F0) -> v $\in$ (Z' F0) -> u = v -> $\sim$ u <_F0 v.
Fact Z'_Ord_Trans : $\forall$ F0 u v w, Arithmetical_ultraFilter F0
  -> u $\in$ (Z' F0) -> v $\in$ (Z' F0) -> w $\in$ (Z' F0) -> u <_F0 v
  -> v <_F0 w -> u <_F0 w.
Fact Z'_Ord_Tri : $\forall$ F0, Arithmetical_ultraFilter F0
  -> Connect (Z'_Ord F0) (Z' F0).
Fact Z'0_lt_Z'1 : $\forall$ F0, Arithmetical_ultraFilter F0 -> (Z'0 F0) <_F0 (Z'1 F0).
\end{lstlisting}

\subsubsection{Operations on $^{*}{\mathbf Z}$}\label{sec5-2-4}

\begin{definition}[Addition on $^{*}{\mathbf Z}$]
Let $u,v$ be members in $^{*}{\mathbf Z}$ and $u = [(m,n)]_{R_{^{*}{\mathbf N}}}, v = [(p,q)]_{R_{^{*}{\mathbf N}}}$
($m,n,p,q \in {^{*}{\mathbf N}}$), then the sum of $u$ and $v$ is defined as follows:
$$u +_{^{*}{\mathbf Z}} v\ =\ [(m,n)]_{R_{^{*}{\mathbf N}}}\ +_{^{*}{\mathbf Z}}\ [(p,q)]_{R_{^{*}{\mathbf N}}}\
=\ [(m +_{^{*}{\mathbf N}} p\ ,\ n +_{^{*}{\mathbf N}} q)]_{R_{^{*}{\mathbf N}}}.$$
\end{definition}

The definition of this addition is reasonable, because for any of $m,m_1,n,n_1,p,p_1,q,q_1$ in ${^{*}{\mathbf N}}$, if
$$u = [(m,n)]_{R_{^{*}{\mathbf N}}} = [(m_1,n_1)]_{R_{^{*}{\mathbf N}}},\
v = [(p,q)]_{R_{^{*}{\mathbf N}}} = [(p_1,q_1)]_{R_{^{*}{\mathbf N}}},$$
then it can be verified that
$$u +_{^{*}{\mathbf Z}} v\ =\ [(m,n)]_{R_{^{*}{\mathbf N}}}\ +_{^{*}{\mathbf Z}}\ [(p,q)]_{R_{^{*}{\mathbf N}}}\ =\
[(m_1,n_1)]_{R_{^{*}{\mathbf N}}}\ +_{^{*}{\mathbf Z}}\ [(p_1,q_1)]_{R_{^{*}{\mathbf N}}}.$$
Namely, the result of $u +_{^{*}{\mathbf Z}} v$ is unique and independent of the selection of representations of equivalence classes.

\begin{lstlisting}
Definition Z'_Plus F0 u v := $\cap$(\{ $\lambda$ a, $\forall$ m n p q, m $\in$ (N' F0)
  -> n $\in$ (N' F0) -> p $\in$ (N' F0) -> q $\in$ (N' F0)
  -> u = \[[m,n]\]_F0 -> v = \[[p,q]\]_F0
  -> a = \[[(m +_F0 p)$\mbox{\color[RGB]{0,100,0}\%}$n',(n +_F0 q)$\mbox{\color[RGB]{0,100,0}\%}$n']\]_F0 \}).
Notation "u +_ F0 v" := (Z'_Plus F0 u v)
  (at level 10, F0 at level 9) : z'_scope.

Corollary Z'_Plus_reasonable : $\forall$ F0 u v, Arithmetical_ultraFilter F0
  -> u $\in$ (Z' F0) -> v $\in$ (Z' F0) -> $\exists$ a, a $\in$ (Z' F0) /\ [a]
    = \{ $\lambda$ a, $\forall$ m n p q, m $\in$ (N' F0)
    -> n $\in$ (N' F0) -> p $\in$ (N' F0) -> q $\in$ (N' F0)
    -> u = \[[m,n]\]_F0 -> v = \[[p,q]\]_F0
    -> a = \[[(m +_F0 p)$\mbox{\color[RGB]{0,100,0}\%}$n',(n +_F0 q)$\mbox{\color[RGB]{0,100,0}\%}$n']\]_F0 \}.
Corollary Z'_Plus_Corollary : $\forall$ F0 m n p q, Arithmetical_ultraFilter F0
  -> m $\in$ (N' F0) -> n $\in$ (N' F0) -> p $\in$ (N' F0) -> q $\in$ (N' F0)
  -> \[[m,n]\]_F0 +_F0 \[[p,q]\]_F0 = \[[(m +_F0 p)$\mbox{\color[RGB]{0,100,0}\%}$n',(n +_F0 q)$\mbox{\color[RGB]{0,100,0}\%}$n']\]_F0.
\end{lstlisting}

\begin{definition}[Multiplication on $^{*}{\mathbf Z}$]
Let $u,v$ be members in $^{*}{\mathbf Z}$ and $u = [(m,n)]_{R_{^{*}{\mathbf N}}}, v = [(p,q)]_{R_{^{*}{\mathbf N}}}$
($m,n,p,q \in {^{*}{\mathbf N}}$), then the product of $u$ and $v$ is defined as follows:
\begin{align}
u \cdot_{^{*}{\mathbf Z}} v\ &=\ [(m,n)]_{R_{^{*}{\mathbf N}}}\ \cdot_{^{*}{\mathbf Z}}\ [(p,q)]_{R_{^{*}{\mathbf N}}}\notag\\
&=\ [((m \cdot_{^{*}{\mathbf N}} p)\ +_{^{*}{\mathbf N}}\ (n \cdot_{^{*}{\mathbf N}} q)\ ,\
(m \cdot_{^{*}{\mathbf N}} q)\ +_{^{*}{\mathbf N}}\ (n \cdot_{^{*}{\mathbf N}} p))]_{R_{^{*}{\mathbf N}}}.\notag
\end{align}
\end{definition}

The definition of this multiplication is reasonable, because for any $m,m_1,n,n_1,p,p_1,q,q_1 \in {^{*}{\mathbf N}}$, if
$$u = [(m,n)]_{R_{^{*}{\mathbf N}}} = [(m_1,n_1)]_{R_{^{*}{\mathbf N}}},\
v = [(p,q)]_{R_{^{*}{\mathbf N}}} = [(p_1,q_1)]_{R_{^{*}{\mathbf N}}},$$
then it can be verified that
$$u \cdot_{^{*}{\mathbf Z}} v\ =\ [(m,n)]_{R_{^{*}{\mathbf N}}}\ \cdot_{^{*}{\mathbf Z}}\ [(p,q)]_{R_{^{*}{\mathbf N}}}\ =\
[(m_1,n_1)]_{R_{^{*}{\mathbf N}}}\ \cdot_{^{*}{\mathbf Z}}\ [(p_1,q_1)]_{R_{^{*}{\mathbf N}}}.$$
Namely, the result of $u \cdot_{^{*}{\mathbf Z}} v$ is unique and independent of the selection of representations of equivalence classes.

\begin{lstlisting}
Definition Z'_Mult F0 u v := $\cap$(\{ $\lambda$ a, $\forall$ m n p q, m $\in$ (N' F0)
  -> n $\in$ (N' F0) -> p $\in$ (N' F0) -> q $\in$ (N' F0)
  -> u = \[[m,n]\]_F0 -> v = \[[p,q]\]_F0
  -> a = \[[((m $\cdot$_F0 p) +_F0 (n $\cdot$_F0 q))$\mbox{\color[RGB]{0,100,0}\%}$n',
    ((m $\cdot$_F0 q) +_F0 (n $\cdot$_F0 p))$\mbox{\color[RGB]{0,100,0}\%}$n']\]_F0 \}).
Notation "u $\cdot$_ F0 v" := (Z'_Mult F0 u v)
  (at level 10, F0 at level 9) : z'_scope.

Corollary Z'_Mult_reasonable : $\forall$ F0 u v, Arithmetical_ultraFilter F0
  -> u $\in$ (Z' F0) -> v $\in$ (Z' F0) -> $\exists$ a, a $\in$ (Z' F0) /\ [a]
    = \{ $\lambda$ a, $\forall$ m n p q, m $\in$ (N' F0) -> n $\in$ (N' F0) -> p $\in$ (N' F0)
    -> q $\in$ (N' F0) -> u = \[[m,n]\]_F0 -> v = \[[p,q]\]_F0
    -> a = \[[((m $\cdot$_F0 p) +_F0 (n $\cdot$_F0 q))$\mbox{\color[RGB]{0,100,0}\%}$n',
      ((m $\cdot$_F0 q) +_F0 (n $\cdot$_F0 p))$\mbox{\color[RGB]{0,100,0}\%}$n']\]_F0 \}.
Corollary Z'_Mult_Corollary : $\forall$ F0 m n p q, Arithmetical_ultraFilter F0
  -> m $\in$ (N' F0) -> n $\in$ (N' F0) -> p $\in$ (N' F0) -> q $\in$ (N' F0)
  -> \[[m,n]\]_F0 $\cdot$_F0 \[[p,q]\]_F0
    = \[[((m $\cdot$_F0 p) +_F0 (n $\cdot$_F0 q))$\mbox{\color[RGB]{0,100,0}\%}$n',((m $\cdot$_F0 q) +_F0 (n $\cdot$_F0 p))$\mbox{\color[RGB]{0,100,0}\%}$n']\]_F0.
\end{lstlisting}

The operations on $^{*}{\mathbf Z}$ possess all the properties that operations on $^{*}{\mathbf N}$ have
and the relevant formalization codes are extremely similar to those on $^{*}{\mathbf N}$, thus we will not present them again.
Nevertheless, the addition on $^{*}{\mathbf Z}$ possesses one more property compared to that on $^{*}{\mathbf N}$:
for every $u$ in $^{*}{\mathbf Z}$,
there exists a unique element $u_0$ of $^{*}{\mathbf Z}$ satisfying $u +_{^{*}{\mathbf Z}} u_0 = {^{*}}$Z$_0$.

\begin{lstlisting}
Fact Z'_neg : $\forall$ F0 u, Arithmetical_ultraFilter F0 -> u $\in$ (Z' F0)
  -> ($\exists$! v, v $\in$ (Z' F0) /\ (u +_F0 v) = (Z'0 F0)).
\end{lstlisting}

This property illuminates that every member of $^{*}{\mathbf Z}$ has its unique corresponding negative element in $^{*}{\mathbf Z}$.
And according to the structure of $^{*}{\mathbf Z}$, here comes the consequence that
if $[(m,n)]_{R_{^{*}{\mathbf N}}} +_{^{*}{\mathbf Z}} u = {^{*}}$Z$_0$, then $u = [(n,m)]_{R_{^{*}{\mathbf N}}}$,
where $m,n$ are members in $^{*}{\mathbf N}$, $u$ in $^{*}{\mathbf Z}$.

\begin{lstlisting}
Fact Z'_neg_Corollary : $\forall$ F0 m n u, Arithmetical_ultraFilter F0
  -> m $\in$ (N' F0) -> n $\in$ (N' F0) -> u $\in$ (Z' F0)
  -> \[[m,n]\]_F0 +_F0 u = Z'0 F0 -> u = \[[n,m]\]_F0.
\end{lstlisting}

\noindent This consequence reveals the structural connection between positive and negative elements in $^{*}{\mathbf Z}$.

\subsubsection{$^{*}{\mathbf N}$ in $^{*}{\mathbf Z}$}\label{sec5-2-5}

Just as $\omega$ can be isomorphically embedded in $^{*}{\mathbf N}$,
$^{*}{\mathbf N}$ can also be isomorphically embedded in $^{*}{\mathbf Z}$.

\begin{theorem}\label{n'-in-z'}
Let $^{*}{\mathbf Z}$ be constructed with $^{*}{\mathbf N}$, the following set:
$${^{*}{\mathbf Z}}_{^{*}{\mathbf N}} = \{u : \exists m, m \in {^{*}{\mathbf N}} \wedge u = [(m,F_0)]_{R_{^{*}{\mathbf N}}}\}$$
is a proper-subset of $^{*}{\mathbf Z}$.
Meanwhile, a 1-1 function whose domain is $^{*}{\mathbf N}$ and range is ${^{*}{\mathbf Z}}_{^{*}{\mathbf N}}$ can be constructed:
$$\varphi_1 = \{(u,v) : u \in {^{*}{\mathbf N}} \wedge v = [(u,F_0)]_{R_{^{*}{\mathbf N}}}\},$$
and for any $m,n$ in ${^{*}{\mathbf N}}$, $\varphi_1$ satisfies:

1) $\varphi_1(F_0) = {^{*}}$Z$_0\ \wedge\ \varphi_1(F_1) = {^{*}}$Z$_1$.

2) Order Preservation: $m <_{^{*}{\mathbf N}} n \Longleftrightarrow \varphi_1(m) <_{^{*}{\mathbf Z}} \varphi_1(n)$.

3) Addition Preservation: $\varphi_1(m +_{^{*}{\mathbf N}} n) = \varphi_1(m) +_{^{*}{\mathbf Z}} \varphi_1(n)$.

4) Multiplication Preservation: $\varphi_1(m \cdot_{^{*}{\mathbf N}} n) = \varphi_1(m) \cdot_{^{*}{\mathbf Z}} \varphi_1(n)$.
\end{theorem}

\begin{lstlisting}
Definition Z'_N' F0 := \{ $\lambda$ u, $\exists$ m, m $\in$ (N' F0) /\ u = \[[m,(F 0)]\]_F0 \}.
Corollary Z'_N'_is_Set : $\forall$ F0, F0 $\in$ $\beta\omega$ -> Ensemble (Z'_N' F0).
Corollary Z'_N'_properSubset_Z' : $\forall$ F0, Arithmetical_ultraFilter F0
  -> (Z'_N' F0) $\subset$ (Z' F0) /\ (Z'_N' F0) <> (Z' F0).

Definition $\varphi$1 F0 := \{\ $\lambda$ u v, u $\in$ (N' F0) /\ v = \[[u,(F 0)]\]_F0 \}\.
Corollary $\varphi$1_is_Function : $\forall$ F0, Arithmetical_ultraFilter F0
  -> Function1_1 ($\varphi$1 F0) /\ dom($\varphi$1 F0) = (N' F0)
    /\ ran($\varphi$1 F0) = (Z'_N' F0).

Corollary $\varphi$1_N'0 : $\forall$ F0, Arithmetical_ultraFilter F0 -> ($\varphi$1 F0)[F 0] = Z'0 F0.
Corollary $\varphi$1_N'1 : $\forall$ F0, Arithmetical_ultraFilter F0 -> ($\varphi$1 F0)[F 1] = Z'1 F0.
Corollary $\varphi$1_PrOrder : $\forall$ F0 m n, Arithmetical_ultraFilter F0 -> m $\in$ (N' F0)
  -> n $\in$ (N' F0) -> (m <_F0 n)$\mbox{\color[RGB]{0,100,0}\%}$n' <-> ($\varphi$1 F0)[m] <_F0 ($\varphi$1 F0)[n].
Corollary $\varphi$1_PrPlus : $\forall$ F0 m n, Arithmetical_ultraFilter F0 -> m $\in$ (N' F0)
  -> n $\in$ (N' F0) -> ($\varphi$1 F0)[(m +_F0 n)$\mbox{\color[RGB]{0,100,0}\%}$n'] = ($\varphi$1 F0)[m] +_F0 ($\varphi$1 F0)[n].
Corollary $\varphi$1_PrMult : $\forall$ F0 m n, Arithmetical_ultraFilter F0 -> m $\in$ (N' F0)
  -> n $\in$ (N' F0) -> ($\varphi$1 F0)[(m $\cdot$_F0 n)$\mbox{\color[RGB]{0,100,0}\%}$n'] = ($\varphi$1 F0)[m] $\cdot$_F0 ($\varphi$1 F0)[n].
\end{lstlisting}

Concerning the set ${^{*}{\mathbf Z}}_{^{*}{\mathbf N}}$, the following four properties hold:

1) ${^{*}{\mathbf Z}} \sim {^{*}{\mathbf Z}}_{^{*}{\mathbf N}} = \{u : \exists m, m \neq F_0 \wedge m \in {^{*}{\mathbf N}} \wedge
u = [(F_0,m)]_{R_{^{*}{\mathbf N}}}\}.$

2) $\forall u \in {^{*}{\mathbf Z}},\ u \in ({^{*}{\mathbf Z}} \sim {^{*}{\mathbf Z}}_{^{*}{\mathbf N}})\ \Longleftrightarrow\ u <_{^{*}{\mathbf Z}} {^{*}}$Z$_0$.

3) $\forall u \in ({^{*}{\mathbf Z}}_{^{*}{\mathbf N}} \sim \{{^{*}}$Z$_0\}),\
\exists! u_0 \in ({^{*}{\mathbf Z}} \sim {^{*}{\mathbf Z}}_{^{*}{\mathbf N}}) \wedge u +_{^{*}{\mathbf N}} u_0 = {^{*}}$Z$_0$.

4) $\forall u \in ({^{*}{\mathbf Z}} \sim {^{*}{\mathbf Z}}_{^{*}{\mathbf N}}),\
\exists! u_0 \in ({^{*}{\mathbf Z}}_{^{*}{\mathbf N}} \sim \{{^{*}}$Z$_0\}) \wedge u +_{^{*}{\mathbf N}} u_0 = {^{*}}$Z$_0$.

\begin{lstlisting}
Fact Z'_N'_Property1 : $\forall$ F0, Arithmetical_ultraFilter F0
  -> (Z' F0) $\sim$ (Z'_N' F0) = \{ $\lambda$ u, $\exists$ m, m $\in$ (N' F0) /\ m <> F 0
    /\ u = \[[(F 0),m]\]_F0 \}.
Fact Z'_N'_Property2 : $\forall$ F0 u, Arithmetical_ultraFilter F0
  -> u $\in$ (Z' F0) -> u $\in$ ((Z' F0) $\sim$ (Z'_N' F0)) <-> u <_F0 (Z'0 F0).
Fact Z'_N'_Property3 : $\forall$ F0 u, Arithmetical_ultraFilter F0
  -> u $\in$ ((Z'_N' F0) $\sim$ [Z'0 F0])
  -> $\exists$! u0, u0 $\in$ ((Z' F0) $\sim$ (Z'_N' F0)) /\ u +_F0 u0 = Z'0 F0.
Fact Z'_N'_Property4 : $\forall$ F0 u, Arithmetical_ultraFilter F0
  -> u $\in$ ((Z' F0) $\sim$ (Z'_N' F0))
  -> $\exists$! u0, u0 $\in$ ((Z'_N' F0) $\sim$ [Z'0 F0]) /\ u +_F0 u0 = Z'0 F0.
\end{lstlisting}

These simple properties indicate that intuitively the relationship between ${^{*}{\mathbf Z}}_{^{*}{\mathbf N}}$
and $^{*}{\mathbf Z}$ is entirely consistent with the relationship between non-negative integers and all integers;
in other words, the members of ${^{*}{\mathbf Z}}_{^{*}{\mathbf N}}$ are all the non-negative elements in $^{*}{\mathbf Z}$.
Therefore, although naturally inheriting the infinity numbers from $^{*}{\mathbf N}$ due to isomorphism,
${^{*}{\mathbf Z}}_{^{*}{\mathbf N}}$ does not include all infinity elements in $^{*}{\mathbf Z}$ (negative infinities are not in ${^{*}{\mathbf Z}}_{^{*}{\mathbf N}}$),
and next we are to construct all of them.

\subsubsection{Infinities in $^{*}{\mathbf Z}$}

We firstly introduce a proper-subset of $^{*}{\mathbf Z}$:
$${^{*}{\mathbf Z}}_{\mathbf Z} = \{u : \exists m,n \in {^{*}{\mathbf N}}_{\mathbf N} \wedge u = [(m,n)]_{R_{^{*}{\mathbf N}}}\}$$

\begin{lstlisting}
Definition Z'_Z F0 := \{ $\lambda$ u, $\exists$ m n, m $\in$ N'_N /\ n $\in$ N'_N
  /\ u = \[[m,n]\]_F0 \}.

Corollary Z'_Z_is_Set : $\forall$ F0, Arithmetical_ultraFilter F0 -> Ensemble (Z'_Z F0).
Corollary Z'_Z_properSubset_Z' : $\forall$ F0, Arithmetical_ultraFilter F0
  -> ($\forall$ n, F0 <> F n) -> (Z'_Z F0) $\subset$ (Z' F0) /\ (Z'_Z F0) <> (Z' F0).
\end{lstlisting}

It is not hard to see that the elements of equivalence classes in ${^{*}{\mathbf Z}}_{\mathbf Z}$ are all from
${{^{*}{\mathbf N}}_{\mathbf N}} \times {{^{*}{\mathbf N}}_{\mathbf N}}$, and ${^{*}{\mathbf N}}_{\mathbf N}$ is the finite section of $^{*}{\mathbf N}$.
This implies that ${^{*}{\mathbf Z}}_{\mathbf Z}$ is also the finite section of $^{*}{\mathbf Z}$.
The reason we choose this symbol is that intuitively ${^{*}{\mathbf Z}}_{\mathbf Z}$
can serve as a set of general integers. On the one hand, as a subset of $^{*}{\mathbf Z}$,
${^{*}{\mathbf Z}}_{\mathbf Z}$ naturally satisfies all the operation and order properties on $^{*}{\mathbf Z}$,
which are consistent with those on general integers; on the other hand, the operations on $^{*}{\mathbf Z}$ are closed in ${^{*}{\mathbf Z}}_{\mathbf Z}$,
meaning that for any $u,v$ in ${^{*}{\mathbf Z}}_{\mathbf Z}$, both $u +_{^{*}{\mathbf Z}} v$ and $u \cdot_{^{*}{\mathbf Z}} v$ are also members of ${^{*}{\mathbf Z}}_{\mathbf Z}$.

\begin{lstlisting}
Fact Z'_Z_Plus_in_Z'_Z : $\forall$ F0 u v, Arithmetical_ultraFilter F0
  -> u $\in$ (Z'_Z F0) -> v $\in$ (Z'_Z F0) -> (u +_F0 v) $\in$ (Z'_Z F0).
Fact Z'_Z_Mult_in_Z'_Z : $\forall$ F0 u v, Arithmetical_ultraFilter F0
  -> u $\in$ (Z'_Z F0) -> v $\in$ (Z'_Z F0) -> (u $\cdot$_F0 v) $\in$ (Z'_Z F0).
\end{lstlisting}

\noindent Undoubtedly, the zero element $^{*}$Z$_0$ and identity element $^{*}$Z$_1$ are members in ${^{*}{\mathbf Z}}_{\mathbf Z}$.

\begin{lstlisting}
Fact Z'0_in_Z'_Z : $\forall$ F0, Arithmetical_ultraFilter F0 -> (Z'0 F0) $\in$ (Z'_Z F0).
Fact Z'1_in_Z'_Z : $\forall$ F0, Arithmetical_ultraFilter F0 -> (Z'1 F0) $\in$ (Z'_Z F0).
\end{lstlisting}

Thus, the set ${^{*}{\mathbf Z}} \sim {{^{*}{\mathbf Z}}_{\mathbf Z}}$ is the infinite section of $^{*}{\mathbf Z}$.
In contrast with ${^{*}{\mathbf N}} \sim {{^{*}{\mathbf N}}_{\mathbf N}}$,
${^{*}{\mathbf Z}} \sim {{^{*}{\mathbf Z}}_{\mathbf Z}}$ additionally includes negative elements,
and here comes the next proposition.

\begin{proposition}
For every member $\tau$ in (${^{*}{\mathbf Z}} \sim {^{*}{\mathbf Z}}_{\mathbf Z}$) and every member $u$ in ${^{*}{\mathbf Z}}_{\mathbf Z}$,
$\tau <_{^{*}{\mathbf Z}} u$ or $u <_{^{*}{\mathbf Z}} \tau$ holds.
\end{proposition}

\begin{lstlisting}
Fact Z'_infinity : $\forall$ F0 u t, Arithmetical_ultraFilter F0
  -> ($\forall$ n, F0 <> F n) -> u $\in$ (Z'_Z F0) -> t $\in$ ((Z' F0) $\sim$ (Z'_Z F0))
  -> t <_F0 u \/ u <_F0 t.
\end{lstlisting}

Therefore, intuitively the members of ${^{*}{\mathbf Z}} \sim {^{*}{\mathbf Z}}_{\mathbf Z}$ are infinity integers,
they are larger or less than all general integers.
Moreover, the operation and order properties on ${^{*}{\mathbf Z}} \sim {^{*}{\mathbf Z}}_{\mathbf Z}$ are entirely same as those on standard (general) integers.
$^{*}{\mathbf Z}$ is a set of non-standard integers, its structure can be intuitively depicted as:
$$\cdots < -\tau < \cdots \rlap{$\overbrace{\phantom{< -n < \cdots
< -2 < -1 < 0 < 1 < 2 < \cdots < n <}}^{{^{*}{\mathbf Z}}_{\mathbf Z}}$} < -n < \cdots
< -2 < -1 < \underbrace{0 < 1 < 2 < \cdots < n <\cdots < \tau < \cdots \notag}_{{^{*}{\mathbf Z}}_{^{*}{\mathbf N}}}$$
where $\tau$ represents an infinity integer.

\subsection{Extension from $^{*}{\mathbf Z}$ to $^{*}{\mathbf Q}$}\label{sec5-3}
\subsubsection{Construction of $^{*}{\mathbf Q}$}\label{sec5-3-1}

Similar to the process from $^{*}{\mathbf N}$ to $^{*}{\mathbf Z}$,
$^{*}{\mathbf Q}$ is obtained by equivalently classifying the elements of ${^{*}{\mathbf Z}} \times ({^{*}{\mathbf Z}} \sim \{{^{*}}$Z$_0\})$
according to the following equivalence relation:
$$R_{^{*}{\mathbf Z}} = \{((a,b),(c,d)) : a,c \in {^{*}{\mathbf Z}}\ \wedge\ b,d \in ({{^{*}{\mathbf Z}} \sim \{{^{*}}\mbox{Z}_0\}})\ \wedge\ a \cdot_{^{*}{\mathbf Z}} d = b \cdot_{^{*}{\mathbf Z}} c\}.$$

\begin{lstlisting}
Definition Z'_De F0 := (Z' F0) $\times$ ((Z' F0) $\sim$ [Z'0 F0]).

Definition R_Z' F0 := \{\ $\lambda$ u v, $\exists$ a b c d, a $\in$ (Z' F0)
  /\ b $\in$ ((Z' F0) $\sim$ [Z'0 F0]) /\ c $\in$ (Z' F0)
  /\ d $\in$ ((Z' F0) $\sim$ [Z'0 F0]) /\ u = [a,b] /\ v = [c,d]
  /\ a $\cdot$_F0 d = b $\cdot$_F0 c \}\.
Corollary R_Z'_is_equRelation : $\forall$ F0, Arithmetical_ultraFilter F0
  -> equRelation (R_Z' F0) (Z'_De F0).
\end{lstlisting}

\noindent Here, the definition ``{\verb+Z'_De+}" in code gives ${^{*}{\mathbf Z}} \times ({^{*}{\mathbf Z}} \sim \{{^{*}}$Z$_0\})$ a shorter name,
which makes the formal description more convenient.

The relationship between two identical equivalence classes regarding $R_{^{*}{\mathbf N}}$ is:
$$\forall a,c \in {^{*}{\mathbf Z}},\ b,d \in ({{^{*}{\mathbf Z}} \sim \{{^{*}}\mbox{Z}_0\}}),\ [(a,b)]_{R_{^{*}{\mathbf Z}}} = [(c,d)]_{R_{^{*}{\mathbf Z}}} \Longleftrightarrow
a \cdot_{^{*}{\mathbf Z}} d = b \cdot_{^{*}{\mathbf Z}} c.$$

\begin{lstlisting}
Notation "\[ u \]_ F0 " := (equClass u (R_Z' F0) (Z'_De F0))
  (at level 5, u at level 0, F0 at level 0) : q'_scope.

Corollary R_Z'_Corollary : $\forall$ F0 a b c d, Arithmetical_ultraFilter F0
  -> a $\in$ (Z' F0) -> b $\in$ ((Z' F0) $\sim$ [Z'0 F0])
  -> c $\in$ (Z' F0) -> d $\in$ ((Z' F0) $\sim$ [Z'0 F0])
  -> \[[a,b]\]_F0 = \[[c,d]\]_F0 <-> a $\cdot$_F0 d = b $\cdot$_F0 c.
\end{lstlisting}

\begin{figure}[h]
\centering
\includegraphics[width=0.8\textwidth]{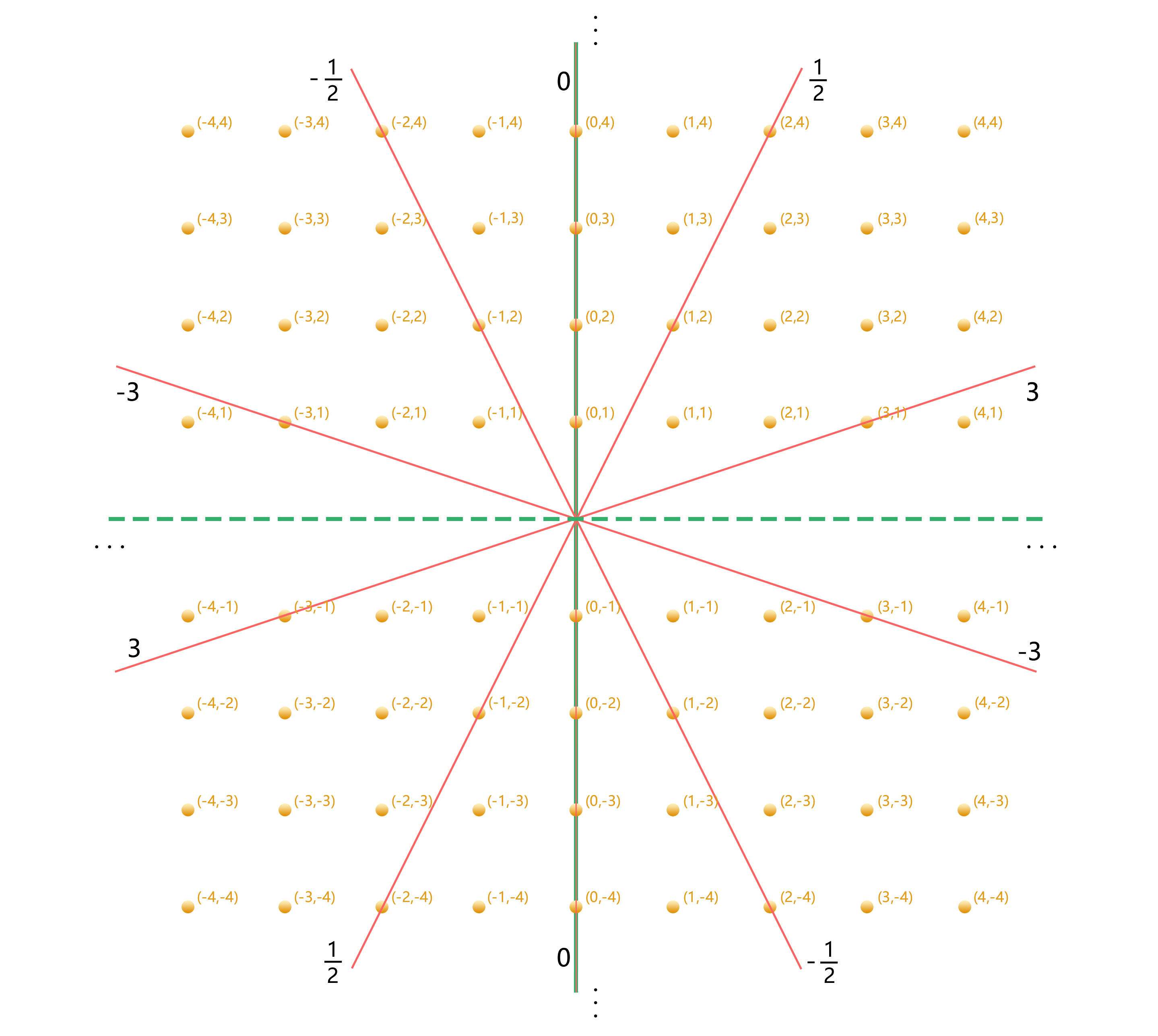}
\caption{equivalence classification of elements in ${^{*}{\mathbf Z}} \times ({^{*}{\mathbf Z}} \sim \{{^{*}}$Z$_0\})$}\label{fig3}
\end{figure}

Figure \ref{fig3} presents how $R_{^{*}{\mathbf Z}}$ classifies ${^{*}{\mathbf Z}} \times ({^{*}{\mathbf Z}} \sim \{{^{*}}$Z$_0\})$.
In the planar dot array composed of elements in ${^{*}{\mathbf Z}} \times ({^{*}{\mathbf Z}} \sim \{{^{*}}$Z$_0\})$,
every line connecting equivalent points represents an equivalence class of ${^{*}{\mathbf Z}} \times ({^{*}{\mathbf Z}} \sim \{{^{*}}$Z$_0\})$, and certainly the lines are dense,
which means that there always exists such a line between any two and the inverse elements such as $\frac{1}{2}$ are contained.
The set consisting of all these lines, namely the quotient set of ${^{*}{\mathbf Z}} \times ({^{*}{\mathbf Z}} \sim \{{^{*}}$Z$_0\})$ regarding $R_{^{*}{\mathbf Z}}$,
is exactly ${^{*}{\mathbf Q}}$.
$${^{*}{\mathbf Q}} = ({^{*}{\mathbf Z}} \times ({^{*}{\mathbf Z}} \sim \{{^{*}}\mbox{Z}_0\})) / {R_{^{*}{\mathbf Z}}} = \{u : \exists x, u = [x]_{R_{^{*}{\mathbf Z}}}\}.$$

\begin{lstlisting}
Definition Q' F0 := (Z'_De F0) / (R_Z' F0).
Corollary Q'_is_Set : $\forall$ F0, F0 $\in$ $\beta\omega$ -> Ensemble (Q' F0).

Corollary Q'_Corollary1 : $\forall$ F0 u, Arithmetical_ultraFilter F0 -> u $\in$ (Q' F0)
  -> $\exists$ a b, a $\in$ (Z' F0) /\ b $\in$ (Z' F0 $\sim$ [Z'0 F0]) /\ b $\in$ (Z' F0)
    /\ b <> Z'0 F0 /\ u = \[[a,b]\]_F0.
Corollary Q'_Corollary2 : $\forall$ F0 a b, Arithmetical_ultraFilter F0
  -> a $\in$ (Z' F0) -> b $\in$ (Z' F0 $\sim$ [Z'0 F0]) -> \[[a,b]\]_F0 $\in$ (Q' F0).
\end{lstlisting}

The zero element $^{*}$Q$_0$ and identity element $^{*}$Q$_1$ of ${^{*}{\mathbf Q}}$ are defined as:
$${^{*}}\mbox{Q}_0 = \{(u,v) : u = {^{*}}\mbox{Z}_0\ \wedge\ v \in {^{*}{\mathbf Z}} \times ({^{*}{\mathbf Z}} \sim \{{^{*}}\mbox{Z}_0\})\}
= [({^{*}}\mbox{Z}_0,{^{*}}\mbox{Z}_1)]_{R_{^{*}{\mathbf Z}}},$$
$${^{*}}\mbox{Q}_1 = \{(u,v) : v \in {^{*}{\mathbf Z}} \times ({^{*}{\mathbf Z}} \sim \{{^{*}}\mbox{Z}_0\})\ \wedge\ u = v \}
= [({^{*}}\mbox{Z}_1,{^{*}}\mbox{Z}_1)]_{R_{^{*}{\mathbf Z}}},$$
where $^{*}$Z$_0$ and $^{*}$Z$_1$ respectively stands for zero and identity element in ${^{*}{\mathbf Z}}$.

\begin{lstlisting}
Definition Q'0 F0 := \{\ $\lambda$ u v, u = (Z'0 F0) /\ v $\in$ ((Z' F0) $\sim$ [Z'0 F0]) \}\.
Corollary Q'0_Corollary : $\forall$ F0, Arithmetical_ultraFilter F0
  -> Q'0 F0 = \[[(Z'0 F0),(Z'1 F0)]\]_F0.
Corollary Q'0_in_Q' : $\forall$ F0, Arithmetical_ultraFilter F0 -> (Q'0 F0) $\in$ (Q' F0).

Definition Q'1 F0 := \{\ $\lambda$ u v, v $\in$ ((Z' F0) $\sim$ [Z'0 F0]) /\ u = v \}\.
Corollary Q'1_Corollary : $\forall$ F0, Arithmetical_ultraFilter F0
  -> Q'1 F0 = \[[(Z'1 F0),(Z'1 F0)]\]_F0.
Corollary Q'1_in_Q' : $\forall$ F0, Arithmetical_ultraFilter F0 -> (Q'1 F0) $\in$ (Q' F0).

Fact Q'0_isn't_Q'1 : $\forall$ F0, Arithmetical_ultraFilter F0 -> Q'0 F0 <> Q'1 F0.
\end{lstlisting}

\subsubsection{Order on $^{*}{\mathbf Q}$}\label{sec5-3-2}

\begin{definition}[Order on $^{*}{\mathbf Q}$]
Let $u,v$ be members in $^{*}{\mathbf Q}$, and $u = [(a,b)]_{R_{^{*}{\mathbf Z}}}, v = [(c,d)]_{R_{^{*}{\mathbf Z}}}$
($a,c \in {^{*}{\mathbf Z}},\ b,d \in ({^{*}{\mathbf Z}} \sim \{{^{*}}$Z$_0\})$), then the order of $u$ and $v$ is defined as follows:
$$u <_{^{*}{\mathbf Q}} v\ \Longleftrightarrow\
[(a,b)]_{R_{^{*}{\mathbf Z}}} <_{^{*}{\mathbf Q}} [(c,d)]_{R_{^{*}{\mathbf Z}}}\ \Longleftrightarrow\
a \cdot_{^{*}{\mathbf Z}} d\ <_{^{*}{\mathbf Z}}\ b \cdot_{^{*}{\mathbf Z}} c.$$
\end{definition}

The definition of this order relation is reasonable,
because for any $a,a_1,c,c_1, \in {^{*}{\mathbf Z}}$ and $b,b_1,d,d_1 \in ({^{*}{\mathbf Z}} \sim \{{^{*}}$Z$_0\})$, if
$$u = [(a,b)]_{R_{^{*}{\mathbf Z}}} = [(a_1,b_1)]_{R_{^{*}{\mathbf Z}}},\
v = [(c,d)]_{R_{^{*}{\mathbf Z}}} = [(c_1,d_1)]_{R_{^{*}{\mathbf Z}}},$$
then it can be verified that
$$u <_{^{*}{\mathbf Q}} v\ \Longleftrightarrow\
a \cdot_{^{*}{\mathbf Z}} d\ <_{^{*}{\mathbf Z}}\ b \cdot_{^{*}{\mathbf Z}} c\ \Longleftrightarrow\
a_1 \cdot_{^{*}{\mathbf Z}} d_1\ <_{^{*}{\mathbf Z}}\ b_1 \cdot_{^{*}{\mathbf Z}} c_1.$$
Namely, the order relation on $^{*}{\mathbf Q}$ is independent of the selection of representations of equivalence classes.

\begin{lstlisting}
Definition Q'_Ord F0 := \{\ $\lambda$ u v, $\forall$ a b c d, a $\in$ (Z' F0)
  -> b $\in$ ((Z' F0) $\sim$ [Z'0 F0])
  -> c $\in$ (Z' F0) -> d $\in$ ((Z' F0) $\sim$ [Z'0 F0])
  -> (Z'0 F0) <_F0 b -> (Z'0 F0) <_F0 d
  -> u = \[[a,b]\]_F0 -> v = \[[c,d]\]_F0
  -> (a $\cdot$_F0 d) <_F0 (b $\cdot$_F0 c) \}\.
Notation "u <_ F0 v" := (Rrelation u (Q'_Ord F0) v)
  (at level 10, F0 at level 9) : q'_scope.

Corollary Q'_Ord_reasonable :
  $\forall$ F0 a b c d a1 b1 c1 d1, Arithmetical_ultraFilter F0
  -> a $\in$ (Z' F0) -> b $\in$ ((Z' F0) $\sim$ [Z'0 F0])
  -> c $\in$ (Z' F0) -> d $\in$ ((Z' F0) $\sim$ [Z'0 F0])
  -> a1 $\in$ (Z' F0) -> b1 $\in$ ((Z' F0) $\sim$ [Z'0 F0])
  -> c1 $\in$ (Z' F0) -> d1 $\in$ ((Z' F0) $\sim$ [Z'0 F0])
  -> \[[a,b]\]_F0 = \[[a1,b1]\]_F0
  -> \[[c,d]\]_F0 = \[[c1,d1]\]_F0
  -> ((Z'0 F0) <_F0 b)$\mbox{\color[RGB]{0,100,0}\%}$z' -> ((Z'0 F0) <_F0 d)$\mbox{\color[RGB]{0,100,0}\%}$z'
  -> ((Z'0 F0) <_F0 b1)$\mbox{\color[RGB]{0,100,0}\%}$z' -> ((Z'0 F0) <_F0 d1)$\mbox{\color[RGB]{0,100,0}\%}$z'
  -> ((a $\cdot$_F0 d) <_F0 (b $\cdot$_F0 c)
    <-> (a1 $\cdot$_F0 d1) <_F0 (b1 $\cdot$_F0 c1))$\mbox{\color[RGB]{0,100,0}\%}$z'.
Corollary Q'_Ord_Corollary :
  $\forall$ F0 a b c d, Arithmetical_ultraFilter F0
  -> a $\in$ (Z' F0) -> b $\in$ ((Z' F0) $\sim$ [Z'0 F0])
  -> c $\in$ (Z' F0) -> d $\in$ ((Z' F0) $\sim$ [Z'0 F0])
  -> ((Z'0 F0) <_F0 b)$\mbox{\color[RGB]{0,100,0}\%}$z' -> ((Z'0 F0) <_F0 d)$\mbox{\color[RGB]{0,100,0}\%}$z'
  -> \[[a,b]\]_F0 <_F0 \[[c,d]\]_F0
    <-> ((a $\cdot$_F0 d) <_F0 (b $\cdot$_F0 c))$\mbox{\color[RGB]{0,100,0}\%}$z'.
\end{lstlisting}

The properties of this order are totally same as those of $^{*}{\mathbf N}$ and $^{*}{\mathbf Z}$.

\begin{lstlisting}
Fact Q'_Ord_irReflex : $\forall$ F0 u v, Arithmetical_ultraFilter F0
  -> u $\in$ (Q' F0) -> v $\in$ (Q' F0) -> u = v -> $\sim$ u <_F0 v.
Fact Q'_Ord_Trans : $\forall$ F0 u v w, Arithmetical_ultraFilter F0
  -> u $\in$ (Q' F0) -> v $\in$ (Q' F0) -> w $\in$ (Q' F0) -> u <_F0 v
  -> v <_F0 w -> u <_F0 w.
Fact Q'_Ord_Tri : $\forall$ F0, Arithmetical_ultraFilter F0
  -> Connect (Q'_Ord F0) (Q' F0).
Fact Q'0_lt_Q'1 : $\forall$ F0, Arithmetical_ultraFilter F0
  -> (Q'0 F0) <_F0 (Q'1 F0).
\end{lstlisting}

\subsubsection{Operations on $^{*}{\mathbf Q}$}\label{sec5-3-3}

\begin{definition}[Addition on $^{*}{\mathbf Q}$]
Let $u,v$ be members in $^{*}{\mathbf Q}$, and $u = [(a,b)]_{R_{^{*}{\mathbf Z}}}, v = [(c,d)]_{R_{^{*}{\mathbf Z}}}$
($a,c \in {^{*}{\mathbf Z}},\ b,d \in ({^{*}{\mathbf Z}} \sim \{{^{*}}$Z$_0\})$), then the sum of $u$ and $v$ is defined as follows:
$$u +_{^{*}{\mathbf Q}} v\ =\ [(a,b)]_{R_{^{*}{\mathbf Z}}}\ +_{^{*}{\mathbf Q}}\ [(c,d)]_{R_{^{*}{\mathbf Z}}}\
=\ [((a \cdot_{^{*}{\mathbf Z}} d) +_{^{*}{\mathbf Z}} (b \cdot_{^{*}{\mathbf Z}} c)\ ,\ b \cdot_{^{*}{\mathbf Z}} d)]_{R_{^{*}{\mathbf Z}}}.$$
\end{definition}

The definition of this addition is reasonable,
because for any $a,a_1,c,c_1$ in ${^{*}{\mathbf Z}}$, and $b,b_1,d,d_1$ in ${^{*}{\mathbf Z}} \sim \{{^{*}}$Z$_0\}$, if
$$u = [(a,b)]_{R_{^{*}{\mathbf Z}}} = [(a_1,b_1)]_{R_{^{*}{\mathbf Z}}},\
v = [(c,d)]_{R_{^{*}{\mathbf Z}}} = [(c_1,d_1)]_{R_{^{*}{\mathbf Z}}},$$
then it can be verified that
$$u +_{^{*}{\mathbf Q}} v\ =\ [(a,b)]_{R_{^{*}{\mathbf Z}}}\ +_{^{*}{\mathbf Q}}\ [(c,d)]_{R_{^{*}{\mathbf Z}}}\ =\
[(a_1,b_1)]_{R_{^{*}{\mathbf Z}}}\ +_{^{*}{\mathbf Q}}\ [(c_1,d_1)]_{R_{^{*}{\mathbf Z}}}.$$
Namely, the result of $u +_{^{*}{\mathbf Q}} v$ is unique and independent of the selection of representations of equivalence classes.

\begin{lstlisting}
Definition Q'_Plus F0 u v := $\cap$(\{ $\lambda$ w,
  $\forall$ a b c d, a $\in$ (Z' F0) -> b $\in$ ((Z' F0) $\sim$ [Z'0 F0])
  -> c $\in$ (Z' F0) -> d $\in$ ((Z' F0) $\sim$ [Z'0 F0])
  -> u = \[[a,b]\]_F0 -> v = \[[c,d]\]_F0
  -> w = \[[((a $\cdot$_F0 d) +_F0 (b $\cdot$_F0 c))$\mbox{\color[RGB]{0,100,0}\%}$z',(b $\cdot$_F0 d)$\mbox{\color[RGB]{0,100,0}\%}$z']\]_F0 \}).
Notation "u +_ F0 v" := (Q'_Plus F0 u v)
  (at level 10, F0 at level 9) : q'_scope.

Corollary Q'_Plus_reasonable : $\forall$ F0 u v, Arithmetical_ultraFilter F0
  -> u $\in$ (Q' F0) -> v $\in$ (Q' F0) -> $\exists$ w, w $\in$ (Q' F0) /\ [w]
    = \{ $\lambda$ w, $\forall$ a b c d, a $\in$ (Z' F0) -> b $\in$ ((Z' F0) $\sim$ [Z'0 F0])
    -> c $\in$ (Z' F0) -> d $\in$ ((Z' F0) $\sim$ [Z'0 F0])
    -> u = \[[a,b]\]_F0 -> v = \[[c,d]\]_F0
    -> w = \[[((a $\cdot$_F0 d) +_F0 (b $\cdot$_F0 c))$\mbox{\color[RGB]{0,100,0}\%}$z',(b $\cdot$_F0 d)$\mbox{\color[RGB]{0,100,0}\%}$z']\]_F0 \}.
Corollary Q'_Plus_Corollary : $\forall$ F0 a b c d, Arithmetical_ultraFilter F0
  -> a $\in$ (Z' F0) -> b $\in$ ((Z' F0) $\sim$ ([Z'0 F0]))
  -> c $\in$ (Z' F0) -> d $\in$ ((Z' F0) $\sim$ ([Z'0 F0]))
  -> \[[a,b]\]_F0 +_F0 \[[c,d]\]_F0
    = \[[((a $\cdot$_F0 d) +_F0 (b $\cdot$_F0 c))$\mbox{\color[RGB]{0,100,0}\%}$z',(b $\cdot$_F0 d)$\mbox{\color[RGB]{0,100,0}\%}$z']\]_F0.
\end{lstlisting}

\begin{definition}[Multiplication on $^{*}{\mathbf Q}$]
Let $u,v$ be members in $^{*}{\mathbf Q}$ and $u = [(a,b)]_{R_{^{*}{\mathbf Z}}}, v = [(c,d)]_{R_{^{*}{\mathbf Z}}}$
$(a,c \in {^{*}{\mathbf Z}},\ b,d \in ({^{*}{\mathbf Z}} \sim \{{^{*}}$Z$_0\}))$, then the product of $u$ and $v$ is defined as follows:
$$u \cdot_{^{*}{\mathbf Q}} v\ =\ [(a,b)]_{R_{^{*}{\mathbf Z}}}\ \cdot_{^{*}{\mathbf Q}}\ [(c,d)]_{R_{^{*}{\mathbf Z}}}\
=\ [(a \cdot_{^{*}{\mathbf Z}} c\ ,\ b \cdot_{^{*}{\mathbf Z}} d)]_{R_{^{*}{\mathbf Z}}}.$$
\end{definition}

The definition of this multiplication is reasonable,
because for any $a,a_1,c,c_1 \in {^{*}{\mathbf Z}}$, and $b,b_1,d,d_1 \in ({^{*}{\mathbf Z}} \sim \{{^{*}}$Z$_0\})$, if
$$u = [(a,b)]_{R_{^{*}{\mathbf Z}}} = [(a_1,b_1)]_{R_{^{*}{\mathbf Z}}},\
v = [(c,d)]_{R_{^{*}{\mathbf Z}}} = [(c_1,d_1)]_{R_{^{*}{\mathbf Z}}},$$
then it can be verified that
$$u \cdot_{^{*}{\mathbf Q}} v\ =\ [(a,b)]_{R_{^{*}{\mathbf Z}}}\ \cdot_{^{*}{\mathbf Q}}\ [(c,d)]_{R_{^{*}{\mathbf Z}}}\ =\
[(a_1,b_1)]_{R_{^{*}{\mathbf Z}}}\ \cdot_{^{*}{\mathbf Q}}\ [(c_1,d_1)]_{R_{^{*}{\mathbf Z}}}.$$
Namely, the result of $u \cdot_{^{*}{\mathbf Q}} v$ is unique and independent of the selection of representations of equivalence classes.

\begin{lstlisting}
Definition Q'_Mult F0 u v := $\cap$(\{ $\lambda$ w,
  $\forall$ a b c d, a $\in$ (Z' F0) -> b $\in$ ((Z' F0) $\sim$ [Z'0 F0])
  -> c $\in$ (Z' F0) -> d $\in$ ((Z' F0) $\sim$ [Z'0 F0])
  -> u = \[[a,b]\]_F0 -> v = \[[c,d]\]_F0
  -> w = \[[(a $\cdot$_F0 c)$\mbox{\color[RGB]{0,100,0}\%}$z',(b $\cdot$_F0 d)$\mbox{\color[RGB]{0,100,0}\%}$z']\]_F0 \}).
Notation "u $\cdot$_ F0 v" := (Q'_Mult F0 u v)
  (at level 10, F0 at level 9) : q'_scope.

Fact Q'_Mult_reasonable : $\forall$ F0 u v, Arithmetical_ultraFilter F0
  -> u $\in$ (Q' F0) -> v $\in$ (Q' F0) -> $\exists$ w, w $\in$ (Q' F0) /\ [w]
    = \{ $\lambda$ w, $\forall$ a b c d, a $\in$ (Z' F0) -> b $\in$ ((Z' F0) $\sim$ [Z'0 F0])
    -> c $\in$ (Z' F0) -> d $\in$ ((Z' F0) $\sim$ [Z'0 F0])
    -> u = \[[a,b]\]_F0 -> v = \[[c,d]\]_F0
    -> w = \[[(a $\cdot$_F0 c)$\mbox{\color[RGB]{0,100,0}\%}$z',(b $\cdot$_F0 d)$\mbox{\color[RGB]{0,100,0}\%}$z']\]_F0 \}.
Corollary Q'_Mult_Corollary : $\forall$ F0 a b c d, Arithmetical_ultraFilter F0
  -> a $\in$ (Z' F0) -> b $\in$ ((Z' F0) $\sim$ [Z'0 F0])
  -> c $\in$ (Z' F0) -> d $\in$ ((Z' F0) $\sim$ [Z'0 F0])
  -> \[[a,b]\]_F0 $\cdot$_F0 \[[c,d]\]_F0
    = \[[(a $\cdot$_F0 c)$\mbox{\color[RGB]{0,100,0}\%}$z',(b $\cdot$_F0 d)$\mbox{\color[RGB]{0,100,0}\%}$z']\]_F0.
\end{lstlisting}

The multiplication on $^{*}{\mathbf Q}$ possesses an additional property compared to that on $^{*}{\mathbf Z}$:
for every $u$ in $^{*}{\mathbf Q}$, where $u \neq {^{*}}$Q$_0$,
there exists a unique element $u_1$ of $^{*}{\mathbf Q}$ such that $u_1 \neq {^{*}}$Q$_0$ and $u \cdot_{^{*}{\mathbf Q}} u_1 = {^{*}}$Q$_1$.
Namely, every non-zero element in $^{*}{\mathbf Q}$ has a unique corresponding non-zero inverse element regarding multiplication on $^{*}{\mathbf Q}$.

\begin{lstlisting}
Fact Q'_inv : $\forall$ F0 u, Arithmetical_ultraFilter F0
  -> u $\in$ (Q' F0) -> u <> Q'0 F0 -> $\exists$! u1, u1 $\in$ (Q' F0)
    /\ u1 <> Q'0 F0 /\ u $\cdot$_F0 u1 = Q'1 F0.
\end{lstlisting}

According to the structure of $^{*}{\mathbf Q}$, here comes the consequence that
if $[(a,b)]_{R_{^{*}{\mathbf Z}}} \cdot_{^{*}{\mathbf Q}} u = {^{*}}$Q$_1$, then $u = [(b,a)]_{R_{^{*}{\mathbf Z}}}$,
where $a,b$ are members in ${^{*}{\mathbf Z}} \sim \{{^{*}}$Z$_0\}$ and $u$ in $^{*}{\mathbf Q}$.

\begin{lstlisting}
Fact Q'_inv_Corollary : $\forall$ F0 a b u, Arithmetical_ultraFilter F0
  -> a $\in$ ((Z' F0) $\sim$ [Z'0 F0]) -> b $\in$ ((Z' F0) $\sim$ [Z'0 F0])
  -> u $\in$ (Q' F0) -> \[[a,b]\]_F0 $\cdot$_F0 u = Q'1 F0 -> u = \[[b,a]\]_F0.
\end{lstlisting}

\noindent This consequence reveals the structural connection between every element and its inverse in ${^{*}{\mathbf Q}} \sim \{{^{*}}$Q$_0\}$.

\begin{definition}[Subtraction on $^{*}{\mathbf Q}$]
Let $u,v,w$ be members in $^{*}{\mathbf Q}$, and $u +_{^{*}{\mathbf Q}} w = v$, then $v$ minus $u$, denoted as $v -_{^{*}{\mathbf Q}} u$, is $w$.
\end{definition}

The definition of this subtraction is reasonable,
because it can be verified that the existence of the very $w$ making $u +_{^{*}{\mathbf Q}} w = v$ hold is unique.

\begin{lstlisting}
Definition Q'_Minus F0 v u := $\cap$(\{ $\lambda$ w, w $\in$ (Q' F0) /\ u +_F0 w = v \}).
Notation "u -_ F0 v" := (Q'_Minus F0 u v)
  (at level 10, F0 at level 9) : q'_scope.

Corollary Q'_Minus_reasonable : $\forall$ F0 u v, Arithmetical_ultraFilter F0
  -> u $\in$ (Q' F0) -> v $\in$ (Q' F0) -> $\exists$! w, w $\in$ (Q' F0) /\ u +_F0 w = v.
Corollary Q'_Minus_Corollary : $\forall$ F0 u w v, Arithmetical_ultraFilter F0
  -> u $\in$ (Q' F0) -> w $\in$ (Q' F0) -> v $\in$ (Q' F0)
  -> u +_F0 w = v <-> v -_F0 u = w.
\end{lstlisting}

\begin{definition}[Division on $^{*}{\mathbf Q}$]
Let $u,v,w$ be members in $^{*}{\mathbf Q}$, $u \neq {^{*}}\mbox{\rm Q}_0$, and $u \cdot_{^{*}{\mathbf Q}} w = v$, then $v$ divided by $u$, denoted as $v /_{^{*}{\mathbf Q}} u$, is $w$.
\end{definition}

The definition of this division is reasonable,
because it can be verified that the existence of the very $w$ making $u /_{^{*}{\mathbf Q}} w = v$ hold is unique.

\begin{lstlisting}
Definition Q'_Divide F0 v u := $\cap$(\{ $\lambda$ w, w $\in$ (Q' F0)
  /\ u <> Q'0 F0 /\ u $\cdot$_F0 w = v \}).
Notation "v /_ F0 u" := (Q'_Divide F0 v u)
  (at level 10, F0 at level 9) : q'_scope.

Corollary Q'_Divide_reasonable : $\forall$ F0 u v, Arithmetical_ultraFilter F0
  -> u $\in$ (Q' F0) -> v $\in$ (Q' F0) -> u <> Q'0 F0
  -> $\exists$! w, w $\in$ (Q' F0) /\ u $\cdot$_F0 w = v.
Corollary Q'_Divide_Corollary : $\forall$ F0 u w v, Arithmetical_ultraFilter F0
  -> u $\in$ (Q' F0) -> u <> Q'0 F0 -> w $\in$ (Q' F0) -> v $\in$ (Q' F0)
  -> u $\cdot$_F0 w = v <-> v /_F0 u = w.
\end{lstlisting}

We present here some implementation for fundamental properties of the subtraction and division.

\begin{lstlisting}
Fact Q'_Minus_in_Q' : $\forall$ F0 v u, Arithmetical_ultraFilter F0
  -> v $\in$ (Q' F0) -> u $\in$ (Q' F0) -> (v -_F0 u) $\in$ (Q' F0).
Fact Q'_Divide_in_Q' : $\forall$ F0 v u, Arithmetical_ultraFilter F0
  -> v $\in$ (Q' F0) -> u $\in$ (Q' F0) -> u <> Q'0 F0 -> (v /_F0 u) $\in$ (Q' F0).

Fact Q'_Minus_Property1 : $\forall$ F0 u, Arithmetical_ultraFilter F0
  -> u $\in$ (Q' F0) -> u -_F0 u = Q'0 F0.
Fact Q'_Minus_Property2 : $\forall$ F0 u, Arithmetical_ultraFilter F0
  -> u $\in$ (Q' F0) -> u -_F0 (Q'0 F0) = u.

Fact Q'_Divide_Property1 : $\forall$ F0 u, Arithmetical_ultraFilter F0
  -> u $\in$ (Q' F0) -> u <> Q'0 F0 -> u /_F0 u = Q'1 F0.
Fact Q'_Divide_Property2 : $\forall$ F0 u, Arithmetical_ultraFilter F0
  -> u $\in$ (Q' F0) -> u /_F0 (Q'1 F0) = u.
Fact Q'_Divide_Property3 : $\forall$ F0 v u, Arithmetical_ultraFilter F0
  -> v $\in$ (Q' F0) -> u $\in$ (Q' F0) -> u <> Q'0 F0 -> u $\cdot$_F0 (v /_F0 u) = v.

Fact Q'_Mix_Association1 : $\forall$ F0 u v w, Arithmetical_ultraFilter F0
  -> u $\in$ (Q' F0) -> v $\in$ (Q' F0) -> w $\in$ (Q' F0)
  -> (u +_F0 v) -_F0 w = u +_F0 (v -_F0 w).
Fact Q'_Mix_Association2 : $\forall$ F0 u v w, Arithmetical_ultraFilter F0
  -> u $\in$ (Q' F0) -> v $\in$ (Q' F0) -> w $\in$ (Q' F0) -> w <> Q'0 F0
  -> u $\cdot$_F0 (v /_F0 w) = (u $\cdot$_F0 v) /_F0 w.
Fact Q'_Mult_Distribution_Minus : $\forall$ F0 u v w, Arithmetical_ultraFilter F0
  -> u $\in$ (Q' F0) -> v $\in$ (Q' F0) -> w $\in$ (Q' F0)
  -> u $\cdot$_F0 (v -_F0 w) = (u $\cdot$_F0 v) -_F0 (u $\cdot$_F0 w).
Fact Q'_Divide_Distribution : $\forall$ F0 u v w, Arithmetical_ultraFilter F0
  -> u $\in$ (Q' F0) -> v $\in$ (Q' F0) -> w $\in$ (Q' F0) -> w <> Q'0 F0
  -> (u +_F0 v) /_F0 w = (u /_F0 w) +_F0 (v /_F0 w).
Fact Q'_Frac_Square : $\forall$ F0 u v, Arithmetical_ultraFilter F0
  -> u $\in$ (Q' F0) -> v $\in$ (Q' F0) -> v <> Q'0 F0
  -> (u $\cdot$_F0 u) /_F0 (v $\cdot$_F0 v) = (u /_F0 v) $\cdot$_F0 (u /_F0 v).
\end{lstlisting}

\begin{definition}[Absolute Value on $^{*}{\mathbf Q}$]
Let $x$ be a member in $^{*}{\mathbf Q}$, then the absolute value of $x$ (denoted as $|x|$) is its image under the following class:
$$\{(u,v) : u \in {^{*}\mathbf Q}\ \wedge\ ((u <_{^{*}\mathbf Q} {^{*}}\mbox{\rm Q}_0\ \wedge\ v = {^{*}}\mbox{\rm Q}_0\ -_{^{*}\mathbf Q}\ u)\
\vee\ ({^{*}}\mbox{\rm Q}_0 <_{^{*}\mathbf Q} u\ \wedge\ v = u)\ \vee\ ({^{*}}\mbox{\rm Q}_0 = u = v))\},$$
which can be verified a function whose domain is $^{*}{\mathbf Q}$ and range is contained in $^{*}{\mathbf Q}$.
\end{definition}

\begin{lstlisting}
Definition Q'_Abs F0 := \{\ $\lambda$ u v, u $\in$ (Q' F0)
  /\ ((u <_F0 (Q'0 F0) /\ v = (Q'0 F0) -_F0 u)
    \/ (u = Q'0 F0 /\ v = Q'0 F0)
    \/ ((Q'0 F0) <_F0 u /\ v = u)) \}\.
Notation "| u |_ F0 " := ((Q'_Abs F0)[u])
  (at level 5, u at level 0, F0 at level 0) : q'_scope.

Corollary Q'Abs_Corollary : $\forall$ F0, Arithmetical_ultraFilter F0
  -> Function (Q'_Abs F0) /\ dom(Q'_Abs F0) = Q' F0
    /\ ran(Q'_Abs F0) $\subset$ (Q' F0).
Corollary Q'Abs_in_Q' : $\forall$ F0 u, Arithmetical_ultraFilter F0
  -> u $\in$ (Q' F0) -> |u|_F0 $\in$ (Q' F0).
\end{lstlisting}

Like the properties of operations, despite straightforward and familiar to us all, the properties of absolute value should be formalized and verified.

\begin{lstlisting}
Fact mt_Q'0_Q'Abs : $\forall$ F0 u, Arithmetical_ultraFilter F0
  -> u $\in$ (Q' F0) -> (Q'0 F0) <_F0 u -> |u|_F0 = u.
Fact eq_Q'0_Q'Abs : $\forall$ F0 u, Arithmetical_ultraFilter F0
  -> u $\in$ (Q' F0) -> |u|_F0 = Q'0 F0 <-> u = Q'0 F0.
Fact lt_Q'0_Q'Abs : $\forall$ F0 u, Arithmetical_ultraFilter F0
  -> u $\in$ (Q' F0) -> u <_F0 (Q'0 F0) -> |u|_F0 = (Q'0 F0) -_F0 u.
Fact Q'0_le_Q'Abs : $\forall$ F0 u, Arithmetical_ultraFilter F0 ->
  u $\in$ (Q' F0) -> |u|_F0 $\in$ (Q' F0) /\
    (|u|_F0 = Q'0 F0 \/ (Q'0 F0) <_F0 |u|_F0).
Fact Self_le_Q'Abs : $\forall$ F0 u, Arithmetical_ultraFilter F0
  -> u $\in$ (Q' F0) -> u = |u|_F0 \/ u <_F0 |u|_F0.
Fact Q'Abs_PrMult : $\forall$ F0 u v, Arithmetical_ultraFilter F0
  -> u $\in$ (Q' F0) -> v $\in$ (Q' F0)
  -> |u $\cdot$_F0 v|_F0 = |u|_F0 $\cdot$_F0 |v|_F0.
Fact Q'Abs_inEqu : $\forall$ F0 u v, Arithmetical_ultraFilter F0
  -> u $\in$ (Q' F0) -> v $\in$ (Q' F0)
  -> |u +_F0 v|_F0 = |u|_F0 +_F0 |v|_F0
    \/ |u +_F0 v|_F0 <_F0 (|u|_F0 +_F0 |v|_F0).
Fact neg_Q'Abs_Equ : $\forall$ F0 u v, Arithmetical_ultraFilter F0
  -> u $\in$ (Q' F0) -> v $\in$ (Q' F0) -> u +_F0 v = Q'0 F0 -> |u|_F0 = |v|_F0.
\end{lstlisting}

\subsubsection{$^{*}{\mathbf N}$ and $^{*}{\mathbf Z}$ in $^{*}{\mathbf Q}$}\label{sec5-3-4}

To clarify that $^{*}{\mathbf Q}$ is a proper extension of $^{*}{\mathbf Z}$,
we present the verification that $^{*}{\mathbf Z}$ along with $^{*}{\mathbf N}$ can be isomorphically embedded in $^{*}{\mathbf Q}$.

\begin{theorem}\label{z'-in-q'}
Let $^{*}{\mathbf Q}$ be constructed with $^{*}{\mathbf Z}$, the following set:
$${^{*}{\mathbf Q}}_{^{*}{\mathbf Z}} = \{u : \exists a, a \in {^{*}{\mathbf Z}} \wedge u = [(a,{^{*}}\mbox{\rm Z}_1)]_{R_{^{*}{\mathbf Z}}}\}$$
is a proper-subset of $^{*}{\mathbf Q}$.
Meanwhile, a 1-1 function whose domain is $^{*}{\mathbf Z}$ and range is ${^{*}{\mathbf Q}}_{^{*}{\mathbf Z}}$ can be constructed:
$$\varphi_2 = \{(u,v) : u \in {^{*}{\mathbf Z}} \wedge v = [(u,{^{*}}\mbox{\rm Z}_1)]_{R_{^{*}{\mathbf Z}}}\},$$
and for any $a,b$ in ${^{*}{\mathbf Z}}$, $\varphi_2$ satisfies:

1) $\varphi_2({^{*}}\mbox{\rm Z}_0) = {^{*}}\mbox{\rm Q}_0\ \wedge\ \varphi_2({^{*}}\mbox{\rm Z}_1) = {^{*}}\mbox{\rm Q}_1$.

2) Order Preservation: $a <_{^{*}{\mathbf Z}} b \Longleftrightarrow \varphi_2(a) <_{^{*}{\mathbf Q}} \varphi_2(b)$.

3) Addition Preservation: $\varphi_2(a +_{^{*}{\mathbf Z}} b) = \varphi_2(a) +_{^{*}{\mathbf Q}} \varphi_2(b)$.

4) Multiplication Preservation: $\varphi_2(a \cdot_{^{*}{\mathbf Z}} b) = \varphi_2(a) \cdot_{^{*}{\mathbf Q}} \varphi_2(b)$.
\end{theorem}

\begin{lstlisting}
Definition Q'_Z' F0 := \{ $\lambda$ u, $\exists$ m, m $\in$ (Z' F0) /\ u = \[[m,(Z'1 F0)]\]_F0 \}.
Corollary Q'_Z'_is_Set : $\forall$ F0, Arithmetical_ultraFilter F0 -> Ensemble (Q'_Z' F0).
Corollary Q'_Z'_properSubset_Q' : $\forall$ F0, Arithmetical_ultraFilter F0
  -> (Q'_Z' F0) $\subset$ (Q' F0) /\ (Q'_Z' F0) <> (Q' F0).

Definition $\varphi$2 F0 := \{\ $\lambda$ u v, u $\in$ (Z' F0)
  /\ v = \[[u,(Z'1 F0)]\]_F0 \}\.
Corollary $\varphi$2_is_Function : $\forall$ F0, Arithmetical_ultraFilter F0
  -> Function1_1 ($\varphi$2 F0) /\ dom($\varphi$2 F0) = (Z' F0) /\ ran($\varphi$2 F0) = (Q'_Z' F0).

Corollary $\varphi$2_Z'0 : $\forall$ F0, Arithmetical_ultraFilter F0
  -> ($\varphi$2 F0)[Z'0 F0] = Q'0 F0.
Corollary $\varphi$2_Z'1 : $\forall$ F0, Arithmetical_ultraFilter F0
  -> ($\varphi$2 F0)[Z'1 F0] = Q'1 F0.
Corollary $\varphi$2_PrOrder : $\forall$ F0 a b, Arithmetical_ultraFilter F0
  -> a $\in$ (Z' F0) -> b $\in$ (Z' F0)
  -> (a <_F0 b)$\mbox{\color[RGB]{0,100,0}\%}$z' <-> ($\varphi$2 F0)[a] <_F0 ($\varphi$2 F0)[b].
Corollary $\varphi$2_PrPlus : $\forall$ F0 a b, Arithmetical_ultraFilter F0
  -> a $\in$ (Z' F0) -> b $\in$ (Z' F0)
  -> ($\varphi$2 F0)[(a +_F0 b)$\mbox{\color[RGB]{0,100,0}\%}$z'] = ($\varphi$2 F0)[a] +_F0 ($\varphi$2 F0)[b].
Corollary $\varphi$2_PrMult : $\forall$ F0 a b, Arithmetical_ultraFilter F0
  -> a $\in$ (Z' F0) -> b $\in$ (Z' F0)
  -> ($\varphi$2 F0)[(a $\cdot$_F0 b)$\mbox{\color[RGB]{0,100,0}\%}$z'] = ($\varphi$2 F0)[a] $\cdot$_F0 ($\varphi$2 F0)[b].
\end{lstlisting}

Given that ${^{*}{\mathbf Q}}_{^{*}{\mathbf Z}}$ is a subset of $^{*}{\mathbf Q}$, to embed $^{*}{\mathbf N}$ in $^{*}{\mathbf Q}$,
it suffices to embed it in ${^{*}{\mathbf Q}}_{^{*}{\mathbf Z}}$.
Reviewing the process from $^{*}{\mathbf N}$ to $^{*}{\mathbf Z}$, it is essentially an expansion by adding negative elements,
and the result of embedding $^{*}{\mathbf N}$ in $^{*}{\mathbf Z}$, denoted as ${^{*}{\mathbf Z}}_{^{*}{\mathbf N}}$,
is composed of all non-negative elements in $^{*}{\mathbf Z}$.
By analogy, when $^{*}{\mathbf N}$ is embedded in ${^{*}{\mathbf Q}}_{^{*}{\mathbf Z}}$,
or in other words, $^{*}{\mathbf N}$ is embedded in $^{*}{\mathbf Q}$, the result, denoted as ${^{*}{\mathbf Q}}_{^{*}{\mathbf N}}$,
is precisely composed of all non-negative elements in ${^{*}{\mathbf Q}}_{^{*}{\mathbf Z}}$:
$${^{*}{\mathbf Q}}_{^{*}{\mathbf N}} = \{u : u \in {^{*}{\mathbf Q}}_{^{*}{\mathbf Z}} \wedge ({^{*}}\mbox{Q}_0 = u \vee {^{*}}\mbox{Q}_0 <_{^{*}{\mathbf Q}} u)\}.$$
${^{*}{\mathbf Q}}_{^{*}{\mathbf N}}$ is a proper-subset of ${^{*}{\mathbf Q}}_{^{*}{\mathbf Z}}$.

\begin{lstlisting}
Definition Q'_N' F0 := \{ $\lambda$ u, u $\in$ (Q'_Z' F0)
  /\ (Q'0 F0 = u \/ (Q'0 F0) <_F0 u) \}.
Corollary Q'_N'_is_Set : $\forall$ F0, Arithmetical_ultraFilter F0 -> Ensemble (Q'_N' F0).
Corollary Q'_N'_properSubset_Q'_Z' : $\forall$ F0, Arithmetical_ultraFilter F0
  -> (Q'_N' F0) $\subset$ (Q'_Z' F0) /\ (Q'_N' F0) <> (Q'_Z' F0).
Corollary Q'_N'_properSubset_Q' : $\forall$ F0, Arithmetical_ultraFilter F0
  -> (Q'_N' F0) $\subset$ (Q' F0) /\ (Q'_N' F0) <> (Q' F0).
\end{lstlisting}

We have utilized $\varphi_1$ in \ref{sec5-2-5} to embed $^{*}{\mathbf N}$ in $^{*}{\mathbf Z}$,
and here $^{*}{\mathbf Z}$ is embedded in $^{*}{\mathbf Q}$ with $\varphi_2$,
thus $^{*}{\mathbf N}$ can be directly embedded in $^{*}{\mathbf Q}$ through composition of $\varphi_1$ and $\varphi_2$:
$$\varphi_3 = \varphi_2 \circ \varphi_1.$$
It can be verified that $\varphi_3$ is an order-preserved, addition-preserved and multiplication-preserved 1-1 function
between $^{*}{\mathbf N}$ and $^{*}{\mathbf Q}_{^{*}{\mathbf N}}$.

\begin{lstlisting}
Definition $\varphi$3 F0 := ($\varphi$2 F0) $\circ$ ($\varphi$1 F0).
Corollary $\varphi$3_is_Function : $\forall$ F0, Arithmetical_ultraFilter F0
  -> Function1_1 ($\varphi$3 F0) /\ dom($\varphi$3 F0) = N' F0 /\ ran($\varphi$3 F0) $\subset$ (Q' F0).

Corollary $\varphi$3_N'0 : $\forall$ F0, Arithmetical_ultraFilter F0 -> ($\varphi$3 F0)[F 0] = Q'0 F0.
Corollary $\varphi$3_N'1 : $\forall$ F0, Arithmetical_ultraFilter F0 -> ($\varphi$3 F0)[F 1] = Q'1 F0.
Corollary $\varphi$3_PrOrder : $\forall$ F0 M N, Arithmetical_ultraFilter F0 -> M $\in$ (N' F0)
  -> N $\in$ (N' F0) -> (M <_F0 N)$\mbox{\color[RGB]{0,100,0}\%}$n' <-> ($\varphi$3 F0)[M] <_F0 ($\varphi$3 F0)[N].
Corollary $\varphi$3_PrPlus : $\forall$ F0 M N, Arithmetical_ultraFilter F0
  -> M $\in$ (N' F0) -> N $\in$ (N' F0) -> ($\varphi$3 F0)[(M +_F0 N)$\mbox{\color[RGB]{0,100,0}\%}$n']
    = ($\varphi$3 F0)[M] +_F0 ($\varphi$3 F0)[N].
Corollary $\varphi$3_PrMult : $\forall$ F0 M N, Arithmetical_ultraFilter F0
  -> M $\in$ (N' F0) -> N $\in$ (N' F0) -> ($\varphi$3 F0)[(M $\cdot$_F0 N)$\mbox{\color[RGB]{0,100,0}\%}$n']
    = ($\varphi$3 F0)[M] $\cdot$_F0 ($\varphi$3 F0)[N].

Corollary $\varphi$3_ran : $\forall$ F0, Arithmetical_ultraFilter F0
  -> ran($\varphi$3 F0) = Q'_N' F0.
\end{lstlisting}

\subsubsection{Rational Numbers in $^{*}{\mathbf Q}$}\label{sec5-3-5}

In this section, we will find the structures that can serve as natural number set, integer set and rational number set contained in ${^{*}{\mathbf Q}}$.

We have already seen two sets -- $\omega$ \& ${^{*}\mathbf N}_{\mathbf N}$ -- that can serve as the natural number set
and either of them can be embedded into ${^{*}{\mathbf Q}}$ to obtain the structure of natural numbers contained in ${^{*}{\mathbf Q}}$.
The embedding of $\omega$ is more concise, because $\omega$ is isomorphically embedded into ${^{*}\mathbf N}$ with $\varphi$ in \ref{sec5-1-4},
${^{*}\mathbf N}$ into ${^{*}\mathbf Z}$ with $\varphi_1$ in \ref{sec5-2-5} and ${^{*}\mathbf Z}$ into ${^{*}\mathbf Q}$ with $\varphi_2$ in \ref{sec5-3-4},
then the composition of $\varphi$, $\varphi_1$ and $\varphi_2$ can directly isomorphically embed $\omega$ into ${^{*}{\mathbf Q}}$.
We denote this composition as $\varphi_4$:
$$\varphi_4 = \varphi_2 \circ (\varphi_1 \circ \varphi).$$
$\varphi_4$ is undoubtedly an order-preserved, addition-preserved and multiplication-preserved 1-1 function whose domain is $\omega$.

\begin{lstlisting}
Definition $\varphi$4 F0 := ($\varphi$2 F0) $\circ$ (($\varphi$1 F0) $\circ$ $\varphi$).
Corollary $\varphi$4_is_Function : $\forall$ F0, Arithmetical_ultraFilter F0
  -> Function1_1 ($\varphi$4 F0) /\ dom($\varphi$4 F0) = $\omega$ /\ ran($\varphi$4 F0) $\subset$ (Q' F0).

Corollary $\varphi$4_0 : $\forall$ F0, Arithmetical_ultraFilter F0 -> ($\varphi$4 F0)[0] = Q'0 F0.
Corollary $\varphi$4_1 : $\forall$ F0, Arithmetical_ultraFilter F0 -> ($\varphi$4 F0)[1] = Q'1 F0.
Corollary $\varphi$4_PrOrder : $\forall$ F0 m n, Arithmetical_ultraFilter F0
  -> m $\in$ $\omega$ -> n $\in$ $\omega$ -> m $\in$ n <-> ($\varphi$4 F0)[m] <_F0 ($\varphi$4 F0)[n].
Corollary $\varphi$4_PrPlus : $\forall$ F0 m n, Arithmetical_ultraFilter F0
  -> m $\in$ $\omega$ -> n $\in$ $\omega$ -> ($\varphi$4 F0)[m + n] = ($\varphi$4 F0)[m] +_F0 ($\varphi$4 F0)[n].
Corollary $\varphi$4_PrMult : $\forall$ F0 m n, Arithmetical_ultraFilter F0
  -> m $\in$ $\omega$ -> n $\in$ $\omega$ -> ($\varphi$4 F0)[m $\cdot$ n] = ($\varphi$4 F0)[m] $\cdot$_F0 ($\varphi$4 F0)[n].
\end{lstlisting}

\noindent Then the structure of natural number set contained in ${^{*}{\mathbf Q}}$ is defined as range of $\varphi_4$:
$${^{*}\mathbf Q}_{\mathbf N} = \mbox{range}\ \varphi_4,$$
it is a proper-subset of ${^{*}\mathbf Q}_{^{*}\mathbf N}$.

\begin{lstlisting}
Definition Q'_N F0 := ran($\varphi$4 F0).
Corollary Q'_N_Subset_Q' : $\forall$ F0, Arithmetical_ultraFilter F0
  -> (Q'_N F0) $\subset$ (Q' F0).
Corollary Q'_N_properSubset_Q'_N' : $\forall$ F0, Arithmetical_ultraFilter F0
  -> ($\forall$ n, F0 <> F n) -> (Q'_N F0) $\subset$ (Q'_N' F0) /\ (Q'_N F0) <> (Q'_N' F0).
\end{lstlisting}

Similarly, the structure of integer set contained in ${^{*}{\mathbf Q}}$ can be obtained through isomorphic embedding.
Since ${^{*}\mathbf Z}$ has already been embedded into ${^{*}{\mathbf Q}}$ with $\varphi_2$ in \ref{sec5-3-4},
${^{*}\mathbf Z}_{\mathbf Z}$, serving as the integer set, is naturally embedded into ${^{*}{\mathbf Q}}$ as a proper-subset of ${^{*}\mathbf Z}$.
The structure of integers contained in ${^{*}{\mathbf Q}}$ can be directly defined as the image set of $\varphi_2$ at ${^{*}\mathbf Z}_{\mathbf Z}$:
$${^{*}\mathbf Q}_{\mathbf Z} = \varphi_2\lceil{{^{*}\mathbf Z}_{\mathbf Z}}\rfloor,$$
it is a proper-subset of ${^{*}\mathbf Q}_{^{*}\mathbf Z}$.

\begin{lstlisting}
Definition Q'_Z F0 := ($\varphi$2 F0)$\lceil$Z'_Z F0$\rfloor$.
Corollary Q'_Z_Subset_Q' : $\forall$ F0, Arithmetical_ultraFilter F0
  -> (Q'_Z F0) $\subset$ (Q' F0).
Corollary Q'_Z_properSubset_Q'_Z' : $\forall$ F0, Arithmetical_ultraFilter F0
  -> ($\forall$ n, F0 <> F n) -> (Q'_Z F0) $\subset$ (Q'_Z' F0) /\ (Q'_Z F0) <> (Q'_Z' F0).
\end{lstlisting}

Meanwhile, it can be proven that every member in ${^{*}\mathbf Q}_{\mathbf Z}$ is structured as an equivalence class $[(a,{^{*}}$Z$_1)]_{R_{^{*}{\mathbf Z}}}$
where $a \in {^{*}\mathbf Z}_{\mathbf Z}$.
That is to say, ${^{*}\mathbf Q}_{\mathbf Z} = \{u : \exists a, a \in {^{*}\mathbf Z}_{\mathbf Z} \wedge u = [(a,{^{*}}$Z$_1)]_{R_{^{*}{\mathbf Z}}}\}$ holds.

\begin{lstlisting}
Corollary Q'_Z_Corollary : $\forall$ F0, Arithmetical_ultraFilter F0
  -> Q'_Z F0 = \{ $\lambda$ u, $\exists$ a, a $\in$ (Z'_Z F0) /\ u = \[[a,(Z'1 F0)]\]_F0 \}.
\end{lstlisting}

The relationship between ${^{*}\mathbf Q}_{\mathbf N}$ and ${^{*}\mathbf Q}_{\mathbf Z}$ must be in line with our intuitive
understanding of the relationship between natural numbers and integers.
Natural numbers are non-negative integers:
$${^{*}\mathbf Q}_{\mathbf N} = \{u : u \in {^{*}\mathbf Q}_{\mathbf Z} \wedge ({^{*}}\mbox{Q}_0 = u \vee {^{*}}\mbox{Q}_0 <_{^{*}\mathbf Q} u)\}.$$
And negative integers are essentially the negative values of all non-zero natural numbers:
$$\{u : u \in {^{*}\mathbf Q}_{\mathbf Z}\ \wedge\ u <_{^{*}\mathbf Q} {^{*}}\mbox{Q}_0\}
= \{u : u \in {^{*}\mathbf Q}\ \wedge\ (\exists v, v \in ({^{*}\mathbf Q}_{\mathbf N} \sim \{{^{*}}\mbox{Q}_0\})\ \wedge\ u = {^{*}}\mbox{Q}_0\ -_{^{*}\mathbf Q}\ v )\}.$$

\begin{lstlisting}
Fact Q'_N_Equ_Q'_Z_me_Q'0 : $\forall$ F0, Arithmetical_ultraFilter F0
  -> Q'_N F0 = \{ $\lambda$ u, u $\in$ (Q'_Z F0) /\ (Q'0 F0 = u \/ (Q'0 F0) <_F0 u) \}.
Fact Q'_N_neg_Equ_Q'_Z_lt_Q'0 : $\forall$ F0, Arithmetical_ultraFilter F0
  -> \{ $\lambda$ u, u $\in$ (Q'_Z F0) /\ u <_F0 (Q'0 F0) \}
    = \{ $\lambda$ u, u $\in$ (Q' F0) /\ $\exists$ v, v $\in$ ((Q'_N F0) $\sim$ [Q'0 F0])
    /\ u = (Q'0 F0) -_F0 v \}.
\end{lstlisting}

Obviously, ${^{*}\mathbf Q}_{\mathbf N}$ is a proper-subset of ${^{*}\mathbf Q}_{\mathbf Z}$.

\begin{lstlisting}
Fact Q'_N_properSubset_Q'_Z : $\forall$ F0, Arithmetical_ultraFilter F0
  -> (Q'_N F0) $\subset$ (Q'_Z F0) /\ (Q'_N F0) <> (Q'_Z F0).
\end{lstlisting}

As for the structure of rational numbers in ${^{*}\mathbf Q}$, we cannot acquire it through isomorphic embedding
because no such structures have been constructed in ${^{*}\mathbf N}$ or ${^{*}\mathbf Z}$.
Nevertheless, inspired by the process of ${^{*}\mathbf Z}$ to ${^{*}\mathbf Q}$ that
${^{*}\mathbf Q}$ is obtained by equivalently classifying the elements in ${^{*}\mathbf Z} \times ({^{*}\mathbf Z} \sim \{^{*}$Z$_0\})$,
we can restrict the classification to ${^{*}\mathbf Z}_{\mathbf Z} \times ({^{*}\mathbf Z}_{\mathbf Z} \sim \{^{*}$Z$_0\})$
and construct the desired structure.
That is to say, the structure of rational numbers in ${^{*}\mathbf Q}$ consists of
the equivalence classes shaped like $[(a,b)]_{R_{^{*}\mathbf Z}}$ where $(a,b) \in {^{*}\mathbf Z}_{\mathbf Z} \times ({^{*}\mathbf Z}_{\mathbf Z} \sim \{^{*}$Z$_0\})$.
$${^{*}\mathbf Q}_{\mathbf Q} = \{u : \exists\ a\ b,\ a \in {^{*}\mathbf Z}_{\mathbf Z} \wedge b \in ({^{*}\mathbf Z}_{\mathbf Z} \sim \{^{*}\mbox{Z}_0\})
\wedge u = [(a,b)]_{R_{^{*}\mathbf Z}}\}.$$

\begin{lstlisting}
Definition Q'_Q F0 := \{ $\lambda$ u, $\exists$ a b, a $\in$ (Z'_Z F0)
  /\ b $\in$ ((Z'_Z F0) $\sim$ [Z'0 F0]) /\ u = \[[a,b]\]_F0 \}.
Corollary Q'_Q_Subset_Q' : $\forall$ F0, Arithmetical_ultraFilter F0
  -> (Q'_Q F0) $\subset$ (Q' F0).
\end{lstlisting}

The relationship between ${^{*}\mathbf Q}_{\mathbf Z}$ and ${^{*}\mathbf Q}_{\mathbf Q}$ must be in line with our intuitive
understanding of the relationship between integers and rational numbers.
Dividing integers yields fractions (i.e., rational numbers):
$${^{*}\mathbf Q}_{\mathbf Q} = \{u : u \in {^{*}\mathbf Q} \wedge
(\exists\ a\ b,\ a \in {^{*}\mathbf Q}_{\mathbf Z} \wedge b \in ({^{*}\mathbf Q}_{\mathbf Z} \sim \{^{*}\mbox{Q}_0\}) \wedge u = a\ /_{^{*}\mathbf Q}\ b)\},$$
namely,
$$\forall u,v,\ (u,v) \in {^{*}\mathbf Q}_{\mathbf Z} \times ({^{*}\mathbf Q}_{\mathbf Z} \sim \{^{*}\mbox{Q}_0\})\ \Longrightarrow\ u\ /_{^{*}\mathbf Q}\ v\ \in\ {^{*}\mathbf Q}_{\mathbf Q}.$$

\begin{lstlisting}
Fact Q'_Q_Equ_Q'_Z_Div : $\forall$ F0, Arithmetical_ultraFilter F0
  -> Q'_Q F0 = \{ $\lambda$ u, u $\in$ (Q' F0) /\ $\exists$ a b, a $\in$ (Q'_Z F0)
    /\ b $\in$ ((Q'_Z F0) $\sim$ [Q'0 F0]) /\ u = a /_F0 b \}.
Fact Q'_Z_Divide_in_Q'_Q : $\forall$ F0 v u, Arithmetical_ultraFilter F0
  -> v $\in$ (Q'_Z F0) -> u $\in$ (Q'_Z F0) -> u <> Q'0 F0 -> (v /_F0 u) $\in$ (Q'_Q F0).
\end{lstlisting}

Similar to that integers are contained in rational numbers, ${^{*}\mathbf Q}_{\mathbf Z}$ must be a proper-subset of ${^{*}\mathbf Q}_{\mathbf Q}$.

\begin{lstlisting}
Fact Q'_Z_properSubset_Q'_Q : $\forall$ F0, Arithmetical_ultraFilter F0
  -> (Q'_Z F0) $\subset$ (Q'_Q F0) /\ (Q'_Z F0) <> (Q'_Q F0).
\end{lstlisting}

Besides, ${^{*}\mathbf Q}_{\mathbf N}$, ${^{*}\mathbf Q}_{\mathbf Z}$ and ${^{*}\mathbf Q}_{\mathbf Q}$ must all satisfy their expected algebraic properties.
Zero and identity element are their elements.

\begin{lstlisting}
Fact Q'0_in_Q'_N : $\forall$ F0, Arithmetical_ultraFilter F0 -> (Q'0 F0) $\in$ (Q'_N F0).
Fact Q'1_in_Q'_N : $\forall$ F0, Arithmetical_ultraFilter F0 -> (Q'1 F0) $\in$ (Q'_N F0).
Fact Q'0_in_Q'_Z : $\forall$ F0, Arithmetical_ultraFilter F0 -> (Q'0 F0) $\in$ (Q'_Z F0).
Fact Q'1_in_Q'_Z : $\forall$ F0, Arithmetical_ultraFilter F0 -> (Q'1 F0) $\in$ (Q'_Z F0).
Fact Q'0_in_Q'_Q : $\forall$ F0, Arithmetical_ultraFilter F0 -> (Q'0 F0) $\in$ (Q'_Q F0).
Fact Q'1_in_Q'_Q : $\forall$ F0, Arithmetical_ultraFilter F0 -> (Q'1 F0) $\in$ (Q'_Q F0).
\end{lstlisting}

\noindent Addition and multiplication are closed in all of them, subtraction in ${^{*}\mathbf Q}_{\mathbf Z}$ and ${^{*}\mathbf Q}_{\mathbf Q}$,
and division in ${^{*}\mathbf Q}_{\mathbf Q}$.

\begin{lstlisting}
Fact Q'_N_Plus_in_Q'_N : $\forall$ F0 u v, Arithmetical_ultraFilter F0
  -> u $\in$ (Q'_N F0) -> v $\in$ (Q'_N F0) -> (u +_F0 v) $\in$ (Q'_N F0).
Fact Q'_Z_Plus_in_Q'_Z : $\forall$ F0 u v, Arithmetical_ultraFilter F0
  -> u $\in$ (Q'_Z F0) -> v $\in$ (Q'_Z F0) -> (u +_F0 v) $\in$ (Q'_Z F0).
Fact Q'_Q_Plus_in_Q'_Q : $\forall$ F0 u v, Arithmetical_ultraFilter F0
  -> u $\in$ (Q'_Q F0) -> v $\in$ (Q'_Q F0) -> (u +_F0 v) $\in$ (Q'_Q F0).

Fact Q'_N_Mult_in_Q'_N : $\forall$ F0 u v, Arithmetical_ultraFilter F0
  -> u $\in$ (Q'_N F0) -> v $\in$ (Q'_N F0) -> (u $\cdot$_F0 v) $\in$ (Q'_N F0).
Fact Q'_Z_Mult_in_Q'_Z : $\forall$ F0 u v, Arithmetical_ultraFilter F0
  -> u $\in$ (Q'_Z F0) -> v $\in$ (Q'_Z F0) -> (u $\cdot$_F0 v) $\in$ (Q'_Z F0).
Fact Q'_Q_Mult_in_Q'_Q : $\forall$ F0 u v, Arithmetical_ultraFilter F0
  -> u $\in$ (Q'_Q F0) -> v $\in$ (Q'_Q F0) -> (u $\cdot$_F0 v) $\in$ (Q'_Q F0).

Fact Q'_Z_Minus_in_Q'_Z : $\forall$ F0 u v, Arithmetical_ultraFilter F0
  -> u $\in$ (Q'_Z F0) -> v $\in$ (Q'_Z F0) -> (u -_F0 v) $\in$ (Q'_Z F0).
Fact Q'_Q_Minus_in_Q'_Q : $\forall$ F0 u v, Arithmetical_ultraFilter F0
  -> u $\in$ (Q'_Q F0) -> v $\in$ (Q'_Q F0) -> (u -_F0 v) $\in$ (Q'_Q F0).

Fact Q'_Q_Divide_in_Q'_Q : $\forall$ F0 u v, Arithmetical_ultraFilter F0
  -> u $\in$ (Q'_Q F0) -> v $\in$ (Q'_Q F0) -> v <> Q'0 F0 -> (u /_F0 v) $\in$ (Q'_Q F0).
\end{lstlisting}

\subsubsection{Infinities and Infinitesimals in $^{*}{\mathbf Q}$}\label{sec5-3-6}

It can be observed that both $^{*}{\mathbf N}$ and $^{*}{\mathbf Z}$ introduce some infinity elements into general number sets,
and $^{*}{\mathbf Q}$, as the extension of $^{*}{\mathbf N}$ and $^{*}{\mathbf Z}$, naturally includes such infinities.
To clarify this point, the elements of $^{*}{\mathbf Q}$ can be divided into two categories.

Elements of one category, called finite elements of $^{*}{\mathbf Q}$, constitute a crucial subset of $^{*}{\mathbf Q}$:
$${^{*}\mathbf Q}_< = \{u : u \in {^{*}\mathbf Q} \wedge \exists\ k \in {^{*}\mathbf Q}_{\mathbf N}, |u| = k \vee |u| <_{^{*}\mathbf Q} k\}.$$
${^{*}\mathbf Q}_<$ is called the Archimedes Subset of ${^{*}\mathbf Q}$ because every of its elements $x$ satisfies Archimedean Property:
$$\exists\ k \in {^{*}\mathbf Q}_{\mathbf N}, |x| = k \vee |x| <_{^{*}\mathbf Q} k.$$
Namely, the absolute value of every element in ${^{*}\mathbf Q}$ is always less than a certain natural number,
which means measurability of elements in ${^{*}\mathbf Q}_<$.

\begin{lstlisting}
Definition Q'_$\nu$ F0 := \{ $\lambda$ u, u $\in$ (Q' F0)
  /\ $\exists$ k, k $\in$ (Q'_N F0) /\ (|u|_F0 = k \/ |u|_F0 <_F0 k) \}.
Corollary Q'_$\nu$_Subset_Q' : $\forall$ F0, (Q'_$\nu$ F0) $\subset$ (Q' F0).
Corollary Q'_$\nu$_is_Set : $\forall$ F0, F0 $\in$ $\beta\omega$ -> Ensemble (Q'_$\nu$ F0).
\end{lstlisting}

Undoubtedly, ${^{*}\mathbf Q}_{\mathbf N}$, ${^{*}\mathbf Q}_{\mathbf Z}$ and ${^{*}\mathbf Q}_{\mathbf Q}$ all should be subset of ${^{*}\mathbf Q}_<$.

\begin{lstlisting}
Fact Q'_N_properSubset_Q'_$\nu$ : $\forall$ F0, Arithmetical_ultraFilter F0
  -> (Q'_N F0) $\subset$ (Q'_$\nu$ F0) /\ (Q'_N F0) <> (Q'_$\nu$ F0).
Fact Q'_Z_properSubset_Q'_$\nu$ : $\forall$ F0, Arithmetical_ultraFilter F0
  -> (Q'_Z F0) $\subset$ (Q'_$\nu$ F0) /\ (Q'_Z F0) <> (Q'_$\nu$ F0).
Fact Q'_Q_Subset_Q'_$\nu$ : $\forall$ F0, Arithmetical_ultraFilter F0
  -> (Q'_Q F0) $\subset$ (Q'_$\nu$ F0).
\end{lstlisting}

For elements of the other category such as
$$\sigma,\ -\frac{3\sigma}{2},\ \frac{\sigma}{100}\ \cdots\ (\sigma \in ({^{*}{\mathbf Q}}_{^{*}{\mathbf Z}} \sim {^{*}{\mathbf Q}}_{\mathbf Z})),$$
they are not members of ${^{*}\mathbf Q}_<$ because their absolute values are larger than any natural number and certainly than all elements in ${^{*}\mathbf Q}_<$,
thus not satisfying Archimedean Property. They represent infinities in ${^{*}\mathbf Q}$ and constitute ${^{*}\mathbf Q} \sim {^{*}\mathbf Q}_<$.

\begin{lstlisting}
Fact Q'_$\nu$_lt_Q'Abs_others : $\forall$ F0 t, Arithmetical_ultraFilter F0
  -> ($\forall$ m, F0 <> (F m)) -> t $\in$ (Q' F0)
  -> t $\in$ ((Q' F0) $\sim$ (Q'_$\nu$ F0)) <-> ($\forall$ u, u $\in$ (Q'_$\nu$ F0) -> u <_F0 (|t|_F0)).
\end{lstlisting}

Concerning the operations on ${^{*}\mathbf Q}$, ${^{*}\mathbf Q}_<$ and ${^{*}\mathbf Q} \sim {^{*}\mathbf Q}_<$ exhibit the following properties:

1) Addition, multiplication and subtraction are all closed in ${^{*}\mathbf Q}_<$:
$$\forall u,v \in {^{*}\mathbf Q}_<,\ u +_{^{*}\mathbf Q} v,\ u \cdot_{^{*}\mathbf Q} v,\ u -_{^{*}\mathbf Q} v\ \in\ {^{*}\mathbf Q}_<.$$

2) An infinity plus a finite element in ${^{*}\mathbf Q}_<$ remains an infinity:
$$\forall u \in ({^{*}\mathbf Q} \sim {^{*}\mathbf Q}_<),\ v \in {^{*}\mathbf Q}_<,\ u +_{^{*}\mathbf Q} v\ \in\ ({^{*}\mathbf Q} \sim {^{*}\mathbf Q}_<).$$

3) Multiplication is closed in ${^{*}\mathbf Q} \sim {^{*}\mathbf Q}_<$:
$$\forall u,v \in ({^{*}\mathbf Q} \sim {^{*}\mathbf Q}_<),\ u \cdot_{^{*}\mathbf Q} v\ \in\ ({^{*}\mathbf Q} \sim {^{*}\mathbf Q}_<).$$

\begin{lstlisting}
Fact Q'_$\nu$_Plus_in_Q'_$\nu$ : $\forall$ F0 u v, Arithmetical_ultraFilter F0
  -> u $\in$ (Q'_$\nu$ F0) -> v $\in$ (Q'_$\nu$ F0) -> (u +_F0 v) $\in$ (Q'_$\nu$ F0).
Fact Q'_$\nu$_Mult_in_Q'_$\nu$ : $\forall$ F0 u v, Arithmetical_ultraFilter F0
  -> u $\in$ (Q'_$\nu$ F0) -> v $\in$ (Q'_$\nu$ F0) -> (u $\cdot$_F0 v) $\in$ (Q'_$\nu$ F0).
Fact Q'_$\nu$_Minus_in_Q'_$\nu$ : $\forall$ F0 u v, Arithmetical_ultraFilter F0
  -> u $\in$ (Q'_$\nu$ F0) -> v $\in$ (Q'_$\nu$ F0) -> (u -_F0 v) $\in$ (Q'_$\nu$ F0).

Fact infinity_Plus_finity : $\forall$ F0 u v, Arithmetical_ultraFilter F0
  -> ($\forall$ m, F0 <> F m) -> u $\in$ ((Q' F0) $\sim$ (Q'_$\nu$ F0))
  -> v $\in$ (Q'_$\nu$ F0) -> (u +_F0 v) $\in$ ((Q' F0) $\sim$ (Q'_$\nu$ F0)).
Fact infinity_Mult_infinity : $\forall$ F0 u v, Arithmetical_ultraFilter F0
  -> ($\forall$ m, F0 <> F m) -> u $\in$ ((Q' F0) $\sim$ (Q'_$\nu$ F0))
  -> v $\in$ ((Q' F0) $\sim$ (Q'_$\nu$ F0)) -> (u $\cdot$_F0 v) $\in$ ((Q' F0) $\sim$ (Q'_$\nu$ F0)).
\end{lstlisting}

Note that division is not closed in ${^{*}\mathbf Q}_<$ due to the existence of a special class of elements in ${^{*}\mathbf Q}_<$ -- infinitesimals.
For example, $^{*}$Q$_1\ /_{^{*}\mathbf Q}\ \sigma\ (\sigma \in ({^{*}\mathbf Q} \sim {^{*}\mathbf Q}_<))$ is an infinitesimal
and $^{*}$Q$_1\ /_{^{*}\mathbf Q}\ (^{*}$Q$_1\ /_{^{*}\mathbf Q}\ \sigma)$ is $\sigma\ (\notin {^{*}\mathbf Q}_<)$.

The absolute value of an infinitesimal is less than the inverse of any non-zero natural number
and the set consisting of infinitesimals is denoted as:
$${\mathbf I} = \{u : u \in {^{*}\mathbf Q}
\wedge (\forall k \in ({^{*}\mathbf Q}_{\mathbf N} \sim \{^{*}\mbox{Q}_0\}),\ |u|\ <_{^{*}\mathbf Q}\ ({^{*}\mbox{Q}_1}\ /_{^{*}\mathbf Q}\ k))\}.$$

\begin{lstlisting}
Definition I F0 := \{ $\lambda$ u, u $\in$ (Q' F0) /\ ($\forall$ k, k $\in$ ((Q'_N F0) $\sim$ [Q'0 F0])
  -> (|u|_F0) <_F0 ((Q'1 F0) /_F0 k)) \}.
Corollary I_Subset_Q' : $\forall$ F0, (I F0) $\subset$ (Q' F0).
Corollary I_is_Set : $\forall$ F0, F0 $\in$ $\beta\omega$ -> Ensemble (I F0).
Corollary I_properSubset_Q'_$\nu$ : $\forall$ F0, Arithmetical_ultraFilter F0
  -> (I F0) $\subset$ (Q'_$\nu$ F0) /\ (I F0) <> (Q'_$\nu$ F0).
\end{lstlisting}

It is obvious that the zero element $^{*}$Q$_0$ is an infinitesimal and $\{^{*}$Q$_0\}$ is a proper-subset of $\mathbf I$.

\begin{lstlisting}
Fact Q'0_in_I : $\forall$ F0, Arithmetical_ultraFilter F0 -> (Q'0 F0) $\in$ (I F0).
Fact Q'0_Singleton_properSubset_I : $\forall$ F0, Arithmetical_ultraFilter F0
  -> ($\forall$ n, F0 <> F n) -> [Q'0 F0] $\subset$ (I F0) /\ [Q'0 F0] <> (I F0).
\end{lstlisting}

Concerning the operations on ${^{*}\mathbf Q}$, $\mathbf I$ exhibit the following properties:

1) Addition is closed in $\mathbf I$:
$$\forall u,v \in {\mathbf I},\ u +_{^{*}\mathbf Q} v\ \in\ {\mathbf I}.$$

2) An infinitesimal times a finite element in ${^{*}\mathbf Q}_<$ remains an infinitesimal:
$$\forall u \in {\mathbf I},\ v \in {^{*}\mathbf Q}_<,\ u \cdot_{^{*}\mathbf Q} v\ \in\ {\mathbf I}.$$

3) Inverse of a non-zero infinitesimal is an infinity and inverse of an infinity is a non-zero infinitesimal:
$$\forall u, v \in {^{*}\mathbf Q},\ if\ u \cdot_{^{*}\mathbf Q} v = {^{*}}\mbox{Q}_1\ then\ u \in ({\mathbf I} \sim \{{^{*}}\mbox{Q}_0\})
\Longleftrightarrow v \in ({^{*}\mathbf Q} \sim {^{*}\mathbf Q}_<).$$

\begin{lstlisting}
Fact I_Plus_in_I : $\forall$ F0 u v, Arithmetical_ultraFilter F0
  -> u $\in$ (I F0) -> v $\in$ (I F0) -> (u +_F0 v) $\in$ (I F0).
Fact I_Mult_in_I : $\forall$ F0 u v, Arithmetical_ultraFilter F0
  -> u $\in$ (I F0) -> v $\in$ (Q'_$\nu$ F0) -> (u $\cdot$_F0 v) $\in$ (I F0).
Fact I_inv_Property1 : $\forall$ F0 u v, Arithmetical_ultraFilter F0
  -> ($\forall$ n, F0 <> F n) -> u $\in$ ((I F0) $\sim$ [Q'0 F0]) -> v $\in$ (Q' F0)
  -> u $\cdot$_F0 v = Q'1 F0 -> v $\in$ ((Q' F0) $\sim$ (Q'_$\nu$ F0)).
Fact I_inv_Property2 : $\forall$ F0 u v, Arithmetical_ultraFilter F0
  -> ($\forall$ n, F0 <> F n) -> u $\in$ ((Q' F0) $\sim$ (Q'_$\nu$ F0)) -> v $\in$ (Q' F0)
  -> u $\cdot$_F0 v = Q'1 F0 -> v $\in$ ((I F0) $\sim$ [Q'0 F0]).
\end{lstlisting}

It is obvious that the structure of ${^{*}\mathbf Q}$ is much more complicated than general rational number set and even than real number set.
Including both infinity and infinitesimal numbers, ${^{*}\mathbf Q}$ is a set of non-standard rational numbers.

\subsection{Real Number Set $\mathbf R$}\label{sec5-4}

\subsubsection{Construction of $\mathbf R$}\label{sec5-4-1}

Given the intricate ${^{*}\mathbf Q}$ structure, real number set $\mathbf R$ can be cleverly derived from it.
Specifically, $\mathbf R$ is constructed by equivalently classifying the elements of ${^{*}\mathbf Q}_<$ based on the following equivalence relation:
$$R_{^{*}{\mathbf Q}} = \{(u,v) : (u,v) \in ({{^{*}\mathbf Q}_<} \times {{^{*}\mathbf Q}_<})\ \wedge\ u -_{^{*}\mathbf Q} v\ \in\ {\mathbf I}\}.$$

\begin{lstlisting}
Definition R_Q' F0 := \{\ $\lambda$ u v, u $\in$ (Q'_$\nu$ F0)
  /\ v $\in$ (Q'_$\nu$ F0) /\ (u -_F0 v) $\in$ (I F0) \}\.
Corollary R_Q'_is_equRelation : $\forall$ F0, Arithmetical_ultraFilter F0
  -> equRelation (R_Q' F0) (Q'_$\nu$ F0).
\end{lstlisting}

According to this equivalence relation, the difference between two elements in the same equivalence class is an infinitesimal.
Thus the relationship between two identical equivalence classes regarding $R_{^{*}{\mathbf Q}}$ is:
$$\forall u,v \in {{^{*}\mathbf Q}_<},\ [u]_{R_{^{*}{\mathbf Q}}} = [v]_{R_{^{*}{\mathbf Q}}}\ \Longleftrightarrow\
u -_{^{*}{\mathbf Q}} v\ \in\ {\mathbf I}.$$

\begin{lstlisting}
Notation "\[ u \]_ F0 " := (equClass u (R_Q' F0) (Q'_$\nu$ F0))
  (at level 5, u at level 0, F0 at level 0) : r_scope.

Corollary R_Q'_Corollary : $\forall$ F0 u v, Arithmetical_ultraFilter F0
  -> u $\in$ (Q'_$\nu$ F0) -> v $\in$ (Q'_$\nu$ F0)
  -> \[u\]_F0 = \[v\]_F0 <-> (u -_F0 v) $\in$ (I F0).
\end{lstlisting}

$\mathbf R$ is the quotient set of ${^{*}\mathbf Q}_<$ regarding $R_{^{*}{\mathbf Q}}$:
$${\mathbf R} = {{^{*}\mathbf Q}_<} / {R_{^{*}{\mathbf Q}}} = \{u : \exists x, u = [x]_{R_{^{*}{\mathbf Q}}}\}.$$

\begin{lstlisting}
Definition $\mathrm R$ F0 := (Q'_$\nu$ F0) / (R_Q' F0).
Corollary R_is_Set : $\forall$ F0, F0 $\in$ $\beta\omega$ -> Ensemble ($\mathrm R$ F0).

Corollary R_Corollary1 : $\forall$ F0 u, Arithmetical_ultraFilter F0
  -> u $\in$ ($\mathrm R$ F0) -> $\exists$ x, x $\in$ (Q'_$\nu$ F0) /\ x $\in$ (Q' F0) /\ u = \[x\]_F0.
Corollary R_Corollary2 : $\forall$ F0 x, Arithmetical_ultraFilter F0
  -> x $\in$ (Q'_$\nu$ F0) -> \[x\]_F0 $\in$ ($\mathrm R$ F0).
\end{lstlisting}

The zero element R$_0$ and identity element R$_1$ are respectively the equivalence class represented by $^{*}$Q$_0$ and $^{*}$Q$_1$:
$${\mbox{R}}_0 = [{^{*}\mbox{Q}}_0]_{R_{^{*}{\mathbf Q}}},$$
$${\mbox{R}}_1 = [{^{*}\mbox{Q}}_1]_{R_{^{*}{\mathbf Q}}}.$$

\begin{lstlisting}
Definition R0 F0 := \[Q'0 F0\]_F0.
Corollary R0_in_R : $\forall$ F0, Arithmetical_ultraFilter F0 -> (R0 F0) $\in$ ($\mathrm R$ F0).

Definition R1 F0 := \[Q'1 F0\]_F0.
Corollary R1_in_R : $\forall$ F0, Arithmetical_ultraFilter F0 -> (R1 F0) $\in$ ($\mathrm R$ F0).

Fact R0_isn't_R1 : $\forall$ F0, Arithmetical_ultraFilter F0 -> R0 F0 <> R1 F0.
\end{lstlisting}

Especially, given that R$_0$ is the equivalence class represented by ${^{*}}$Q$_0$,
the difference between every member of R$_0$ and ${^{*}}$Q$_0$ (i.e., this member itself) is an infinitesimal.
Therefore, R$_0$ is exactly the infinitesimal set $\mathbf I$.

\begin{lstlisting}
Corollary R0_Corollary : $\forall$ F0, Arithmetical_ultraFilter F0 -> R0 F0 = I F0.
\end{lstlisting}

\subsubsection{Order on $\mathbf R$}\label{sec5-4-2}

\begin{definition}[Order on ${\mathbf R}$]
Let $u,v$ be members in ${\mathbf R}$, and $u = [x]_{R_{^{*}{\mathbf Q}}}, v = [y]_{R_{^{*}{\mathbf Q}}}$
($x,y \in {^{*}{\mathbf Q}_<}$), then the order of $u$ and $v$ is defined as follows:
$$u <_{\mathbf R} v\ \Longleftrightarrow\
[x]_{R_{^{*}{\mathbf Q}}} <_{\mathbf R} [y]_{R_{^{*}{\mathbf Q}}}\ \Longleftrightarrow\
x <_{^{*}{\mathbf Q}} y\ \wedge\ [x]_{R_{^{*}{\mathbf Q}}} \neq [y]_{R_{^{*}{\mathbf Q}}}.$$
\end{definition}

The definition of this order relation is reasonable,
because for any $x,x_1,y,y_1, \in {^{*}{\mathbf Q}_<}$, if
$$u = [x]_{R_{^{*}{\mathbf Q}}} = [x_1]_{R_{^{*}{\mathbf Q}}},\
v = [y]_{R_{^{*}{\mathbf Q}}} = [y_1]_{R_{^{*}{\mathbf Q}}},$$
then it can be verified that
$$u <_{\mathbf R} v\ \Longleftrightarrow\
x <_{^{*}{\mathbf Q}} y \Longleftrightarrow\
x_1 <_{^{*}{\mathbf Q}} y_1.$$
Namely, the order relation on ${\mathbf R}$ is independent of the selection of representations of equivalence classes.

\begin{lstlisting}
Definition R_Ord F0 := \{\ $\lambda$ u v, u <> v /\ $\forall$ x y, x $\in$ (Q'_$\nu$ F0)
  -> y $\in$ (Q'_$\nu$ F0) -> u = \[x\]_F0 -> v = \[y\]_F0 -> x <_F0 y \}\.
Notation "u <_ F0 v" := (Rrelation u (R_Ord F0) v)
  (at level 10, F0 at level 9) : r_scope.

Corollary R_Ord_reasonable : $\forall$ F0 x x1 y y1, Arithmetical_ultraFilter F0
  -> x $\in$ (Q'_$\nu$ F0) -> x1 $\in$ (Q'_$\nu$ F0) -> y $\in$ (Q'_$\nu$ F0) -> y1 $\in$ (Q'_$\nu$ F0)
  -> (\[x\]_F0 = \[x1\]_F0)$\mbox{\color[RGB]{0,100,0}\%}$r -> (\[y\]_F0 = \[y1\]_F0)$\mbox{\color[RGB]{0,100,0}\%}$r
  -> (\[x\]_F0 <> \[y\]_F0)$\mbox{\color[RGB]{0,100,0}\%}$r -> x <_F0 y <-> x1 <_F0 y1.
Corollary R_Ord_Corollary : $\forall$ F0 x y, Arithmetical_ultraFilter F0
  -> x $\in$ (Q'_$\nu$ F0) -> y $\in$ (Q'_$\nu$ F0) -> \[x\]_F0 <_F0 \[y\]_F0
    <-> ((x <_F0 y)$\mbox{\color[RGB]{0,100,0}\%}$q' /\ \[x\]_F0 <> \[y\]_F0).
\end{lstlisting}

The properties of this order are entirely the same as those on $^{*}{\mathbf N}$, $^{*}{\mathbf Z}$ and $^{*}{\mathbf Q}$,
aligning with the conventional understanding of the order among real numbers.

\begin{lstlisting}
Fact R_Ord_irReflex : $\forall$ F0 u v, Arithmetical_ultraFilter F0 -> u = v
  -> $\sim$ u <_F0 v.
Fact R_Ord_Trans : $\forall$ F0 u v w, Arithmetical_ultraFilter F0
  -> u $\in$ ($\mathrm R$ F0) -> v $\in$ ($\mathrm R$ F0) -> w $\in$ ($\mathrm R$ F0)
  -> u <_F0 v -> v <_F0 w -> u <_F0 w.
Fact R_Ord_Tri : $\forall$ F0, Arithmetical_ultraFilter F0 -> Connect (R_Ord F0) ($\mathrm R$ F0).
Fact R0_lt_R1 : $\forall$ F0, Arithmetical_ultraFilter F0 -> (R0 F0) <_F0 (R1 F0).
\end{lstlisting}

\subsubsection{Operations on $\mathbf R$}\label{sec5-4-3}

\begin{definition}[Addition on ${\mathbf R}$]
Let $u,v$ be members in ${\mathbf R}$, and $u = [x]_{R_{^{*}{\mathbf Q}}}, v = [y]_{R_{^{*}{\mathbf Q}}}$
($x,y \in {^{*}{\mathbf Q}_<}$), then the sum of $u$ and $v$ is defined as follows:
$$u +_{\mathbf R} v\ =\ [x]_{R_{^{*}{\mathbf Q}}}\ +_{\mathbf R}\ [y]_{R_{^{*}{\mathbf Q}}}\
=\ [x +_{^{*}{\mathbf Q}} y]_{R_{^{*}{\mathbf Q}}}.$$
\end{definition}

The definition of this addition is reasonable,
because for any $x,x_1,y,y_1$ in ${^{*}{\mathbf Q}_<}$, if
$$u = [x]_{R_{^{*}{\mathbf Q}}} = [x_1]_{R_{^{*}{\mathbf Q}}},\
v = [y]_{R_{^{*}{\mathbf Q}}} = [y_1]_{R_{^{*}{\mathbf Q}}},$$
then it can be verified that
$$u +_{\mathbf R} v\ =\ [x]_{R_{^{*}{\mathbf Q}}}\ +_{\mathbf R}\ [y]_{R_{^{*}{\mathbf Q}}}\ =\
[x_1]_{R_{^{*}{\mathbf Q}}}\ +_{\mathbf R}\ [y_1]_{R_{^{*}{\mathbf Q}}}.$$
Namely, the result of $u +_{\mathbf R} v$ is unique and independent of the selection of representations of equivalence classes.

\begin{lstlisting}
Definition R_Plus F0 u v := $\cap$(\{ $\lambda$ w, $\forall$ x y, x $\in$ (Q'_$\nu$ F0)
  -> y $\in$ (Q'_$\nu$ F0) -> u = \[x\]_F0 -> v = \[y\]_F0
  -> w = \[(x +_F0 y)$\mbox{\color[RGB]{0,100,0}\%}$q'\]_F0 \}).
Notation "u +_ F0 v" := (R_Plus F0 u v)
  (at level 10, F0 at level 9) : r_scope.

Corollary R_Plus_reasonable : $\forall$ F0 u v, Arithmetical_ultraFilter F0
  -> u $\in$ ($\mathrm R$ F0) -> v $\in$ ($\mathrm R$ F0) -> $\exists$ w, w $\in$ ($\mathrm R$ F0) /\ [w]
    = \{ $\lambda$ w, $\forall$ x y, x $\in$ (Q'_$\nu$ F0) -> y $\in$ (Q'_$\nu$ F0)
    -> u = \[x\]_F0 -> v = \[y\]_F0  -> w = \[(x +_F0 y)$\mbox{\color[RGB]{0,100,0}\%}$q'\]_F0 \}.
Corollary Property_R_Plus : $\forall$ F0 x y, Arithmetical_ultraFilter F0
  -> x $\in$ (Q'_$\nu$ F0) -> y $\in$ (Q'_$\nu$ F0)
  -> (\[x\]_F0) +_F0 (\[y\]_F0) = \[(x +_F0 y)$\mbox{\color[RGB]{0,100,0}\%}$q'\]_F0.
\end{lstlisting}

\begin{definition}[Multiplication on ${\mathbf R}$]
Let $u,v$ be members in ${\mathbf R}$, and $u = [x]_{R_{^{*}{\mathbf Q}}}, v = [y]_{R_{^{*}{\mathbf Q}}}$
($x,y \in {^{*}{\mathbf Q}_<}$), then the product of $u$ and $v$ is defined as follows:
$$u \cdot_{\mathbf R} v\ =\ [x]_{R_{^{*}{\mathbf Q}}}\ \cdot_{\mathbf R}\ [y]_{R_{^{*}{\mathbf Q}}}\
=\ [x \cdot_{^{*}{\mathbf Q}} y]_{R_{^{*}{\mathbf Q}}}.$$
\end{definition}

The definition of this multiplication is reasonable,
because for any $x,x_1,y,y_1$ in ${^{*}{\mathbf Q}_<}$, if
$$u = [x]_{R_{^{*}{\mathbf Q}}} = [x_1]_{R_{^{*}{\mathbf Q}}},\
v = [y]_{R_{^{*}{\mathbf Q}}} = [y_1]_{R_{^{*}{\mathbf Q}}},$$
then it can be verified that
$$u \cdot_{\mathbf R} v\ =\ [x]_{R_{^{*}{\mathbf Q}}}\ \cdot_{\mathbf R}\ [y]_{R_{^{*}{\mathbf Q}}}\ =\
[x_1]_{R_{^{*}{\mathbf Q}}}\ \cdot_{\mathbf R}\ [y_1]_{R_{^{*}{\mathbf Q}}}.$$
Namely, the result of $u \cdot_{\mathbf R} v$ is unique and independent of the selection of representations of equivalence classes.

\begin{lstlisting}
Definition R_Mult F0 u v := $\cap$(\{ $\lambda$ w, $\forall$ x y, x $\in$ (Q'_$\nu$ F0)
  -> y $\in$ (Q'_$\nu$ F0) -> u = \[x\]_F0 -> v = \[y\]_F0
  -> w = \[(x $\cdot$_F0 y)$\mbox{\color[RGB]{0,100,0}\%}$q'\]_F0 \}).
Notation "u $\cdot$_ F0 v" := (R_Mult F0 u v)
  (at level 10, F0 at level 9) : r_scope.

Corollary R_Mult_reasonable : $\forall$ F0 u v, Arithmetical_ultraFilter F0
  -> u $\in$ ($\mathrm R$ F0) -> v $\in$ ($\mathrm R$ F0) -> $\exists$ w, w $\in$ ($\mathrm R$ F0) /\ [w]
    = \{ $\lambda$ w, $\forall$ x y, x $\in$ (Q'_$\nu$ F0) -> y $\in$ (Q'_$\nu$ F0)
    -> u = \[x\]_F0 -> v = \[y\]_F0 -> w = \[(x $\cdot$_F0 y)$\mbox{\color[RGB]{0,100,0}\%}$q'\]_F0 \}.
Corollary R_Mult_Corollary : $\forall$ F0 x y, Arithmetical_ultraFilter F0
  -> x $\in$ (Q'_$\nu$ F0) -> y $\in$ (Q'_$\nu$ F0)
  -> \[x\]_F0 $\cdot$_F0 \[y\]_F0 = \[(x $\cdot$_F0 y)$\mbox{\color[RGB]{0,100,0}\%}$q'\]_F0.
\end{lstlisting}

\begin{definition}[Subtraction on ${\mathbf R}$]
Let $u,v,w$ be members in ${\mathbf R}$, and $u +_{\mathbf R} w = v$, then $v$ minus $u$, denoted as $v -_{\mathbf R} u$, is $w$.
\end{definition}

The definition of this subtraction is reasonable,
because it can be verified that the existence of the very $w$ making $u +_{\mathbf R} w = v$ hold is unique.

\begin{lstlisting}
Definition R_Minus F0 v u := $\cap$(\{ $\lambda$ w, w $\in$ ($\mathrm R$ F0) /\ u +_F0 w = v \}).
Notation "u -_ F0 v" := (R_Minus F0 u v)(at level 10, F0 at level 9) : r_scope.

Corollary R_Minus_reasonable : $\forall$ F0 u v, Arithmetical_ultraFilter F0
  -> u $\in$ ($\mathrm R$ F0) -> v $\in$ ($\mathrm R$ F0) -> $\exists$! u0, u0 $\in$ ($\mathrm R$ F0) /\ u +_F0 u0 = v.
Corollary R_Minus_Corollary : $\forall$ F0 u w v, Arithmetical_ultraFilter F0
  -> u $\in$ ($\mathrm R$ F0) -> w $\in$ ($\mathrm R$ F0) -> v $\in$ ($\mathrm R$ F0)
  -> u +_F0 w = v <-> v -_F0 u = w.
\end{lstlisting}

\begin{definition}[Division on ${\mathbf R}$]
Let $u,v,w$ be members in ${\mathbf R}$, $u \neq$ {\rm R}$_0$, and $u \cdot_{\mathbf R} w = v$, then $v$ divided by $u$, denoted as $v /_{\mathbf R} u$, is $w$.
\end{definition}

The definition of this division is reasonable,
because it can be verified that the existence of the very $w$ making $u /_{\mathbf R} w = v$ hold is unique.

\begin{lstlisting}
Definition R_Divide F0 v u := $\cap$(\{ $\lambda$ w, w $\in$ ($\mathrm R$ F0)
  /\ u <> R0 F0 /\ u $\cdot$_F0 w = v \}).
Notation "v /_ F0 u" := (R_Divide F0 v u)(at level 10, F0 at level 9) : r_scope.

Corollary R_Divide_reasonable : $\forall$ F0 u v, Arithmetical_ultraFilter F0
  -> u $\in$ ($\mathrm R$ F0) -> v $\in$ ($\mathrm R$ F0) -> u <> R0 F0
  -> $\exists$! w, w $\in$ ($\mathrm R$ F0) /\ u $\cdot$_F0 w = v.
Corollary R_Divide_Corollary : $\forall$ F0 u w v, Arithmetical_ultraFilter F0
  -> u $\in$ ($\mathrm R$ F0) -> u <> R0 F0 -> w $\in$ ($\mathrm R$ F0) -> v $\in$ ($\mathrm R$ F0)
  -> u $\cdot$_F0 w = v <-> v /_F0 u = w.
\end{lstlisting}

${\mathbf R}$ possesses entirely the same operational properties as those on ${^{*}}{\mathbf Q}$,
and they are in line with our conventional understanding about real numbers.

Firstly, these operations are closed in ${\mathbf R}$.

\begin{lstlisting}
Fact R_Plus_in_R : $\forall$ F0 u v, Arithmetical_ultraFilter F0
  -> u $\in$ ($\mathrm R$ F0) -> v $\in$ ($\mathrm R$ F0) -> (u +_F0 v) $\in$ ($\mathrm R$ F0).
Fact R_Mult_in_R : $\forall$ F0 u v, Arithmetical_ultraFilter F0
  -> u $\in$ ($\mathrm R$ F0) -> v $\in$ ($\mathrm R$ F0) -> (u $\cdot$_F0 v) $\in$ ($\mathrm R$ F0).
Fact R_Minus_in_R : $\forall$ F0 v u, Arithmetical_ultraFilter F0
  -> v $\in$ ($\mathrm R$ F0) -> u $\in$ ($\mathrm R$ F0) -> (v -_F0 u) $\in$ ($\mathrm R$ F0).
Fact R_Divide_in_R : $\forall$ F0 v u, Arithmetical_ultraFilter F0
  -> v $\in$ ($\mathrm R$ F0) -> u $\in$ ($\mathrm R$ F0) -> u <> R0 F0 -> (v /_F0 u) $\in$ ($\mathrm R$ F0).
\end{lstlisting}

Secondly, the zero element R$_0$ and identity element R$_1$ satisfy: for every $u \in {\mathbf R}$,
$u\ +_{\mathbf R}$ R$_0 = u$, $u\ \cdot_{\mathbf R}$ R$_0 =$ R$_0$, $u\ \cdot_{\mathbf R}$ R$_1 = u$,
$u\ -_{\mathbf R} u =$ R$_0$, $u\ -_{\mathbf R}$ R$_0 = u$, $u\ /_{\mathbf R}$ R$_1 = u$
and $u\ /_{\mathbf R} u =$ R$_1$ if $u \neq$ R$_0$. And there exist no zero divisors in $\mathbf R$;
in other words, if $u \cdot_{\mathbf R} v =$ R$_0$ then at least one of $u$ and $v$ is R$_0$.

\begin{lstlisting}
Fact R_Plus_Property : $\forall$ F0 u, Arithmetical_ultraFilter F0
  -> u $\in$ ($\mathrm R$ F0) -> u +_F0 (R0 F0) = u.
Fact R_Mult_Property1 : $\forall$ F0 u, Arithmetical_ultraFilter F0
  -> u $\in$ ($\mathrm R$ F0) -> u $\cdot$_F0 (R0 F0) = R0 F0.
Fact R_Mult_Property2 : $\forall$ F0 u, Arithmetical_ultraFilter F0
  -> u $\in$ ($\mathrm R$ F0) -> u $\cdot$_F0 (R1 F0) = u.
Fact R_Mult_Property3 : $\forall$ F0 u v, Arithmetical_ultraFilter F0
  -> u $\in$ ($\mathrm R$ F0) -> v $\in$ ($\mathrm R$ F0) -> u $\cdot$_F0 v = (R0 F0)
  -> u = R0 F0 \/ v = R0 F0.
Fact R_Minus_Property1 : $\forall$F0 u, Arithmetical_ultraFilter F0
  -> u $\in$ ($\mathrm R$ F0) -> u -_F0 u = R0 F0.
Fact R_Minus_Property2 : $\forall$F0 u, Arithmetical_ultraFilter F0
  -> u $\in$ ($\mathrm R$ F0) -> u -_F0 (R0 F0) = u.
Fact R_Divide_Property1 : $\forall$F0 u, Arithmetical_ultraFilter F0
  -> u $\in$ ($\mathrm R$ F0) -> u <> R0 F0 -> u /_F0 u = R1 F0.
Fact R_Divide_Property2 : $\forall$F0 u, Arithmetical_ultraFilter F0
  -> u $\in$ ($\mathrm R$ F0) -> u /_F0 (R1 F0) = u.
\end{lstlisting}

Thirdly, these operations respectively satisfy their own expected operation laws such as
commutative law, associative law, cancellation law, distributive law, etc.

\begin{lstlisting}
Fact R_Plus_Commutation : $\forall$ F0 u v, Arithmetical_ultraFilter F0
  -> u $\in$ ($\mathrm R$ F0) -> v $\in$ ($\mathrm R$ F0) -> u +_F0 v = v +_F0 u.
Fact R_Plus_Association : $\forall$ F0 u v w, Arithmetical_ultraFilter F0
  -> u $\in$ ($\mathrm R$ F0) -> v $\in$ ($\mathrm R$ F0) -> w $\in$ ($\mathrm R$ F0)
  -> (u +_F0 v) +_F0 w = u +_F0 (v +_F0 w).
Fact R_Plus_Cancellation : $\forall$ F0 u v w, Arithmetical_ultraFilter F0
  -> u $\in$ ($\mathrm R$ F0) -> v $\in$ ($\mathrm R$ F0) -> w $\in$ ($\mathrm R$ F0)
  -> u +_F0 v = u +_F0 w -> v = w.
Fact R_Plus_PrOrder : $\forall$ F0 u v w, Arithmetical_ultraFilter F0
  -> u $\in$ ($\mathrm R$ F0) -> v $\in$ ($\mathrm R$ F0) -> w $\in$ ($\mathrm R$ F0)
  -> u <_F0 v <-> (w +_F0 u) <_F0 (w +_F0 v).
Fact R_neg_Property : $\forall$ F0 u, Arithmetical_ultraFilter F0
  -> u $\in$ ($\mathrm R$ F0) -> $\exists$! v, v $\in$ ($\mathrm R$ F0) /\ u +_F0 v = R0 F0.

Fact R_Mult_Commutation : $\forall$ F0 u v, Arithmetical_ultraFilter F0
  -> u $\in$ ($\mathrm R$ F0) -> v $\in$ ($\mathrm R$ F0) -> u $\cdot$_F0 v = v $\cdot$_F0 u.
Fact R_Mult_Association : $\forall$ F0 u v w, Arithmetical_ultraFilter F0
  -> u $\in$ ($\mathrm R$ F0) -> v $\in$ ($\mathrm R$ F0) -> w $\in$ ($\mathrm R$ F0)
  -> (u $\cdot$_F0 v) $\cdot$_F0 w = u $\cdot$_F0 (v $\cdot$_F0 w).
Fact R_Mult_Distribution : $\forall$ F0 u v w, Arithmetical_ultraFilter F0
  -> u $\in$ ($\mathrm R$ F0) -> v $\in$ ($\mathrm R$ F0) -> w $\in$ ($\mathrm R$ F0)
  -> u $\cdot$_F0 (v +_F0 w) = (u $\cdot$_F0 v) +_F0 (u $\cdot$_F0 w).
Fact R_Mult_Cancellation : $\forall$ F0 w u v, Arithmetical_ultraFilter F0
  -> w $\in$ ($\mathrm R$ F0) -> u $\in$ ($\mathrm R$ F0) -> v $\in$ ($\mathrm R$ F0)
  -> w <> R0 F0 -> w $\cdot$_F0 u = w $\cdot$_F0 v -> u = v.
Fact R_Mult_PrOrder : $\forall$ F0 u v w, Arithmetical_ultraFilter F0
  -> u $\in$ ($\mathrm R$ F0) -> v $\in$ ($\mathrm R$ F0) -> w $\in$ ($\mathrm R$ F0) -> (R0 F0) <_F0 w
  -> u <_F0 v <-> (w $\cdot$_F0 u) <_F0 (w $\cdot$_F0 v).
Fact R_inv : $\forall$ F0 u, Arithmetical_ultraFilter F0 -> u $\in$ ($\mathrm R$ F0)
  -> u <> R0 F0 -> $\exists$! v, v $\in$ ($\mathrm R$ F0) /\ v <> R0 F0 /\ u $\cdot$_F0 v = R1 F0.

Fact R_Mix_Association1 : $\forall$F0 u v w, Arithmetical_ultraFilter F0
  -> u $\in$ ($\mathrm R$ F0) -> v $\in$ ($\mathrm R$ F0) -> w $\in$ ($\mathrm R$ F0)
  -> (u +_F0 v) -_F0 w = u +_F0 (v -_F0 w).
Fact R_Mix_Association2 : $\forall$F0 u v w, Arithmetical_ultraFilter F0
  -> u $\in$ ($\mathrm R$ F0) -> v $\in$ ($\mathrm R$ F0) -> w $\in$ ($\mathrm R$ F0) -> w <> R0 F0
  -> u $\cdot$_F0 (v /_F0 w) = (u $\cdot$_F0 v) /_F0 w.
Fact R_Mult_Distribution_Minus : $\forall$F0 u v w, Arithmetical_ultraFilter F0
  -> u $\in$ ($\mathrm R$ F0) -> v $\in$ ($\mathrm R$ F0) -> w $\in$ ($\mathrm R$ F0)
  -> u $\cdot$_F0 (v -_F0 w) = (u $\cdot$_F0 v) -_F0 (u $\cdot$_F0 w).
\end{lstlisting}

\begin{definition}[Absolute Value on ${\mathbf R}$]
Let $x$ be a member in ${\mathbf R}$, then the absolute value of $x$ (denoted as $|x|$) is its image under the following class:
$$\{(u,v) : u \in {\mathbf R}\ \wedge\ ((u <_{\mathbf R} \mbox{\rm R}_0\ \wedge\ v = \mbox{\rm R}_0\ -_{\mathbf R}\ u)\
\vee\ (\mbox{\rm R}_0 <_{\mathbf R} u\ \wedge\ v = u)\ \vee\ (\mbox{\rm R}_0 = u = v))\},$$
which can be verified a function whose domain is ${\mathbf R}$ and range is contained in ${\mathbf R}$.
\end{definition}

\begin{lstlisting}
Definition R_Abs F0 := \{\ $\lambda$ u v, u $\in$ ($\mathrm R$ F0)
  /\ ((u <_F0 (R0 F0) /\ v = (R0 F0) -_F0 u)
    \/ (u = R0 F0 /\ v = R0 F0) \/ ((R0 F0) <_F0 u /\ v = u)) \}\.
Notation "| u |_ F0" := (R_Abs F0)[u]
  (at level 5, u at level 0, F0 at level 0) : r_scope.

Corollary RAbs_Corollary : $\forall$ F0, Arithmetical_ultraFilter F0
  -> Function (R_Abs F0) /\ dom(R_Abs F0) = $\mathrm R$ F0 /\ ran(R_Abs F0) $\subset$ ($\mathrm R$ F0).
Corollary RAbs_in_R : $\forall$ F0 u, Arithmetical_ultraFilter F0
  -> u $\in$ ($\mathrm R$ F0) -> |u|_F0 $\in$ ($\mathrm R$ F0).
\end{lstlisting}

The properties of this absolute value are same as those of $^{*}\mathbf Q$.

\begin{lstlisting}
Fact mt_R0_RAbs : $\forall$ F0 u, Arithmetical_ultraFilter F0
  -> u $\in$ ($\mathrm R$ F0) -> (R0 F0) <_F0 u -> |u|_F0 = u.
Fact eq_R0_RAbs : $\forall$ F0 u, Arithmetical_ultraFilter F0
  -> u $\in$ ($\mathrm R$ F0) -> |u|_F0 = R0 F0 <-> u = R0 F0.
Fact lt_R0_RAbs : $\forall$ F0 u, Arithmetical_ultraFilter F0
  -> u $\in$ ($\mathrm R$ F0) -> u <_F0 (R0 F0) -> |u|_F0 = (R0 F0) -_F0 u.
Fact R0_le_RAbs : $\forall$ F0 u, Arithmetical_ultraFilter F0
  -> u $\in$ ($\mathrm R$ F0) -> |u|_F0 $\in$ ($\mathrm R$ F0)
    /\ (|u|_F0 = R0 F0 \/ (R0 F0) <_F0 |u|_F0).
Fact Self_le_RAbs : $\forall$ F0 u, Arithmetical_ultraFilter F0
  -> u $\in$ ($\mathrm R$ F0) -> u = |u|_F0 \/ u <_F0 |u|_F0.
Fact RAbs_PrMult : $\forall$ F0 u v, Arithmetical_ultraFilter F0
  -> u $\in$ ($\mathrm R$ F0) -> v $\in$ ($\mathrm R$ F0)
  -> |u $\cdot$_F0 v|_F0 = |u|_F0 $\cdot$_F0 |v|_F0.
Fact RAbs_inEqu : $\forall$ F0 u v, Arithmetical_ultraFilter F0
  -> u $\in$ ($\mathrm R$ F0) -> v $\in$ ($\mathrm R$ F0)
  -> |u +_F0 v|_F0 = (|u|_F0 +_F0 |v|_F0)
    \/ |u +_F0 v|_F0 <_F0 (|u|_F0 +_F0 |v|_F0).
Fact neg_RAbs_equ : $\forall$ F0 u v, Arithmetical_ultraFilter F0
  -> u $\in$ ($\mathrm R$ F0) -> v $\in$ ($\mathrm R$ F0) -> u +_F0 v = R0 F0 -> |u|_F0 = |v|_F0.
\end{lstlisting}

\subsubsection{Rational Numbers}\label{sec5-4-4}

In \ref{sec5-3-5}, the structures of rational numbers ${^{*}\mathbf Q}_{\mathbf N}$, ${^{*}\mathbf Q}_{\mathbf Z}$ and ${^{*}\mathbf Q}_{\mathbf Q}$
are constructed, with which natural number set $\mathbf N$, integer set $\mathbf Z$ and rational number set $\mathbf Q$ can be conveniently established.

The process from ${^{*}\mathbf Q}_<$ to $\mathbf R$ is to classify the elements of ${^{*}\mathbf Q}_<$
through the equivalence relation that the difference between two elements is an infinitesimal.
Thus it is not hard to come out that $\mathbf N$ consists of equivalence classes represented by members of ${^{*}\mathbf Q}_{\mathbf N}$,
and the same idea applies to $\mathbf Z$ and $\mathbf Q$.

\begin{align}
{\mathbf N} &= \{u : \exists x \in {^{*}\mathbf Q}_{\mathbf N}, u = [x]_{R_{^{*}\mathbf Q}}\},\notag\\
{\mathbf Z} &= \{u : \exists x \in {^{*}\mathbf Q}_{\mathbf Z}, u = [x]_{R_{^{*}\mathbf Q}}\},\notag\\
{\mathbf Q} &= \{u : \exists x \in {^{*}\mathbf Q}_{\mathbf Q}, u = [x]_{R_{^{*}\mathbf Q}}\}. \notag
\end{align}

\begin{lstlisting}
Definition $\mathrm N$ F0 := \{ $\lambda$ u, $\exists$ n, n $\in$ (Q'_N F0) /\ u = \[n\]_F0 \}.
Corollary N_is_Set : $\forall$ F0, Arithmetical_ultraFilter F0 -> Ensemble ($\mathrm N$ F0).
Corollary N_properSubset_R : $\forall$ F0, Arithmetical_ultraFilter F0
  -> ($\mathrm N$ F0) $\subset$ ($\mathrm R$ F0) /\ ($\mathrm N$ F0) <> ($\mathrm R$ F0).

Definition $\mathrm Z$ F0 := \{ $\lambda$ u, $\exists$ z, z $\in$ (Q'_Z F0) /\ u = \[z\]_F0 \}.
Corollary Z_is_Set : $\forall$ F0, Arithmetical_ultraFilter F0 -> Ensemble ($\mathrm Z$ F0).
Corollary Z_properSubset_R : $\forall$ F0, Arithmetical_ultraFilter F0
  -> ($\mathrm Z$ F0) $\subset$ ($\mathrm R$ F0) /\ ($\mathrm Z$ F0) <> ($\mathrm R$ F0).

Definition $\mathrm Q$ F0 := \{ $\lambda$ u, $\exists$ q, q $\in$ (Q'_Q F0) /\ u = \[q\]_F0 \}.
Corollary Q_is_Set : $\forall$ F0, Arithmetical_ultraFilter F0 -> Ensemble ($\mathrm Q$ F0).
Corollary Q_Subset_R : $\forall$ F0, Arithmetical_ultraFilter F0 -> ($\mathrm Q$ F0) $\subset$ ($\mathrm R$ F0).
\end{lstlisting}

The relationships among the sets $\mathbf N$, $\mathbf Z$ and $\mathbf Q$
are definitely in line with the relationships among natural numbers, integers and rational numbers.

1) Natural numbers are non-negative integers:
$${\mathbf N} = \{u : u \in {\mathbf Z} \wedge (\mbox{R}_0 = u \vee \mbox{R}_0 <_{\mathbf R} u)\}.$$

2) Negative integers are negative values of all non-zero natural numbers:
$$\{u : u \in {\mathbf Z}\ \wedge\ u <_{\mathbf R} \mbox{R}_0\}
= \{u : u \in {\mathbf R}\ \wedge\ (\exists v \in ({\mathbf N} \sim \{\mbox{R}_0\}),\ u = \mbox{R}_0\ -_{\mathbf R}\ v )\}.$$

3) Dividing integers yields rational numbers:
$${\mathbf Q} = \{u : u \in {\mathbf R} \wedge
(\exists\ a\ b,\ a \in {\mathbf Z} \wedge b \in ({\mathbf Z} \sim \{\mbox{R}_0\}) \wedge u = a\ /_{\mathbf R}\ b)\},$$
namely,
$$\forall u,v,\ (u,v) \in {\mathbf Z} \times ({\mathbf Z} \sim \{\mbox{R}_0\})\ \Longrightarrow\ u\ /_{\mathbf R}\ v\ \in\ {\mathbf Q}.$$

4) Natural number set is a proper-subset of integer set; and integer set is a proper-subset of rational number set.

\begin{lstlisting}
Fact N_Equ_Z_me_R0 : $\forall$ F0, Arithmetical_ultraFilter F0
  -> $\mathrm N$ F0 = \{ $\lambda$ u, u $\in$ ($\mathrm Z$ F0) /\ (R0 F0 = u \/ (R0 F0) <_F0 u) \}.
Fact N_neg_Equ_Z_lt_R0 : $\forall$ F0, Arithmetical_ultraFilter F0
  -> \{ $\lambda$ u, u $\in$ ($\mathrm Z$ F0) /\ u <_F0 (R0 F0) \}
    = \{ $\lambda$ u, u $\in$ ($\mathrm R$ F0) /\ $\exists$ v, v $\in$ (($\mathrm N$ F0 ) $\sim$ [R0 F0]) /\ u +_F0 v = R0 F0 \}.
Fact Q_Equ_Z_Div : $\forall$ F0, Arithmetical_ultraFilter F0
  -> $\mathrm Q$ F0 = \{ $\lambda$ u, u $\in$ ($\mathrm R$ F0) /\ $\exists$ a b, a $\in$ ($\mathrm Z$ F0)
    /\ b $\in$ (($\mathrm Z$ F0) $\sim$ ([R0 F0])) /\ u = a /_F0 b \}.
Fact Z_Divide_in_Q : $\forall$F0 v u, Arithmetical_ultraFilter F0
  -> v $\in$ ($\mathrm Z$ F0) -> u $\in$ ($\mathrm Z$ F0) -> u <> R0 F0 -> (v /_F0 u) $\in$ ($\mathrm Q$ F0).
Fact N_properSubset_Z : $\forall$ F0, Arithmetical_ultraFilter F0
  -> ($\mathrm N$ F0) $\subset$ ($\mathrm Z$ F0) /\ ($\mathrm N$ F0) <> ($\mathrm Z$ F0).
Fact Z_properSubset_Q : $\forall$ F0, Arithmetical_ultraFilter F0
  -> ($\mathrm Z$ F0) $\subset$ ($\mathrm Q$ F0) /\ ($\mathrm Z$ F0) <> ($\mathrm Q$ F0).
\end{lstlisting}

Besides, the three sets each includes zero and identity element respectively, and operations are closed in them.

\begin{lstlisting}
Fact R0_in_N : $\forall$ F0, Arithmetical_ultraFilter F0 -> (R0 F0) $\in$ ($\mathrm N$ F0).
Fact R1_in_N : $\forall$ F0, Arithmetical_ultraFilter F0 -> (R1 F0) $\in$ ($\mathrm N$ F0).
Fact R0_in_Z : $\forall$ F0, Arithmetical_ultraFilter F0 -> (R0 F0) $\in$ ($\mathrm Z$ F0).
Fact R1_in_Z : $\forall$ F0, Arithmetical_ultraFilter F0 -> (R1 F0) $\in$ ($\mathrm Z$ F0).
Fact R0_in_Q : $\forall$ F0, Arithmetical_ultraFilter F0 -> (R0 F0) $\in$ ($\mathrm Q$ F0).
Fact R1_in_Q : $\forall$ F0, Arithmetical_ultraFilter F0 -> (R1 F0) $\in$ ($\mathrm Q$ F0).

Fact N_Plus_in_N : $\forall$ F0 u v, Arithmetical_ultraFilter F0
  -> u $\in$ ($\mathrm N$ F0) -> v $\in$ ($\mathrm N$ F0) -> (u +_F0 v) $\in$ ($\mathrm N$ F0).
Fact N_Mult_in_N : $\forall$ F0 u v, Arithmetical_ultraFilter F0
  -> u $\in$ ($\mathrm N$ F0) -> v $\in$ ($\mathrm N$ F0) -> (u $\cdot$_F0 v) $\in$ ($\mathrm N$ F0).

Fact Z_Plus_in_Z : $\forall$ F0 u v, Arithmetical_ultraFilter F0
  -> u $\in$ ($\mathrm Z$ F0) -> v $\in$ ($\mathrm Z$ F0) -> (u +_F0 v) $\in$ ($\mathrm Z$ F0).
Fact Z_Mult_in_Z : $\forall$ F0 u v, Arithmetical_ultraFilter F0
  -> u $\in$ ($\mathrm Z$ F0) -> v $\in$ ($\mathrm Z$ F0) -> (u $\cdot$_F0 v) $\in$ ($\mathrm Z$ F0).
Fact Z_Minus_in_Z : $\forall$F0 u v, Arithmetical_ultraFilter F0
  -> u $\in$ ($\mathrm Z$ F0) -> v $\in$ ($\mathrm Z$ F0) -> (u -_F0 v) $\in$ ($\mathrm Z$ F0).

Fact Q_Plus_in_Q : $\forall$ F0 u v, Arithmetical_ultraFilter F0
  -> u $\in$ ($\mathrm Q$ F0) -> v $\in$ ($\mathrm Q$ F0) -> (u +_F0 v) $\in$ ($\mathrm Q$ F0).
Fact Q_Mult_in_Q : $\forall$ F0 u v, Arithmetical_ultraFilter F0
  -> u $\in$ ($\mathrm Q$ F0) -> v $\in$ ($\mathrm Q$ F0) -> (u $\cdot$_F0 v) $\in$ ($\mathrm Q$ F0).
Fact Q_Minus_in_Q : $\forall$F0 v u, Arithmetical_ultraFilter F0
  -> v $\in$ ($\mathrm Q$ F0) -> u $\in$ ($\mathrm Q$ F0) -> (v -_F0 u) $\in$ ($\mathrm Q$ F0).
Fact Q_Divide_in_Q : $\forall$ F0 u v, Arithmetical_ultraFilter F0
  -> u $\in$ ($\mathrm Q$ F0) -> v $\in$ ($\mathrm Q$ F0) -> v <> R0 F0 -> (u /_F0 v) $\in$ ($\mathrm Q$ F0).
\end{lstlisting}

\subsubsection{Archimedean Property and Completeness}\label{sec5-4-5}

It seems that $\mathbf R$ possesses the same algebraic properties as $^{*}\mathbf Q$,
but the next contents guarantee that $\mathbf R$ is the expected real number set rather than $^{*}\mathbf Q$.

First of all, obtained by classifying Archimedes Subset $^{*}{\mathbf Q}_<$, $\mathbf R$ undoubtedly satisfies Archimedean Property:
$$\forall r \in {\mathbf R},\ \exists k \in {\mathbf N}, |r| = k \vee |r| <_{\mathbf R} k.$$

\begin{lstlisting}
Theorem Archimedes : $\forall$ F0 r, Arithmetical_ultraFilter F0 -> r $\in$ ($\mathrm R$ F0)
  -> $\exists$ n, n $\in$ ($\mathrm N$ F0) /\ (|r|_F0 = n \/ |r|_F0 <_F0 n).
\end{lstlisting}

Then with Archimedean Property, it can be proven that rational numbers are dense in real numbers, namely,
$$\forall u,v \in {\mathbf R},\ u <_{\mathbf R} v \Longrightarrow \exists q \in {\mathbf Q}, u <_{\mathbf R} q <_{\mathbf R} v.$$

\begin{lstlisting}
Theorem Q_Density : $\forall$ F0 u v, Arithmetical_ultraFilter F0
  -> u $\in$ ($\mathrm R$ F0) -> v $\in$ ($\mathrm R$ F0) -> u <_F0 v
  -> $\exists$ q, q $\in$ ($\mathrm Q$ F0) /\ u <_F0 q /\ q <_F0 v.
\end{lstlisting}

However, it is not hard to think of that $^{*}\mathbf Q$, including both infinities and infinitesimals,
does not possess the two properties above.

Moreover, completeness is what ultimately ensures that $\mathbf R$ is the general real number set.
As is known to us all, there are various methods to describe the completeness of real numbers, and the sequence method is chosen here.

Intuitively, a sequence of ${\mathbf R}$ is exactly a function whose domain is the natural number set and range is contained in ${\mathbf R}$.
Up to now we have seen several sets that can serve as the natural number set: $\omega$, ${^{*}\mathbf N}_{\mathbf N}$, ${^{*}\mathbf Q}_{\mathbf N}$ and $\mathbf N$.
Considering its concise set structure, we choose $\omega$ as the foundation of definitions about sequence.

\begin{definition}[Sequence and Subsequence of ${\mathbf R}$]
$\forall f \in {\mathbf R}^{\omega}$ is a sequence of ${\mathbf R}$. Meanwhile, if there exists $g \in {\omega}^{\omega}$ satisfying:
$$\forall m,n \in \omega,\ m < n\ \footnote{Namely m $\in$ n.} \Longrightarrow g(m) <_{\mathbf R} g(n),$$
then the following set
$$h = \{(u,v) : u \in \omega \wedge v = f(g(u))\}$$
is a sequence of $\mathbf R$, named subsequence of $f$.
\end{definition}

\begin{lstlisting}
Definition $\omega$_Seq f := Function f /\ dom(f) = $\omega$ /\ ran(f) $\subset$ $\omega$.
Definition R_Seq F0 f := Function f /\ dom(f) = $\omega$ /\ ran(f) $\subset$ ($\mathrm R$ F0).
Definition R_subSeq F0 h f := R_Seq F0 f /\ $\exists$ g, $\omega$_Seq g
  /\ ($\forall$ m n, m $\in$ $\omega$ -> n $\in$ $\omega$ -> m $\in$ n -> g[m] $\in$ g[n])
  /\ h = \{\ $\lambda$ u v, u $\in$ $\omega$ /\ v = f[g[u]] \}\.
Corollary R_subSeq_Corollary : $\forall$ F0 h f, R_subSeq F0 h f -> R_Seq F0 h.
\end{lstlisting}

\begin{definition}[Upper Bound and Least Upper Bound]
If sequence $f$ of $\mathbf R$ and real number $r \in \mathbf R$ satisfy:
$$\forall n \in \omega,\ f(n) <_{\mathbf R} r \vee f(n) = r,$$
then $r$ is called an upper bound of $f$; and if $r$ further satisfies:
$$\forall r_1 \in {\mathbf R},\ r_1 <_{\mathbf R} r \Longrightarrow \exists n \in \omega, r_1 <_{\mathbf R} f(n),$$
$r$ is called the least upper bound or supremum of $f$.
\end{definition}

\begin{lstlisting}
Definition R_UP F0 f r := r $\in$ ($\mathrm R$ F0) /\ ($\forall$ n, n $\in$ $\omega$ -> f[n] = r \/ f[n] <_F0 r).
Definition R_miniUP F0 f r := R_UP F0 f r
  /\ ($\forall$ r1, r1 $\in$ ($\mathrm R$ F0) -> r1 <_F0 r -> $\exists$ n, n $\in$ $\omega$ /\ r1 <_F0 f[n]).
\end{lstlisting}

\begin{theorem}[Completeness]
A monotonically increasing and bounded sequence of $\mathbf R$ must have a supremum. Namely, if sequence $f$ has an upper bound and satisfies:
$$\forall m,n \in \omega,\ f(m) <_{\mathbf R} f(n)\ \vee\ f(m) = f(n),$$
then there must exist the supremum of $f$.
\end{theorem}

\begin{lstlisting}
Definition R_monoIncrease F0 f := $\forall$ m n, m $\in$ $\omega$ -> n $\in$ $\omega$
  -> m $\in$ n -> f[m] = f[n] \/ f[m] <_F0 f[n].

Theorem R_Completeness : $\forall$ F0 f r, Arithmetical_ultraFilter F0
  -> ($\forall$ m, F0 <> F m) -> R_Seq F0 f -> R_monoIncrease F0 f
  -> R_UP F0 f r -> $\exists$ r1, R_miniUP F0 f r1.
\end{lstlisting}

\subsubsection{Irrational Number $\sqrt{2}$}\label{sec5-4-6}

At last, as the ending of FMCR, we introduce the formalization for a classic instance -- the existence of the irrational number $\sqrt{2}$.
Firstly, the concept of square root ``$\sqrt{\ }$" should be defined.

\begin{proposition}
For any positive real number $r$ (i.e., $r \in {\mathbf R} \wedge $R$_0 <_{\mathbf R} r$),
there exists a unique positive real number $r_1$ that satisfies: $r_1 \cdot_{\mathbf R} r_1 = r$.
\end{proposition}

\begin{lstlisting}
Theorem Square_Root_T1 : $\forall$ F0 r, Arithmetical_ultraFilter F0
  -> ($\forall$ n, F0 <> F n) -> r $\in$ ($\mathrm R$ F0) -> (R0 F0) <_F0 r
  -> $\exists$ r1, r1 $\in$ ($\mathrm R$ F0) /\ (R0 F0) <_F0 r1 /\ r = r1 $\cdot$_F0 r1.
Theorem Square_Root_T2 : $\forall$ F0 r, Arithmetical_ultraFilter F0
  -> ($\forall$ n, F0 <> F n) -> r $\in$ ($\mathrm R$ F0) -> (R0 F0) <_F0 r
  -> $\exists$! r1, r1 $\in$ ($\mathrm R$ F0) /\ (R0 F0) <_F0 r1 /\ r = r1 $\cdot$_F0 r1.
\end{lstlisting}

With the proposition above, for every positive real number $r \in {\mathbf R}$, its squire root can be defined as:
$$\sqrt{r} = \bigcap\{u : u \in {\mathbf R} \wedge {\mbox{R}}_0 <_{\mathbf R} u \wedge u \cdot_{\mathbf R} u = r\},$$
thus $\sqrt{r} \cdot_{\mathbf R} \sqrt{r} = r$ naturally holds.

\begin{lstlisting}
Definition Square_Root F0 r := $\cap$(\{ $\lambda$ u, u $\in$ ($\mathrm R$ F0)
  /\ (R0 F0) <_F0 u /\  u $\cdot$_F0 u = r \}).
Notation "$\sqrt{}$_ F0 r" := (Square_Root F0 r)(at level 5, F0 at level 4) : r_scope.
Corollary Square_Root_Corollary : $\forall$ F0 r, Arithmetical_ultraFilter F0
  -> ($\forall$ n, F0 <> F n) -> r $\in$ ($\mathrm R$ F0) -> (R0 F0) <_F0 r
  -> ($\sqrt{}$_F0 r) $\in$ ($\mathrm R$ F0) /\ (R0 F0) <_F0 ($\sqrt{}$_F0 r) /\ ($\sqrt{}$_F0 r) $\cdot$_F0 ($\sqrt{}$_F0 r) = r.
\end{lstlisting}

``$2$" in $\mathbf R$ (denoted as R$_2$) can be simply defined as R$_1\ +_{\mathbf R}$\ R$_1$,
and it can be proven that $\sqrt{2}$ is an irrational number, namely $\sqrt{\mbox{R}_2} \in ({\mathbf R} \sim {\mathbf Q})$.

\begin{lstlisting}
Definition R2 F0 := (R1 F0) +_F0 (R1 F0).
Corollary R2_Corollary : $\forall$ F0, Arithmetical_ultraFilter F0
  -> (R2 F0) $\in$ ($\mathrm N$ F0) /\ (R0 F0) <_F0 (R2 F0).

Fact Existence_of_irRational_Numbers : $\forall$ F0, Arithmetical_ultraFilter F0
  -> ($\forall$ n, F0 <> F n) -> ($\sqrt{}$_F0 (R2 F0)) $\in$ (($\mathrm R$ F0) $\sim$ ($\mathrm Q$ F0)).
\end{lstlisting}

The formal verification concerning $\sqrt{2}$ actually takes thousands of lines of code
and makes full use of the specific structures of $^{*}\mathbf Q$ and $\mathbf R$.

Up to now, we have completed the formalization of Filter-Method Construction for Reals.

\section{Conclusion}\label{sec6}

\subsection{What Are Real Numbers}\label{sec6-1}

The process of FMCR is accomplished in several steps: extending the set $\omega$ to the non-standard natural number set $^{*}\mathbf N$ through
an non-principal arithmetical ultrafilter; extending $^{*}\mathbf N$ to the non-standard integer set $^{*}\mathbf Z$ and
the non-standard rational number set $^{*}\mathbf Q$ through the idea of equivalence classification.
Figure \ref{fig4} illustrates the entire process from $\omega$ to $^{*}\mathbf Q$.

\begin{figure}[h]
\centering
\includegraphics[width=0.8\textwidth]{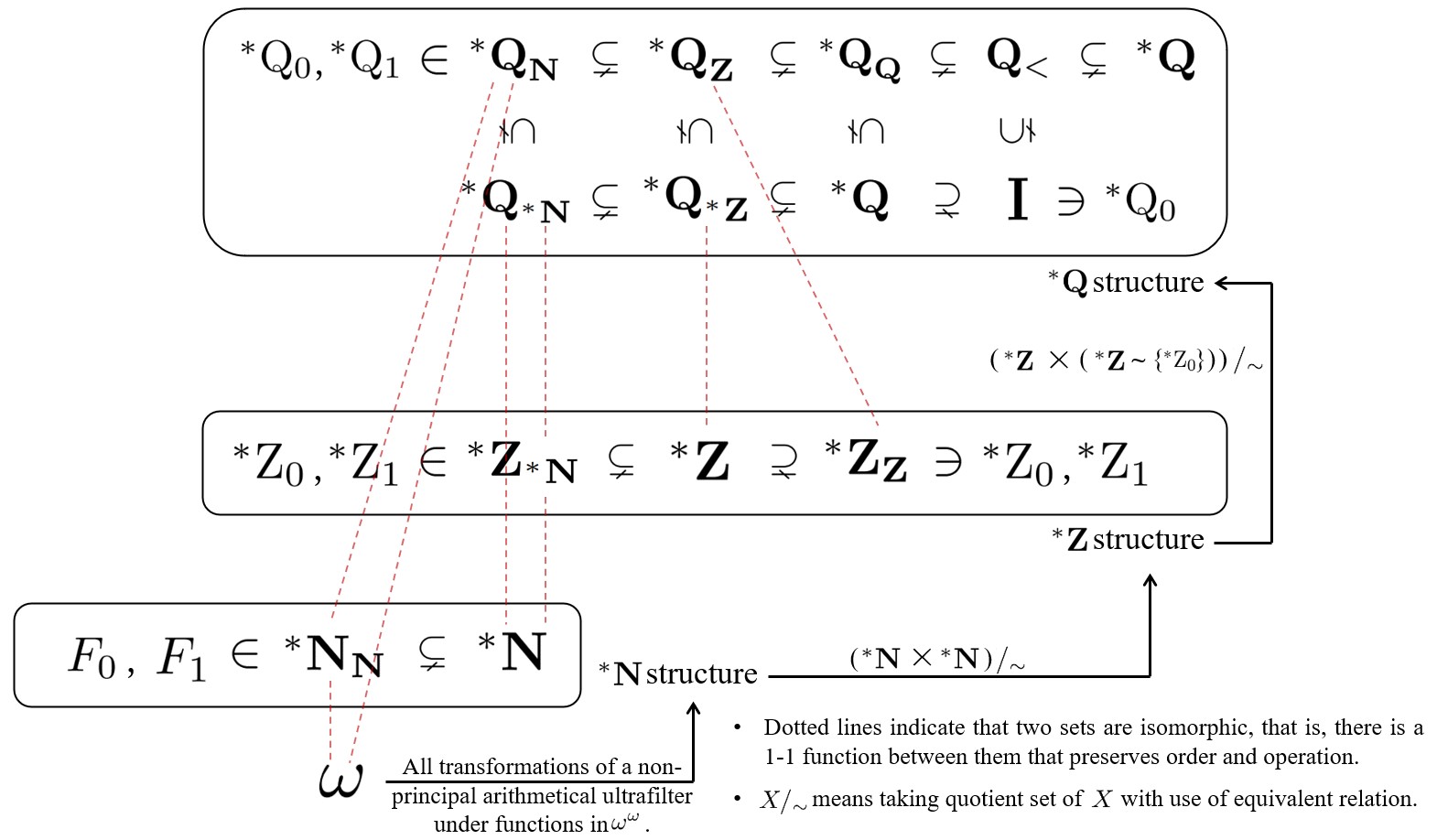}
\caption{the process from $\omega$ to $^{*}\mathbf Q$}\label{fig4}
\end{figure}

$^{*}\mathbf Q$ possesses the following algebraic properties:\\

(I) Properties of Addition

1) $u +_{^{*}\mathbf Q} {^{*}}$Q$_0 = u$.

2) existence of negative element: $\exists! u_0, u +_{^{*}\mathbf Q} u_0 = $Q$_0$.

3) associative law: $(u +_{^{*}\mathbf Q} v) +_{^{*}\mathbf Q} w = u +_{^{*}\mathbf Q} (v +_{^{*}\mathbf Q} w)$.

4) commutative law: $u +_{^{*}\mathbf Q} v = v +_{^{*}\mathbf Q} u$.\\

(II) Properties of Multiplication

1) $u \cdot_{^{*}\mathbf Q} {^{*}}$Q$_1 = u$.

2) existence of inverse element: $\exists! u_1, u \cdot_{^{*}\mathbf Q} u_1 = $Q$_1$.

3) associative law: $(u \cdot_{^{*}\mathbf Q} v) \cdot_{^{*}\mathbf Q} w = u \cdot_{^{*}\mathbf Q} (v \cdot_{^{*}\mathbf Q} w)$.

4) commutative law: $u \cdot_{^{*}\mathbf Q} v = v \cdot_{^{*}\mathbf Q} u$.\\

(I,II) Connection between Addition and Multiplication

distributive law: $u \cdot_{^{*}\mathbf Q} (v +_{^{*}\mathbf Q} w) = (u \cdot_{^{*}\mathbf Q} v) +_{^{*}\mathbf Q} (u \cdot_{^{*}\mathbf Q} w)$.\\

(III) Properties of Order

1) ${^{*}}$Q$_0 <_{^{*}\mathbf Q} {^{*}}$Q$_1$.

2) reflexivity: $\sim u <_{^{*}\mathbf Q} u$.

3) transitivity: $u <_{^{*}\mathbf Q} v \wedge v <_{^{*}\mathbf Q} w \Longrightarrow u <_{^{*}\mathbf Q} w$.

4) trichotomy: $u <_{^{*}\mathbf Q} v \vee v <_{^{*}\mathbf Q} u \vee u = v$.\\

(I,III) Connection between Addition and Order: $u <_{^{*}\mathbf Q} v \Longleftrightarrow (u +_{^{*}\mathbf Q} w) <_{^{*}\mathbf Q} (v +_{^{*}\mathbf Q} w)$.\\

(II,III) Connection between Multiplication and Order: $u <_{^{*}\mathbf Q} v \Longleftrightarrow (u \cdot_{^{*}\mathbf Q} w) <_{^{*}\mathbf Q} (v \cdot_{^{*}\mathbf Q} w)$,
where $w \neq {^{*}}$Q$_0$.\\

The structure satisfying the properties mentioned above is called ordered field,
thus $^{*}\mathbf Q$ and ${^{*}\mathbf Q}_{\mathbf Q}$ are two specific ordered field structure.
The real number set is required not only to fulfill the criteria of being an ordered field but also to satisfy:\\

(IV) Completeness: A monotonically increasing and bounded sequence must possess a supremum.\\

What mathematical analysis needs is exactly the complete ordered field.
Besides, there is one more crutial property of real numbers -- Archimedean Property:
$$\forall r \in {\mathbf R},\ \exists k \in {\mathbf N}, |r| <_{\mathbf R} k \vee |r| = k,$$
which means the measurability of quantity. $^{*}\mathbf Q$ does not satisfy this property;
although $^{*}{\mathbf Q}_<$, as a subset of $^{*}\mathbf Q$, satisfies Archimedean Property, it does not possess all properties of ordered field
(infinitesimals do not have their corresponding inverse elements in $^{*}{\mathbf Q}_<$).

Nevertheless, ${\mathbf R}$ still can be obtained by equivalently classifying $^{*}{\mathbf Q}_<$ through the equivalence relation $R_{^{*}\mathbf Q}$:
$$xR_{^{*}\mathbf Q}y \Longleftrightarrow (x -_{^{*}{\mathbf Q}} y \in {\mathbf I}).$$
The real number set $\mathbf R$ is another ordered field following $^{*}{\mathbf Q}$ and $^{*}{\mathbf Q}_{\mathbf Q}$, and is a complete ordered field.
As a subset of $\mathbf R$, $\mathbf Q$ is also an ordered field but lacks completeness, which is the essential difference between rational and real numbers.

Figure \ref{fig5} shows the process from $^{*}\mathbf Q$ to $\mathbf R$.

\begin{figure}[h]
\centering
\includegraphics[width=0.8\textwidth]{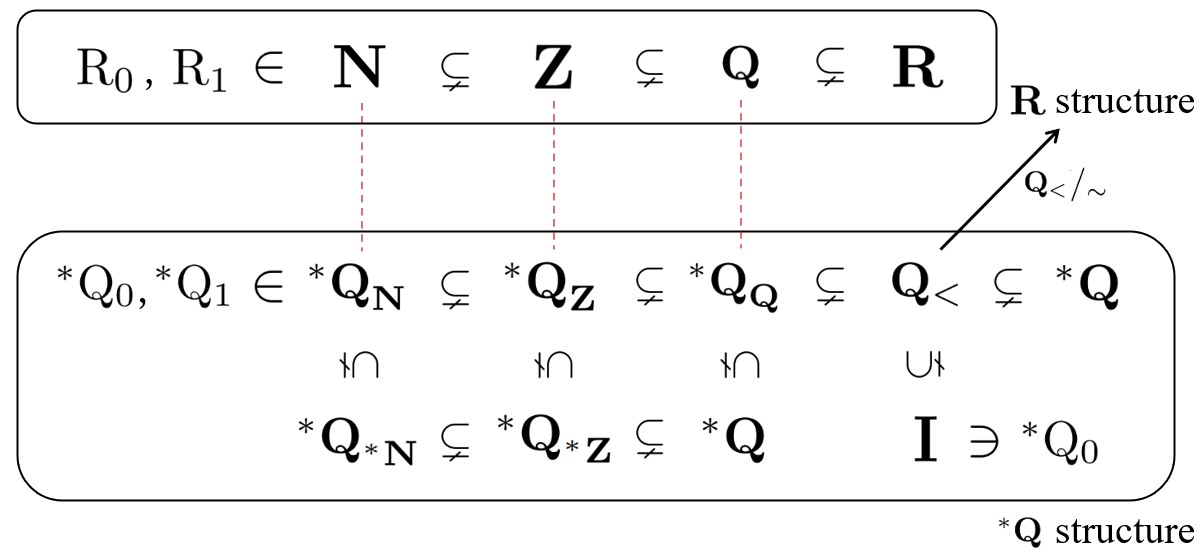}
\caption{the process from $^{*}\mathbf Q$ to $\mathbf R$}\label{fig5}
\end{figure}

Various theories on real numbers have been developed, such as Dedekind cut method\cite{Lan1966} by Dedekind, Cauchy sequence method by Cantor,
axiomatic method\cite{Zor2004} by Hilbert, and filter method introduced in this paper.
The complete ordered fields (i.e., real number fields) established by all different methods and satisfying Archimedean Property,
although structurally different, are isomorphic to each other.\cite{Zor2004}
And regardless of the method used to construct the real number set, the difficulty of introducing some form of actual infinity must be overcome.
Therefore, it is highly natural that infinity numbers enter mathematics as the identity of arithmetical ultrafilter
and the human comprehension of the infinite is boundless.

\subsection{Notes on the Implementation}\label{sec6-1}

As introduced in \ref{sec2-1}, substantial work has been done in formalizing real number theories.
However, the human quest for exploring the knowledge hidden within our world remains boundless.
We do not require our work to surpass the distinguished achievements of our predecessors
but hope to garner valuable insights from the field of formal verification and contribute our modest efforts to its advancement.

Previously, our team formally established the Dedekind cut real numbers within the framework of type theory.
This implementation is guided by Landau's {\it Foundations of Analysis}\cite{Lan1966}
and differs from the one built in the Coq standard library\cite{YFG2022}.
Since finishing the formalization of Morse-Kelley axiomatic set theory\cite{SY2020,YSF2020},
our focus has been dedicated to the formal verification for mathematical theorems grounded in set theory.

Concerning real number theory and based on set theory, we reconstructed Dedekind cut real number (still guided by Landau's {\it Foundations of Analysis}),
established an axiomatic real number system (guided by Zorich's {\it Mathematical Analysis}\cite{Zor2004}) and verified its compatibility and categoricalness in MK,
and the work presented in this paper is another advancement.

FMCR, as mentioned in introduction, is a distinctive real number theory.
As its precondition, the existence of non-principal arithmetical ultrafilter is quite not easy to prove within general set theories.
Wang introduces several proof methods in \cite{Wan2016}, respectively using different additional set hypothesis independent of ZFC (and MK).
Following this we have also verified the existence of non-principal arithmetical ultrafilter
with the assistance of the Continuum Hypothesis based on MK in Coq. This will be presented in another paper.

The formalization of FMCR encompasses 12 files (4 for filters and 8 for other components) and comprises nearly 20,000 lines of code.
It is self-contained and does not rely on any libraries in Coq.
This formalization includes the fundamental concepts and properties of filters,
as well as the non-standard number sets $^{*}\mathbf N$, $^{*}\mathbf Z$, $^{*}\mathbf Q$,
and the standard real number set $\mathbf R$.
The definitions, theorems, and propositions presented in this paper, which constitute the FMCR theory,
are primarily derived from \cite{CN1974,Dag1979,Wan1998,Wan2000,Wan2016,Wan2018}.
We present in this paper all their formal descriptions, all of which have been rigorously proven correct in Coq.

To advance further, it is possible to formalize additional algebraic properties, such as exponentiation and summation operations,
in preparation for the formalization of mathematical analysis.
Alternatively, one can enhance the implementation of filters to prepare for the formalization of non-standard analysis.
Both of these paths are viable and we are actively exploring them at present.

\bmhead{Acknowledgments}

Our work is funded by National Natural Science Foundation (NNSF) of China under Grant 61936008.

\end{document}